\documentclass[pdflatex,sn-mathphys-num]{sn-jnl}

\usepackage{amsmath,amssymb,amsfonts}%
\usepackage{amsthm}%
\usepackage{mathrsfs}%
\usepackage[title]{appendix}%
\usepackage{xcolor}%
\usepackage{textcomp}%
\usepackage{manyfoot}%
\usepackage{booktabs}%
\usepackage{algorithm}%
\usepackage{algorithmicx}%
\usepackage{algpseudocode}%
\usepackage{amsmath, amsfonts} 
\usepackage{listings}%
\usepackage{booktabs}

\newtheorem{lemma}{Lemma}

\usepackage{multirow}  
\usepackage{makecell}  
\usepackage{booktabs}  
\usepackage{graphicx}    
\usepackage{caption}     
\usepackage{float}       
\usepackage{geometry}    
\usepackage{tikz}        

\usepackage{xcolor}     

\usepackage{hyperref}  
\hypersetup{colorlinks=true, linkcolor=blue}  

\usepackage{subcaption}  

\usepackage{enumitem}  
\setlist[enumerate,1]{label=(\roman*), leftmargin=2em, topsep=0.5ex, itemsep=0.3ex}

\theoremstyle{thmstyleone}%
\newtheorem{theorem}{Theorem}
\theoremstyle{thmstyletwo}%
\newtheorem{example}{Example}%
\newtheorem{remark}{Remark}%

\theoremstyle{thmstylethree}%
\newtheorem{definition}{Definition}%

\begin{document}

\title[Article Title]{Large-Scale Tensor Trace-Ratio Problem: Randomized Algorithm and Sparse Multilinear Discriminant Analysis via Closed-Form Solution}


\author[1]{\fnm{Yuhang} \sur{Li}}\email{tb24080014a41@cumt.edu.cn}

\author*[1]{\fnm{Gang} \sur{Wu}}\email{gangwu@cumt.edu.cn}

\affil*[1]{\orgdiv{School of Mathematics}, \orgname{China University of Mining and Technology}, \orgaddress{\city{ Xuzhou}, \postcode{221116}, \state{Jiangsu}, \country{P.R. China}}}

\abstract{Trace-ratio optimization problems from high-dimensional tensor data appear across numerous research domains including machine learning, signal processing, as well as pattern recognition. Up to now, analytical closed-form expressions remain unavailable for tensor trace-ratio tasks; prevailing techniques mostly adopt inner-outer iterative schemes to tackle such optimization tasks. However, when dealing with high-dimensional tensor data sets, existing multilinear discriminant analysis methods often encounter limitations such as heavy computational costs, insufficient memory space, and interference from redundant features. This paper yields two core contributions.
First, we present a {\it closed-form} solution under the full row rank condition, and construct a randomized solver targeting large-scale tensor trace-ratio optimization. In contrast to prior approaches, the devised algorithm achieves substantial computational speedup. Its core mechanism is to reformulate the original large-scale tensor trace-ratio optimization into a tensor least-squares formulation, instead of the much more complicated (generalized) tensor eigenproblems. We derive theoretical error-bound results which verify the reliability of the newly developed algorithm. Second, we build a {\it sparse} multilinear discriminant analysis model leveraging the above closed-form result combined with tensor $\ell_{2,0}$-norm regularization, and propose a sparse multilinear discriminant method for tensor feature selection. The sparse method can identify key discriminative features from high-dimensional data and suppress redundant information. Comprehensive numerical evaluations are conducted over multiple publicly available real-world benchmarks.
The proposed methods demonstrate significant advantages in computational efficiency and scalability for large-scale problems, while achieving competitive classification accuracy on some datasets.
}


\maketitle

\vspace{-26pt}
\noindent \textbf{\fontsize{10}{12}\selectfont ~~~~~~Mathematics Subject Classification (2020)} \fontsize{10}{12}\selectfont 65F15 · 65F10

\section{Introduction}\label{sec1}
Driven by growing big-data resources, large-scale tensor-form datasets are increasingly encountered across diverse research areas including signal processing \cite{singp1,singp2}, pattern recognition \cite{pr}, as well as machine-learning investigations \cite{Webb2011Statistical}.
With the increase of data dimensionality, one has to face with the challenge of curse of dimensionality. Dimensionality reduction approaches serve as an effective remedy for this issue. Such techniques aim to project raw data onto a low-dimensional embeddings from its high-dimensional representations while retaining the intrinsic information embedded in raw observations to the greatest extent \cite{FukunagaLDA}.

Typical dimensionality reduction schemes reduce to the trace-ratio problem \cite{LDA,NgoLDA,Matrix_Fish}.
When performing dimensionality reduction on high-dimensional observations such as facial data, image data and text data, traditional methods first convert the data into (1-D) vectors and then process them using 1-D techniques such as LDA \cite{LDA} and PCA \cite{PCA}.

Nowadays, representing data and formulating optimization problems using higher-order tensors present us significant challenges for contemporary machine learning methods.
With the development of data sciences, more and more dimensionality reduction methods have been investigated for tensor data \cite{Ye, Liu, Lu,MLPA1,Yan}.
For instance, Lu {\it et al.} \cite{Lu} developed an approach that yields uncorrelated discriminative features via a tensor-to-vector projection. In \cite{Yan}, Yan {\it et al.} introduced a discriminant analysis approach built upon tensor representation (DATER) for tensor data classification via a tensor-to-tensor projection.
However, DATER fails to converge during iterative updates. This makes it difficult to specify a proper stopping iteration count, a factor that directly affects the overall performance of DATER.
Recently, some multilinear discriminant analysis methods based on different tensor products were proposed for dimensionality reduction \cite{Tensor_product1, Guide2026, Tracevs, Zahir}. The work reported in \cite{Tensor_product1} by Dufrenois {\it et al.} put forward a transform domain (TD)-based approach to solve the multilinear discriminant analysis (MDA) problem.
Zahir {\it et al.} \cite{Tracevs} generalized linear discriminant analysis to higher-order tensors built upon the tensor $m$-mode product framework.
Nevertheless, both methods adopt inner-outer iterative schemes to tackle tensor trace-ratio optimization. Repeated high-order tensor calculations in iterations bring prohibitive computational expenses, making them inefficient for large-scale data scenarios.

Furthermore, researchers have extended a variety of classic dimensionality reduction methods to the tensor trace-ratio model framework.
To name a few, El Guide {\it et al.} \cite{Guide2026} provided a unified t-product trace-optimization framework for multidimensional reduction, generalized mainstream matrix methods to third-order tensors with linear/kernel extensions, and unified manifold learning techniques under this paradigm for vectorization-free tensor computation.
Using the Einstein product, Zahir {\it et al.} \cite{Zahir} generalized some supervised methods, unsupervised methods and kernel methods to the tensor trace-ratio model. While both works enrich multidimensional tensor reduction techniques, they still suffer from the limitation of prohibitive computational costs when handling large-scale high-dimensional data.

Some existing methods rewrite the tensor {\it trace-ratio} problem as a tensor {\it ratio-trace} problem, and compute approximate solution via solving a large-scale generalized tensor eigenvalue problem \cite{Lu,MLPA1,Yan,Zahir}. Mathematically, the tensor trace-ratio formulation and tensor ratio-trace formulation are fundamentally inequivalent.
Thus, these methods risk generating solutions that deviate from the original target solutions. Recently, some {\it inner-outer} iterative methods were utilized to solve the tensor {trace-ratio} problem  \cite{Tensor_product1, Tracevs}. However, they may incur substantial computational overhead and memory demands, especially on high-dimensional large-sample tensor datasets. Thus, this motivates the development of new approaches to accelerate the numerical computation for large-scale tensor trace-ratio problems.

On the other hand, another practical challenge arises from the fast expansion of high-dimensional large-scale tensor data, where conventional multilinear discriminant analysis schemes suffer from the curse of dimensionality and adverse impacts caused by redundant information. How to extract key discriminative information during the dimensionality reduction process has become a critically desirable research objective. Applying sparse representation or sparse learning methods to feature selection can achieve this objective. As far as we are aware, current feature selection approaches built upon multilinear discriminant analysis work with a reformulated {\it ratio-trace} formulation in place of the canonical {\it trace-ratio} optimization. Such treatment brings about large-scale generalized tensor eigenvalue computations \cite{Lai_Sparse,Wang_Sparse,Sparse1,Sparse2}. As was mentioned before, this procedure can result in discarding valuable information, moreover, large-scale generalized tensor eigenvalue problems are difficult to solve in general.

The primary contributions of the present work can be outlined as below:
\begin{itemize}
\item{To start with, we derive a closed-form solution to the tensor trace-ratio problem under the full row-rank requirement. More precisely, $\mathrm{bcirc}(\mathcal{X}^T)$ possesses full row rank, where $\mathcal{X}$ stands for training samples.
    With the aid of the derived closed-form solution, the multilinear discriminant analysis method avoids employing the inner-outer iteration framework. Consequently, substantial reductions in computational complexity can be achieved, alongside remarkable improvements in computational efficiency for conventional algorithms.

    Instead of solving large-scale tensor eigenvalue problems, we cast the tensor trace-ratio problem into the form of a large-scale tensor least-squares formulation, and develop a randomized algorithm for its calculation. Theoretical analysis is provided to validate both rationality and practical performance of our proposed framework.}
\item{Considering the tensor trace-ratio problem together with tensor $\ell_{2,0}$-norm regularization, this work constructs a new sparse multilinear discriminant analysis model (SMDA) targeted at tensor feature selection. Moreover, we rewrite the SMDA model as a constrained tensor trace-optimization problem, and present an efficient way for it.
    Extensive numerical tests performed on multiple real-world datasets illustrate practicability and competitive advantages of our presented approach against several existing state-of-the-art algorithms for tensor dimensionality reduction and feature selection.
    }
\end{itemize}

This paper is structured as follows. Section \ref{sec2} briefly reviews essential tensor preliminaries and outlines the tensor trace-ratio formulation tailored for tensor dimensionality reduction. Section \ref{sec3} derives a closed-form solution to tensor trace-ratio problems, develops a randomized algorithm to calculate the projection tensor, and provides theoretical analysis supporting the presented algorithm. Section \ref{sec4} establishes a new sparse multilinear discriminant analysis model alongside its algorithm dedicated to tensor feature selection. Section \ref{sec5} reports our numerical evaluations, while concluding remarks are presented in Section \ref{sec6}.

\section{Preliminaries}\label{sec2}
This section introduces notations and definitions used in this paper.
\subsection{Notations}
A tensor is defined as a high-dimensional array extending the concepts of vectors and matrices. In this work, bold Euler-script letters $\mathcal{A}\in\mathbb{R}^{I\times J \times K}$ denote third-order tensors, with $\mathcal{A}(i,j,k)$ referring to its $(i,j,k)$-th entry. Indeed, tensor can also be characterized in the form of slices \cite{t-product2}. We adopt $\mathcal{A}(i,:,:)$, $\mathcal{A}(:,j,:)$, and $\mathcal{A}(:,:,k)$ to indicate horizontal, lateral and frontal slices correspondingly, while $\mathcal{A}^{(i)}$ represents the $i$-th frontal slice.
According to \cite{Kolda_ten}, a tube associated with a third-order tensor can be obtained by fixing the first two indices and varying the third index, written as $\mathcal{A}(i,j,:)$. Matrices correspond to bold uppercase symbols such as $\mathbf{A}$. Vectors adopt bold-form lowercase notations such as $\mathbf{a}$. Scalars are denoted by ordinary lowercase letters, for instance $a$.
Moreover, $\mathbf{0}$ stands for the zero vector, $\mathbf{O}$ stands for the zero matrix and $\mathcal{O}$ stands for the zero tensor, respectively. Table \ref{Ta1} provides a summary of these notations.
\begin{table}[h]
	\centering
	\caption{List of notations adopted throughout this work} \label{Ta1}
	\begin{tabular}{{@{}llll@{}}}
		\toprule
		\textbf{Notations} & \textbf{Description} \\
		\midrule
		$\mathcal{X}$ & Training samples \\
		$\overrightarrow{\mathcal{X}_j}$ & The $j$-th lateral slice belonging to tensor $\mathcal{X}$\\
        $\mathcal{X}^{(i)}$ & The $i$-th frontal slice belonging to tensor $\mathcal{X}$\\
		$c$ & Total number of classes \\
		$d$ & Target dimension after dimensionality reduction \\
		$N, m$ & Quantities for overall samples and training instances \\
		$\boldsymbol{0}, \mathbf{O},\mathcal{O}$ & Zero-valued vector, zero-valued matrix and zero-valued tensor \\
		$ \mathcal{I},\mathbf{I},\boldsymbol{1}_i$ & Identity tensor, identity matrix, $i$-dimensional all-ones vector \\
		$\text{rank}(\mathbf{A})$ & Rank associated with matrix $\mathbf{A}$ \\
		$\text{tr}(\mathcal{A})$ & Trace computed for tensor $\mathcal{A}$ \\
        $\sigma_i$ & $i$-th largest singular value yielded by matrix $\mathbf{X}$ \\
		$\mathcal{A}^T, \mathcal{A}^H, \mathcal{A}^\dagger$ & Transpose, conjugate transpose as well as Moore-Penrose inverse of tensor $\mathcal{A}$ \\
		$\text{span}\{\mathcal{A}\}$ & Subspace spanned by lateral slices from tensor $\mathcal{A}$ \\
		$\mathcal{N}(\mathcal{A})$ & Null subspace corresponding to tensor $\mathcal{A}$ \\
		$\|\cdot\|_2$ & 2-norm defined for matrices or tensors \\
        $\|\cdot\|_F$ & Frobenius norm for a given tensor \\
        $\|\cdot\|_{2,0}$, ${|||\cdot|||}_{2,0}$ &  $\ell_{2,0}$-norm for matrices or tensors \\
		$\mathcal{N}(\mathcal{A}) \setminus \mathcal{N}(\mathcal{B})$ & Subspace lying inside $\mathcal{N}(\mathcal{A})$ while excluded from $\mathcal{N}(\mathcal{B})$ \\
		$\mathcal{P}_\mathcal{A}$ & Orthogonal projector onto subspace $\text{span}\{\mathcal{A}\}$ \\
		$\sin\angle(\mathcal{V}, \tilde{\mathcal{V}})$ & Sine of the angle between the two subspaces $\mathcal{V}$ and $\tilde{\mathcal{V}}$\\
		\botrule
	\end{tabular}
\end{table}

\subsection{The t-product}
For tensor $\mathcal{A}\in\mathbb{R}^{n_1\times n_2\times n_3}$, to elaborate on t-product operation, we begin by presenting several block-matrix-related definitions \cite{t-product2,t-product1}.

$\bullet$ For tensor $\mathcal{A}$, its block circulant matrix ${\rm bcirc}\left(\mathcal{A}\right)$ \cite{t-product2} is
\begin{align*}
	{\rm bcirc}\left(\mathcal{A}\right)=\begin{bmatrix}
		\mathcal{A}^{(1)}&\mathcal{A}^{(n_3)}&\cdots&\mathcal{A}^{(2)}\\
		\mathcal{A}^{(2)}&\mathcal{A}^{(1)}&\cdots&\mathcal{A}^{(3)}\\
		\vdots&\vdots&\ddots&\vdots\\
		\mathcal{A}^{(n_3)}&\mathcal{A}^{(n_3-1)}&\cdots&\mathcal{A}^{(1)}
	\end{bmatrix} \in\mathbb{R}^{n_1n_3\times n_2n_3}.
\end{align*}
Making use of discrete Fourier matrix $\mathbf{F}_{n_3} \in \mathbb{C}^{n_3 \times n_3}$, one can achieve the block-diagonalization for ${\rm bcirc}\left(\mathcal{A}\right)$ using the discrete Fourier transform \cite{t-bdiag}
\begin{equation*}
	\left(\mathbf{F}_{n_3}^H\otimes \mathbf{I}_{n_1}\right){\rm bcirc}\left(\mathcal{A}\right)\left(\mathbf{F}_{n_3}\otimes \mathbf{I}_{n_2}\right)={\rm bdiag}\left(\widehat{\mathcal{A}}\right).
\end{equation*}

$\bullet$ The operator ${\rm unfold}\left(\mathcal{A}\right)$ \cite{t-product2} is used to convert the tensor $\mathcal{A}$ into a matrix with dimensions $n_1n_3\times n_2$:
\begin{align*}
	{\rm unfold}\left(\mathcal{A}\right)=\begin{bmatrix}
		\mathcal{A}^{(1)}\\
		\mathcal{A}^{(2)}\\
		\vdots\\
		\mathcal{A}^{(n_3)}
	\end{bmatrix} \in\mathbb{R}^{n_1n_3\times n_2},
\end{align*}
and the operator ${\rm fold}$ represents the inverse operator, i.e., ${\rm fold}\left({\rm unfold}\left(\mathcal{A}\right)\right)=\mathcal{A}$.

$\bullet$ Tensor block-diagonal matrix form \cite{t-product2} is
\begin{align*}
	{\rm bdiag}\left(\mathcal{A}\right)=\begin{bmatrix}
		\mathcal{A}^{(1)}&&&\\
		&\mathcal{A}^{(2)}&&\\
		&&\ddots &\\
		&&&\mathcal{A}^{(n_3)}
	\end{bmatrix} \in\mathbb{R}^{n_1n_3\times n_2n_3}.
\end{align*}

Under the {\sc Matlab} notations, $\widehat{\mathcal{A}}={\rm fft}\left(\mathcal{A},[~],3\right)$ denotes a tensor yielded by performing fast Fourier transform over each tube fiber of $\mathcal{A}$, while its inverse transformation reads $\mathcal{A}={\rm ifft}\left(\widehat{\mathcal{A}},[~],3\right)$.

\begin{definition}\label{def1}
	(t-product \cite{t-product2, t-pro_Algorithm}). The t-product of tensor $\mathcal{A}\in\mathbb{R}^{n_1\times n_2\times n_3}$ and $\mathcal{B}\in\mathbb{R}^{n_2\times l\times n_3}$ is defined as
	\begin{equation*}
		\mathcal{M}=\mathcal{A}\ast\mathcal{B}={\rm fold}\left({\rm bcirc}\left(\mathcal{A}\right){\rm unfold}\left(\mathcal{B}\right)\right)\in\mathbb{R}^{n_1\times l\times n_3}.
	\end{equation*}
\end{definition}
As reported in \cite{t-product1}, the t-product may be computed efficiently via the fast Fourier transform, namely,
\begin{equation*}
	\widehat{\mathcal{M}}^{(i)}=\widehat{\mathcal{A}}^{(i)}\widehat{\mathcal{B}}^{(i)},~i=1,2,\ldots, n_3.
\end{equation*}

\begin{definition}\label{def2}
	(Orthogonal tensor \cite{t-orth_tensor}) An orthogonal tensor $\mathcal{Q} \in\mathbb{R}^{n\times n\times l}$ under t-product is defined as
	\begin{equation*}
		\mathcal{Q}^T\ast\mathcal{Q}=\mathcal{Q}\ast\mathcal{Q}^T =\mathcal{I}.
	\end{equation*}
\end{definition}
\begin{definition}\label{def3}
	(Identity tensor \cite{t-product1}). An identity tensor $\mathcal{I}$ is a tensor whose first frontal slice is an identity matrix and the rest are zero matrices.
\end{definition}

\begin{definition}
   (Positive (semi-) definite tensor \cite{Bouallala}). Consider the tensor $\mathcal{A} \in \mathbb{R}^{n \times n \times n_3}$, the tensor $\mathcal{A}$ is positive (semi-) definite if and only if each frontal slice $\widehat{\mathcal{A}}^{(i)}$ is positive (semi-) definite.
\end{definition}

\begin{definition}\label{def4}
	(Tensor inverse \cite{t-product2}). An $n\times n \times n_3$ tensor $\mathcal{A}$ has an inverse $\mathcal{B}$, if
	\begin{equation*}
		\mathcal{A}\ast \mathcal{B}=\mathcal{I}, \quad and\quad \mathcal{B}\ast \mathcal{A}=\mathcal{I}.
	\end{equation*}
\end{definition}
\begin{definition}\label{def5}
	(Tensor Moore-Penrose inverse \cite{tgen_inv1,tgen_inv2}). A tensor $\mathcal{A}^\dagger \in\mathbb{R}^{n_2 \times n_1 \times n_3}$ is called the Moore-Penrose inverse of $\mathcal{A}\in\mathbb{R}^{n_1 \times n_2 \times n_3}$, if $\mathcal{A}^\dagger$ satisfies the following four equations
	\begin{equation*}
		\mathcal{A}\ast \mathcal{A}^\dagger\ast \mathcal{A}=\mathcal{A},~\mathcal{A}^\dagger\ast \mathcal{A}\ast\mathcal{A}^\dagger=\mathcal{A}^\dagger,~{\left(\mathcal{A}\ast \mathcal{A}^\dagger\right)}^H=\mathcal{A}\ast \mathcal{A}^\dagger,~{\left(\mathcal{A}^\dagger \ast\mathcal{A}\right)}^H=\mathcal{A}^\dagger \ast\mathcal{A}.
	\end{equation*}
\end{definition}
\begin{remark}
	{\it In particular, if ${\rm bcirc}(\mathcal{A})$ has full-row rank, one obtains the pseudoinverse
	\begin{equation*}
		\mathcal{A}^\dagger={\rm fold } \left({{\rm bcirc}(\mathcal{A})}^T{\left[{\rm bcirc}\left(\mathcal{A}\right){{\rm bcirc}\left(\mathcal{A}\right)}^T\right]}^{-1}{\rm unfold}\left(\mathcal{I}\right)\right),
	\end{equation*}
	and $\mathcal{A}\ast \mathcal{A}^\dagger=\mathcal{I}$.}
\end{remark}
\begin{definition}\label{def6}
	\cite{t-product2} Let the tensor $\mathcal{A}\in\mathbb{C}^{n_1\times n_2 \times n_3}$, then the conjugate transpose of $\mathcal{A}^H\in\mathbb{C}^{n_2\times n_1 \times n_3}$, ${\left(\mathcal{A}^H\right)}^{(1)}={\left({\mathcal{A}}^{(1)}\right)}^H$ and ${\left(\mathcal{A}^H\right)}^{(i)}={\left({\mathcal{A}}^{(n_3-i+2)}\right)}^H$, $i=2,3,\ldots,n_3$. Specifically, there holds
	${\left(\widehat{\mathcal{A}^H}\right)}^{(i)}={\left(\widehat{\mathcal{A}}^{(i)}\right)}^H$ in the Fourier-transformed domain.
\end{definition}
\begin{definition}\label{def7}
	(Tensor trace \cite{Bouallala}). The trace of a tensor $\mathcal{A}\in \mathbb{R}^{n\times n\times n_3}$ is defined as
	\begin{equation*}
		{\rm tr}\left(\mathcal{A}\right)=\frac{1}{n_3}\sum_{i=1}^{n_3}{\rm tr}\left({\widehat{\mathcal{A}}}^{(i)}\right).
	\end{equation*}
\end{definition}
Moreover, it follows from \cite{t-svd} that
\begin{equation*}
	{\mathcal{A}}^{(1)}=\frac{1}{n_3}\sum_{i=1}^{n_3}{\widehat{\mathcal{A}}}^{(i)}.
\end{equation*}
Accordingly, one may also define the trace for a third-order tensor as
\begin{equation*}
	{\rm tr}\left(\mathcal{A}\right)={\rm tr}\left(\mathcal{A}^{(1)}\right).
\end{equation*}

\begin{theorem}\label{theo1}
	(t-SVD \cite{t-svd}) Let $\mathcal{A}\in \mathbb{R}^{n_1\times n_2\times n_3}$. Then it can be factorized as
	\begin{equation*}
		\mathcal{A}=\mathcal{U}\ast\mathcal{S}\ast\mathcal{V}^T,
	\end{equation*}
	where $\mathcal{U}\in \mathbb{R}^{n_1\times n_1\times n_3}$, $\mathcal{V}\in \mathbb{R}^{n_2\times n_2\times n_3}$ are orthogonal and $\mathcal{S}\in \mathbb{R}^{n_1\times n_2\times n_3}$ is called the singular tensor of $\mathcal{A} $.
\end{theorem}

\begin{definition}\label{def8}
	(Tensor tubal rank \cite{t-tubal_rank}). Let $\mathcal{A}\in\mathbb{R}^{n_1\times n_2 \times n_3}$ be a third-order tensor, its tubal rank is defined as
	\begin{equation*}
		{\rm rank_t}\left(\mathcal{A}\right)=\#\left\{i:\mathcal{S}\left(i,i,:\right)\ne0\right\},~i=1,2,\ldots,\min\left\{n_1,n_2\right\},
	\end{equation*}
	where $\#$ denotes the number of nonzero tubes in the singular tensor $\mathcal{S}$.
\end{definition}
\begin{definition}\label{def10}
	\cite{Spectral} The null space of a tensor $\mathcal{A} \in \mathbb{R}^{n_1\times n_2\times n_3}$ is defined as
\begin{equation*}
		\mathcal{N}(\mathcal{A}) = \left\{ \vec{\mathcal{X}}_1, \ldots, \vec{\mathcal{X}}_r \in \mathcal{R}^{n_2\times 1\times n_3}: \widehat{\vec{\mathcal{X}}}_j^{(i)} \in \mathcal{N} \left( \widehat{\mathcal{A}}^{(i)} \right) \text{ for } i = 1, 2,\ldots,n_3 \right\}.
	\end{equation*}
	It is seen that if $\mathcal{X} \in \mathbb{R}^{n_2\times 1\times n_3}$ belongs to the null space space associated with tensor $\mathcal{A} \in \mathbb{R}^{n_1\times n_2\times n_3}$, then $\mathcal{A}\ast \mathcal{X}=\mathcal{O}$.
\end{definition}
\begin{definition} \label{tensubspace}
  (Tensor subspace \cite{t-product2, Tensubspace}) Given a tensor $\mathcal{A} \in \mathbb{R}^{m \times p \times n}$, where the lateral slices $\overrightarrow{\mathcal{A}_i} \in \mathbb{R}^{m \times 1 \times n},~i=1,2,\ldots,p$, are linearly independent. Then we call
\[
{\rm span}\left\{\mathcal{A}\right\} = \left\{ \mathcal{Y} \,\middle|\, \mathcal{Y} = \sum_{i=1}^{p}\overrightarrow{\mathcal{A}_i}\ast\mathbf{c}_i= \mathcal{A} * \mathcal{C},\quad \forall\, \mathcal{C} \in \mathbb{R}^{p \times 1 \times n} \right\}
\]
a tensor subspace of dimension $p$, with $\mathbf{c}_i=\mathcal{C}(i,1,:),~i=1,2,\ldots,p$, being the tubal-scalars.
\end{definition}
\begin{theorem}\label{theo2}
	(t-QR \cite{t-product2}). Let $\mathcal{A}\in\mathbb{R}^{n_1\times n_2 \times n_3}$. Then it can be factorized as
	\begin{equation*}
		\mathcal{A}=\mathcal{Q}\ast\mathcal{R},
	\end{equation*}
	where $\mathcal{Q}\in\mathbb{R}^{n_1\times n_1 \times n_3}$ is orthogonal, and $\mathcal{R}\in\mathbb{R}^{n_1\times n_2 \times n_3}$ is an $f$-upper triangular tensor, i.e., each frontal slice is an upper triangular matrix.
\end{theorem}
\begin{definition} \label{TenFnorm}
(Tensor Frobenius norm \cite{Bouallala, Kolda_ten}). The Frobenius norm of a tensor $\mathcal{A}\in \mathbb{R}^{I \times J \times K}$ is defined as
\begin{equation*}
  {\left \|\mathcal{A} \right \|}_F=\sqrt{
    \sum_{i=1}^{I}
    \sum_{j=1}^{J}
    \sum_{k=1}^{K}
    {\left(\mathcal{A}\left(i,j,k\right)\right)}^2
}
=\operatorname{tr}\bigl(\mathcal{A}^T \ast \mathcal{A}\bigr).
\end{equation*}
\end{definition}
\begin{definition}\label{def9}
	(Tensor spectral norm \cite{t-specnorm}). The spectral norm of a tensor $\mathcal{A}\in\mathbb{R}^{n_1\times n_2 \times n_3}$ is defined as
	\begin{equation*}
		{\|\mathcal{A}\|}_2={\|{\rm bcirc}(\mathcal{A})\|}_2={\|{\rm bdiag}(\widehat{\mathcal{A}})\|}_2.
	\end{equation*}
\end{definition}
\begin{remark}\label{rem1}
{\it From ${\|{\rm bdiag}(\widehat{\mathcal{A}})\|}_2=\sigma_1({\rm bdiag}(\widehat{\mathcal{A}}))$, it follows that ${\|\mathcal{A}\|}_2={\|{\rm bdiag}(\widehat{\mathcal{A}})\|}_2=\max_{k=1,2,\ldots,n_3}{\|\widehat{\mathcal{A}}^{(k)}\|}_2$ {\rm\cite{t-specnorm2}}.}
\end{remark}

\subsection{Multilinear discriminant analysis}
This subsection gives a concise introduction to the multilinear discriminant analysis (MDA) approach built upon tensor t-product \cite{Bouallala}, together with a brief overview regarding the tensor trace-ratio problem.
Multilinear discriminant analysis serves as one generalization of classic linear discriminant analysis (LDA) \cite{LDA} toward tensor-form inputs. As a widely adopted dimensionality reduction tool, linear discriminant analysis mathematically corresponds to the well-studied trace-ratio optimization task. The core objective lies in seeking a projection matrix that enlarges the between-class scatter matrix and restricts the within-class scatter matrix \cite{LDA}. Similar to LDA, multilinear discriminant analysis pursues a projection tensor. Its target is to boost the between-class scatter tensor whilst lowering the magnitude of the within-class scatter tensor.

More detailed descriptions are given below. For a data tensor $\mathcal{X}\in \mathbb{R}^{n_1 \times m \times n_3}$ containing $m$ points samples together with label vector $\mathbf{y}$, multilinear discriminant analysis method aims at seeking a projection tensor $\mathcal{V}\in \mathbb{R}^{n_1 \times d \times n_3}$, where $d$ stands for the target feature dimension after original-data dimensionality reduction.
For notation purposes, the tensor $\mathcal{X}$ is expressed as $\mathcal{X}=\left[\overrightarrow{\mathcal{X}_1}, \overrightarrow{\mathcal{X}_2}, \dots, \overrightarrow{\mathcal{X}_m}\right] \in \mathbb{R}^{n_1 \times m \times n_3}$  to represent the collection of $m$ training instances. Suppose the whole data tensor is split into $c$ groups such that $ \mathcal{X} = \left[\mathcal{X}_1, \dots, \mathcal{X}_c \right] $, in which $\mathcal{X}_j\in\mathbb{R}^{n_1\times m_j \times n_3}$ corresponds to the $j$-th group and $m_j$ denotes its sample quantity. $\mathcal{C}_j \in\mathbb{R}^{n_1\times 1 \times n_3}$ represents the centroid tensor associated with $ \mathcal{X}_j $, while $\mathcal{C} \in\mathbb{R}^{n_1\times 1 \times n_3}$  denotes the overall centroid tensor over the training set. Define $\mathcal{E}_{m_j} \in\mathbb{R}^{m_j\times 1 \times n_3}$ satisfying $\mathcal{E}_{m_j}^{(1)}=\mathbf{1}_{m_j}$, $\mathcal{E}_{m_j}^{(i)}=\mathbf{0}_{m_j}$, $i=2,3,\ldots,n_3$, where $\mathbf{1}_{m_j}$ and $\mathbf{0}_{m_j}$ are the $m_j$-dimensional one-vector and zero-vector, respectively.

Accordingly, one can write the expression of the between-class scatter tensor $\mathcal{S}_b$ as
\begin{equation} \label{Sb}
	\mathcal{S}_b=\sum_{j=1}^{c}m_j\left(\mathcal{C}_j-\mathcal{C}\right)\ast{\left(\mathcal{C}_j-\mathcal{C}\right)}^T=\mathcal{H}_b\ast\mathcal{H}_b^T \in\mathbb{R}^{n_1 \times n_1 \times n_3},
\end{equation}
in which $\mathcal{H}_b=\left[\sqrt{m_1}(\mathcal{C}_1-\mathcal{C}),\sqrt{m_2}(\mathcal{C}_2-\mathcal{C}),\dots,\sqrt{m_c}(\mathcal{C}_c-\mathcal{C})\right]\in\mathbb{R}^{n_1 \times c \times n_3}$.
For the within-class scatter tensor $\mathcal{S}_w$, its expression is given by
\begin{equation*}
	\mathcal{S}_w=\sum_{j}^{c}\sum_{\overrightarrow{\mathcal{X}_i} \in \mathcal{X}_j}\left(\overrightarrow{\mathcal{X}_i}-\mathcal{C}_j \right)\ast{\left(\overrightarrow{\mathcal{X}_i}-\mathcal{C}_j \right)}^T=\mathcal{H}_w\ast\mathcal{H}_w^T \in\mathbb{R}^{n_1 \times n_1 \times n_3},
\end{equation*}
where $\mathcal{H}_w=\left[\mathcal{X}_1-\mathcal{C}_1\ast\mathcal{E}_{m_1}^T,\mathcal{X}_2-\mathcal{C}_2\ast\mathcal{E}_{m_2}^T,\dots,\mathcal{X}_c-\mathcal{C}_c\ast \mathcal{E}_{m_c}^T \right]\in\mathbb{R}^{n_1 \times m \times n_3}$. The total scatter tensor is defined as
\begin{equation}\label{ST}
	\mathcal{S}_t=\sum_{j}^{m}\left(\overrightarrow{\mathcal{X}_j}-\mathcal{C} \right)\ast{\left(\overrightarrow{\mathcal{X}_j}-\mathcal{C} \right)}^T=\mathcal{H}_t\ast\mathcal{H}_t^T \in\mathbb{R}^{n_1 \times n_1 \times n_3},
\end{equation}
where $\mathcal{H}_t=\left[\overrightarrow{\mathcal{X}_1}-\mathcal{C},\overrightarrow{\mathcal{X}_2}-\mathcal{C},\dots,\overrightarrow{\mathcal{X}_m}-\mathcal{C} \right]\in\mathbb{R}^{n_1 \times m \times n_3}$. Moreover, it holds that
\begin{equation}\label{St}
	\mathcal{S}_t=\mathcal{S}_w+\mathcal{S}_b.
\end{equation}

Then the multilinear discriminant analysis reduces to the following {\it tensor trace-ratio problem} \cite{Bouallala}:
\begin{equation}\label{equa2.3.1}
	\mathcal{V}^*=\arg \max_{\mathcal{V} \in \mathbb{R}^{n_1\times d\times n_3},\mathcal{V}^T\ast\mathcal{V}=\mathcal{I}_d} \frac{{\rm tr}\left(\mathcal{V}^T\ast\mathcal{S}_b\ast\mathcal{V}\right)}{{\rm tr}\left(\mathcal{V}^T\ast\mathcal{S}_w\ast\mathcal{V}\right)},
\end{equation}
where $\mathcal{I}_d$ denotes an identity tensor of size $d \times d \times n_3$.

In \cite{Bouallala}, the problem \eqref{equa2.3.1} was replaced by the following {\it tensor ratio-trace problem}:
\begin{equation} \label{equa2.3.2}
	\mathcal{V}^*=\arg \max_{\mathcal{V} \in \mathbb{R}^{n_1\times d\times n_3},\mathcal{V}^T\ast\mathcal{V}=\mathcal{I}_d}{\rm tr}\left({\left(\mathcal{V}^T\ast\mathcal{S}_w\ast\mathcal{V}\right)}^{-1}\ast\mathcal{V}^T\ast\mathcal{S}_b\ast\mathcal{V}\right) .
\end{equation}
Solving the subsequent generalized tensor eigenvalue problem yields the desired solution
\begin{equation}\label{equa2.3.3}
	\mathcal{S}_b\ast\mathcal{V}=\mathcal{S}_w\ast\mathcal{V}\ast\mathcal{D},
\end{equation}
where $\mathcal{D}\in \mathbb{R}^{d\times d\times n_3}$ represents an $f$-diagonal tensor. The optimal result for \eqref{equa2.3.2} corresponds to the $d$ eigenslices associated with the $d$ largest eigentube of \eqref{equa2.3.3}. However, the tensor trace-ratio problem \eqref{equa2.3.1} and the tensor ratio-trace problem \textcolor{red}{\eqref{equa2.3.2}} do not hold mathematical equivalence under general conditions.

\section{A randomized algorithm for tensor trace-ratio problem with applications to tensor dimensionality reduction}\label{sec3}

Within this section, we derive a closed-form solution and propose a randomized algorithm to tackle the tensor trace-ratio optimization problem \eqref{equa2.3.1}.
Combining \eqref{St} and \eqref{equa2.3.1}, we rewrite \eqref{equa2.3.1} as
\begin{equation}\label{equa3.1}
	\mathcal{V}^*=\arg \max_{\mathcal{V} \in \mathbb{R}^{n_1\times d\times n_3},\mathcal{V}^T\ast\mathcal{V}=\mathcal{I}_d} \frac{{\rm tr}\left(\mathcal{V}^T\ast\mathcal{S}_b\ast\mathcal{V}\right)}{{\rm tr}\left(\mathcal{V}^T\ast\mathcal{S}_t\ast\mathcal{V}\right)}.
\end{equation}

We mention that \eqref{equa2.3.1} and \eqref{equa3.1} share mathematical equivalence. Since $\mathcal{S}_w$, $\mathcal{S}_b$ and $\mathcal{S}_t$ are positive semi-definite tensors, \eqref{equa2.3.1} and \eqref{equa3.1} can be equivalently rewritten as
\begin{equation}\label{e4r}
\mathcal{V}^*=\arg \min_{\mathcal{V} \in \mathbb{R}^{n_1\times d\times n_3},\mathcal{V}^T\ast\mathcal{V}=\mathcal{I}_d} \frac{{\rm tr}\left(\mathcal{V}^T\ast\mathcal{S}_w\ast\mathcal{V}\right)}{{\rm tr}\left(\mathcal{V}^T\ast\mathcal{S}_b\ast\mathcal{V}\right)},
\end{equation}
and
\begin{equation}\label{e7r}
\mathcal{V}^*=\arg \min_{\mathcal{V} \in \mathbb{R}^{n_1\times d\times n_3},\mathcal{V}^T\ast\mathcal{V}=\mathcal{I}_d} \frac{{\rm tr}\left(\mathcal{V}^T\ast\mathcal{S}_t\ast\mathcal{V}\right)}{{\rm tr}\left(\mathcal{V}^T\ast\mathcal{S}_b\ast\mathcal{V}\right)},
\end{equation}
respectively. It follows from \eqref{St} that
\begin{equation*} \label{etr}
\frac{{\rm tr}\left(\mathcal{V}^T\ast\mathcal{S}_t\ast\mathcal{V}\right)}{{\rm tr}\left(\mathcal{V}^T\ast\mathcal{S}_b\ast\mathcal{V}\right)}
=1+
\frac{{\rm tr}\left(\mathcal{V}^T\ast\mathcal{S}_w\ast\mathcal{V}\right)}{{\rm tr}\left(\mathcal{V}^T\ast\mathcal{S}_b\ast\mathcal{V}\right)},
\end{equation*}
which shows the equivalence between \eqref{e4r} and \eqref{e7r}, as well as that between \eqref{equa2.3.1} and \eqref{equa3.1}.
Indeed, if $\mathcal{V} \in \mathcal{N}\left(\mathcal{S}_w\right) \setminus \mathcal{N}\left(\mathcal{S}_b\right)$, then \eqref{e4r} and \eqref{e7r} attain their minimum values simultaneously, which are $0$ and $1$, respectively.

Subsequently, we investigate efficient approaches for solving \eqref{equa3.1}, and provide several theoretical analyses to verify the rationality behind our strategy.

\subsection{A closed-form solution }\label{Sec3.1}
Let $\mathcal{E}_m \in\mathbb{R}^{m\times 1 \times n_3}$, with $\mathcal{E}_m^{(1)}=\mathbf{1}_m$, and $\mathcal{E}_m^{(i)}=\mathbf{0}_m$, $i=2,3,\ldots,n_3$. Consider
the centered data tensor
$$
\overline{\mathcal{X}}=\mathcal{X}\ast\left(\mathcal{I}_m-\mathcal{E}_m\ast\mathcal{E}_m^{T}/m\right)=\mathcal{X}\ast\mathcal{H},
$$
where $\mathcal{H}=\mathcal{I}_m-\mathcal{E}_m\ast\mathcal{E}_m^{T}/m \in \mathbb{R}^{m\times m \times n_3}$. Define the block diagonal tensor $\mathcal{W} \in \mathbb{R}^{m\times m \times n_3}$, with
\begin{equation}
\label{eqnw}
{W}^{(1)}= \begin{bmatrix}
	\frac{1}{m_1} \mathbf{1}_{m_1} \mathbf{1}_{m_1}^T & & & \\
	& \frac{1}{m_2} \mathbf{1}_{m_2} \mathbf{1}_{m_2}^T & & \\
	& & \ddots & \\
	& & & \frac{1}{m_c} \mathbf{1}_{m_c} \mathbf{1}_{m_c}^T
\end{bmatrix} ,
\quad \mathcal{W}^{(i)}=\mathbf{O},~ i=2,\dots,n_3,
\end{equation}
where $m_j,~j=1,2,\ldots,n_3$, denote the sample count belonging to the $j$-th class.
The lemma presented below reveals several key properties associated with $\mathcal{H}$ and $\mathcal{W}$.

\begin{lemma} \label{Lem1}
	 Using the notation introduced above, the following results hold
  \begin{enumerate}
    \item $\mathcal{H} \ast \mathcal{H} = \mathcal{H}$.
    \item $\mathcal{W} \ast \mathcal{E}_m = \mathcal{E}_m$.
    \item $\mathcal{H} \ast \mathcal{E}_m = \mathcal{O}$.
  \end{enumerate}
\end{lemma}
\begin{proof}
	(i) From $\mathcal{H}^{(1)}=\mathbf{I}_m-\mathbf{1}_m\mathbf{1}_m^T/m$ and $\mathcal{H}^{(i)}=\mathbf{O}$, $i=2,3,\ldots,n_3$, we obtain
	\begin{align*}
		\mathcal{H}\ast\mathcal{H}\!&=\! {\rm fold } \left( {\rm bcirc}\left(\mathcal{H}\right) {\rm unfold}\left(\mathcal{H}\right)\right)
		\!=\! {\rm fold}\left(\begin{bmatrix}
			\mathcal{H}^{(1)}& \mathbf{O} & \cdots & \mathbf{O}\\
			\mathbf{O}& \mathcal{H}^{(1)} & \mathbf{O} \\
			\vdots& \vdots& \ddots & \vdots \\
			\mathbf{O} & \mathbf{O}& \cdots& \mathcal{H}^{(1)}
		\end{bmatrix}
		\begin{bmatrix}
			\mathcal{H}^{(1)} \\
			\mathbf{O} \\
			\vdots \\
			\mathbf{O}
		\end{bmatrix}\right)\\
		\!&=\! {\rm fold}\left({\begin{bmatrix}
				\mathbf{I}_m-\mathbf{1}_m\mathbf{1}_m^T/m &\mathbf{O}&\cdots&\mathbf{O}\end{bmatrix}}^T\right)\\
\!&=\! {\rm fold}\left({\begin{bmatrix}
				\mathcal{H}^{(1)} &\mathbf{O}&\cdots&\mathbf{O}\end{bmatrix}}^T\right)
		=\mathcal{H}.
	\end{align*}
	(ii) We note that $\mathcal{W}^{(i)}\ast \mathcal{E}_m^{(j)}=\mathbf{O},~ i \neq 1 \text{ or } j \neq 1$, and
	\begin{align*}
		\mathcal{W}^{(1)}\ast \mathcal{E}_m^{(1)}=\begin{bmatrix}
			\frac{1}{m_1} \mathbf{1}_{m_1} \mathbf{1}_{m_1}^T & \mathbf{O}& \cdots& \mathbf{O}\\
			\mathbf{O}& \frac{1}{m_2} \mathbf{1}_{m_2} \mathbf{1}_{m_2}^T & \cdots&\mathbf{O} \\
			\vdots& \vdots& \ddots & \vdots\\
			\mathbf{O}& \mathbf{O}& \ldots& \frac{1}{m_c} \mathbf{1}_{m_c} \mathbf{1}_{m_c}^T
		\end{bmatrix}
		\mathbf{1}_m
		=\mathbf{1}_m.
	\end{align*}
	So we have that
	\begin{align*}
		\mathcal{W} \ast \mathcal{E}_m \!&=\! {\rm fold}\left( {\rm bcirc}\left(\mathcal{W}\right) {\rm unfold}\left(\mathcal{E}_m\right)\right)\\
		\!&=\! {\rm fold}\left(\begin{bmatrix}
			\mathcal{W}^{(1)}& \mathcal{W}^{(n_3)} & \cdots & \mathcal{W}^{(2)}\\
			\mathcal{W}^{(2)} & \mathcal{W}^{(1)} & \cdots & \mathcal{W}^{(3)} \\
			\vdots& \vdots& \ddots & \vdots \\
			\mathcal{W}^{(n_3)}& \mathcal{W}^{(n_3-1)}& \cdots& \mathcal{W}^{(1)}
		\end{bmatrix}
		\begin{bmatrix}
			\mathcal{E}_m^{(1)} \\
			\mathcal{E}_m^{(2)} \\
			\vdots \\
			\mathcal{E}_m^{(n_3)}
		\end{bmatrix}\right) \\
		\!&=\! {\rm fold}\left({\begin{bmatrix}
				\mathbf{1}_m^T &\mathbf{0}^T&\cdots&\mathbf{0}^T \end{bmatrix}}^T\right)
		= \mathcal{E}_m.
	\end{align*}
	(iii) From $\mathcal{H}^{(i)}\ast \mathcal{E}_m^{(j)}=\mathbf{O}~ (i \neq 1 \text{ or } j \neq 1)$, and
	\begin{align*}
		\begin{bmatrix}
\mathbf{I}_m-\mathbf{1}_m\mathbf{1}_m^T/m& \mathbf{O} & \cdots & \mathbf{O}\\
			\mathbf{O}& \mathbf{I}_m-\mathbf{1}_m\mathbf{1}_m^T/m & \cdots & \mathbf{O} \\
			\vdots& \vdots& \ddots & \vdots \\
			\mathbf{O} & \mathbf{O}& \cdots& \mathbf{I}_m-\mathbf{1}_m\mathbf{1}_m^T/m
		\end{bmatrix}
\begin{bmatrix}
\mathbf{1}_m \\
\mathbf{0} \\
\vdots \\
\mathbf{0}
\end{bmatrix}
=\mathbf{0},
	\end{align*}
	it follows that
	\begin{align*}
		\mathcal{H} \ast\mathcal{E}_m \!&=\! {\rm fold}\left( {\rm bcirc}(\mathcal{H}){\rm unfold}(\mathcal{E}_m)\right)\\
		\!&=\! {\rm fold}\left(\begin{bmatrix}
			\mathcal{H}^{(1)}& \mathcal{H}^{(n_3)} & \cdots & \mathcal{H}^{(2)}\\
			\mathcal{H}^{(2)} & \mathcal{H}^{(1)} & \cdots & \mathcal{H}^{(3)} \\
			\vdots& \vdots& \ddots & \vdots \\
			\mathcal{H}^{(n_3)}& \mathcal{H}^{(n_3-1)}& \cdots& \mathcal{H}^{(1)}
		\end{bmatrix}
		\begin{bmatrix}
			\mathcal{E}_m^{(1)} \\
			\mathcal{E}_m^{(2)} \\
			\vdots \\
			\mathcal{E}_m^{(n_3)}
		\end{bmatrix}\right)\\
		\!&=\! {\rm fold}\left({\begin{bmatrix}
				\mathbf{0}^T&\mathbf{0}^T &\ldots&\mathbf{0}^T \end{bmatrix}}^T\right)
		= \mathcal{O}.
	\end{align*}
\end{proof}

\begin{theorem}
   The three scatter tensors $\mathcal{S}_b$, $\mathcal{S}_w$, together with $\mathcal{S}_t$, admit the equivalent expressions
 \begin{equation}\label{eqsb}
  \mathcal{S}_b=\overline{\mathcal{X}} \ast \mathcal{W} \ast \overline{\mathcal{X}}^T,
  \end{equation}
 \begin{equation}\label{eqsw}
  \mathcal{S}_w=\overline{\mathcal{X}} \ast (\mathcal{I}_m-\mathcal{W}) \ast \overline{\mathcal{X}}^T,
\end{equation}
  \begin{equation}\label{eqst}
  \mathcal{S}_t=\overline{\mathcal{X}} \ast \overline{\mathcal{X}}^T.
  \end{equation}
\end{theorem}
\begin{proof}
It suffices to establish these results within the Fourier domain. Observe that the centered-data tensor admits the form
\begin{align*}
\overline{\mathcal{X}}
\!&=\!
\begin{bmatrix}
			\overrightarrow{\mathcal{X}_1}-\mathcal{C},& \overrightarrow{\mathcal{X}_2}-\mathcal{C},& \ldots,&\overrightarrow{\mathcal{X}_m}-\mathcal{C}
		\end{bmatrix}\\
\!&=\!
\begin{bmatrix}
			\overline{\overrightarrow{\mathcal{X}_1}},& \overline{\overrightarrow{\mathcal{X}_2}},& \ldots,&\overline{\overrightarrow{\mathcal{X}_m}}
		\end{bmatrix}\\
\!&=\!
\begin{bmatrix}
			\overline{\mathcal{X}_1},& \overline{\mathcal{X}_2},& \ldots,&\overline{\mathcal{X}_c}
		\end{bmatrix},
\end{align*}
where one may interpret $\overline{\overrightarrow{\mathcal{X}_i}} (i=1,2,\ldots,m)$ as the $i$-th centered sample tensor, and $\overline{\mathcal{X}_j} (j=1,2,\ldots,c)$ as the centered data tensor associated with the $j$-th class.

By virtue of \eqref{Sb}, \eqref{eqnw}, along with the inherent properties of the t-product, we obtain
  \begin{align*}
    \widehat{\mathcal{S}}_b^{(k)}\!&=\!
    \sum_{j=1}^{c}m_j\left( \widehat{\mathcal{C}}_j^{(k)}-\widehat{\mathcal{C}}^{(k)} \right) {\left( \widehat{\mathcal{C}}_j^{(k)}-\widehat{\mathcal{C}}^{(k)} \right)}^H \\
    \!&=\!
    \sum_{j=1}^{c} m_j \left(\frac{1}{m_j} \sum_{i=1}^{m_j}  \left( \overrightarrow{\left[ \widehat{\mathcal{X}}_j^{(k)} \right]}_i - \widehat{\mathcal{C}}^{(k)}\right)\right) {\left(\sum_{i=1}^{m_j} \frac{1}{m_j} \left( \overrightarrow{\left[ \widehat{\mathcal{X}}_j^{(k)} \right]}_i - \widehat{\mathcal{C}}^{(k)}\right)\right)}^H\\
    \!&=\!
    \sum_{j=1}^{c} \frac{1}{m_j}  \left(\sum_{i=1}^{m_j} \overline{\overrightarrow{\left[ \widehat{\mathcal{X}}_j^{(k)} \right]}_i}
    \sum_{i=1}^{m_j} {\left(\overline{\overrightarrow{\left[ \widehat{\mathcal{X}}_j^{(k)} \right]}_i}\right)}^H \right) \\
    \!&=\!
    \sum_{j=1}^{c} \widehat{\overline{\mathcal{X}}}_j^{(k)} \left(\mathbf{1}_{m_j}\mathbf{1}_{m_j}^T/{m_j} \right) {\left(\widehat{{\overline{\mathcal{X}_j}}}^H\right)}^{(k)}\\
  \!&=\!\widehat{\overline{\mathcal{X}}}^{(k)}\widehat{\mathcal{W}}^{(k)}{\left(\widehat{{\overline{\mathcal{X}}}}^H\right)}^{(k)},\quad k=1,2,\ldots,n_3,
  \end{align*}
  which yields \eqref{eqsb}. Here $\overrightarrow{\left[ \widehat{\mathcal{X}}_j^{(k)} \right]}_i,~i=1,2,\ldots,m_j,~j=1,2,\ldots,c,~k=1,2,\ldots,n_3$, are the $i$-th sample data in the $k$-th frontal slice of the $j$-th class data tensor in the Fourier domain, while $\overline{\overrightarrow{\left[ \widehat{\mathcal{X}}_j^{(k)} \right]}_i}$ stands for the centered form associated with $\overrightarrow{\left[ \widehat{\mathcal{X}}_j^{(k)} \right]}_i$.

 For \eqref{eqst}, it follows from \eqref{ST} that
  \begin{align*}
    \widehat{\mathcal{S}}_t^{(k)}\!&=\!
    \sum_{j=1}^{m}\left(\widehat{\overrightarrow{\mathcal{X}_j}}^{(k)}-\mathcal{C}^{(k)}\right){\left(\widehat{\overrightarrow{\mathcal{X}_j}}^{(k)}-\mathcal{C}^{(k)}\right)}^H\\
\!&=\!
\sum_{j=1}^{m} \widehat{\overline{\mathcal{X}_j}}^{(k)} {\left(\widehat{{\overline{\mathcal{X}_j}}}^H\right)}^{(k)}\\
\!&=\! \widehat{\overline{\mathcal{X}}}^{(k)}
{\left(\widehat{{\overline{\mathcal{X}}}}^H\right)}^{(k)},~k=1,2,\ldots,n_3.
  \end{align*}
Finally, from \eqref{St}, we arrive at
$$
\mathcal{S}_w=\mathcal{S}_t-\mathcal{S}_b=\overline{\mathcal{X}} \ast (\mathcal{I}_m-\mathcal{W}) \ast \overline{\mathcal{X}}^T,
$$
which yields \eqref{eqsw}.
\end{proof}

In terms of Lemma \ref{Lem1}, we can reformulate $\mathcal{S}_b$ and $\mathcal{S}_t$ as
\begin{equation}\label{equa3.2}
	\mathcal{S}_b=\mathcal{X} \ast \mathcal{H} \ast \mathcal{W} \ast \mathcal{H} \ast \mathcal{X}^T,\quad
	\mathcal{S}_t=\mathcal{X} \ast \mathcal{H} \ast \mathcal{H} \ast \mathcal{X}^T=\mathcal{X}  \ast \mathcal{H} \ast \mathcal{X}^T,
\end{equation}
in which we make use of the identity $ \mathcal{H} \ast \mathcal{H}=\mathcal{H}$. Moreover, we have from $\mathcal{W} \ast \mathcal{E}_m=\mathcal{E}_m$ that
\begin{equation}\label{equa3.3}
	\mathcal{S}_w=\mathcal{S}_t-\mathcal{S}_b=\mathcal{X}  \ast \left(\mathcal{H}- \mathcal{H} \ast \mathcal{W} \ast \mathcal{H}\right) \ast \mathcal{X}^T=\mathcal{X}  \ast \left(\mathcal{I}_m-\mathcal{W}\right) \ast \mathcal{X}^T.
\end{equation}
\begin{remark}
{\it Without loss of generality, as was done in the matrix case \cite{LDA}, we fix the reduced dimension $d$ equal to $c-1$ in this paper.
It is seen from \eqref{St} that
$$
{\rm tr}\left(\mathcal{V}^T\ast\mathcal{S}_t\ast\mathcal{V}\right)={\rm tr}\left(\mathcal{V}^T\ast\mathcal{S}_b\ast\mathcal{V}\right)+{\rm tr}\left(\mathcal{V}^T\ast\mathcal{S}_w\ast\mathcal{V}\right).
$$
Motivated by \cite{Shi}, suppose ${\rm tr}\left(\mathcal{V}^T\ast\mathcal{S}_w\ast\mathcal{V}\right)=0$ while ${\rm tr}\left(\mathcal{V}^T\ast\mathcal{S}_b\ast\mathcal{V}\right)\ne 0$. Under such conditions, \eqref{equa3.1} attains its maximum value equal to $1$. As a result, any orthogonal basis in $\mathcal{N}\left(\mathcal{S}_w\right) \setminus \mathcal{N}\left(\mathcal{S}_t\right)$ of dimension $(c-1)$ is a solution to \eqref{equa3.1}.}
\end{remark}
Motivated by this observation, we first aim to construct a basis spanning a $c$-dimension subspace within $\mathcal{N}(\mathcal{S}_w)$.
\begin{theorem}
	Denote by $\mathcal{Y} \in \mathbb{R}^{m \times c \times n_3}$, with
	\begin{eqnarray} \label{equa3.4}
		\begin{array}{cll}
			\mathcal{Y}^{(1)}&=&  \begin{bmatrix}
				\mathbf {1}_{m_1} & \mathbf {0} & \cdots & \mathbf {0}\\
				\mathbf {0} & \mathbf {1}_{m_2}  & \cdots & \mathbf {0} \\
				\vdots& \vdots& \ddots & \vdots \\
				\mathbf {0}& \mathbf {0}& \cdots& \mathbf {1}_{m_c} \end{bmatrix},\quad\mathcal{Y}^{(i)}=\mathbf{O},~i=2,\dots,n_3.
		\end{array}
	\end{eqnarray}
Provided that ${\rm bcirc}\left(\mathcal{X}^T\right)$ has full row rank, $\mathcal{F}={\left(\mathcal{X}^T\right)}^\dagger \ast \mathcal{Y}$ forms the basis of a $c$-dimension subspace lying in $\mathcal{N}\left(\mathcal{S}_w\right)$. Here, $c$ represents the tensor tubal rank, and ${\left(\mathcal{X}^T\right)}^\dagger $ stands for the Moore-Penrose inverse associated with $\mathcal{X}^T$.
\end{theorem}
\begin{proof}
Let
	\begin{equation*}
		\mathcal{Y}=\begin{bmatrix}
			\overrightarrow{\mathcal{Y}_1},& \overrightarrow{\mathcal{Y}_2},& \ldots,&\overrightarrow{\mathcal{Y}_c}
		\end{bmatrix},
	\end{equation*}
	in which $\overrightarrow{\mathcal{Y}_i}\in \mathbb{R}^{m \times 1 \times n_3},~i=1,2,\ldots,c$, with
$$
\overrightarrow{\mathcal{Y}_i}^{(1)}={\left[\mathbf {0}^T,\ldots,\mathbf {1}_{m_i}^T,\ldots, \mathbf {0}^T\right]}^T,
$$
and
$$
\overrightarrow{\mathcal{Y}_i}^{(j)}={\left[\mathbf {0}^T,\ldots,\mathbf {0}^T\right]}^T,~j=2,3,\ldots,n_3.
$$
Since ${\rm bcirc}\left(\mathcal{X}^T\right)$  possesses full row-rank property, together with the relation $\mathcal{W}\ast \mathcal{Y}=\mathcal{Y}$, one obtains
	\begin{equation*}\label{equa3.5}
		\mathcal{S}_w \ast \mathcal{F}=\mathcal{X} \ast \left(\mathcal{I}_m-\mathcal{W}\right)\ast \mathcal{X}^T\ast \left({\left(\mathcal{X}^T\right)}^\dagger\ast\mathcal{Y}\right)=\mathcal{X} \ast\left(\mathcal{I}_m-\mathcal{W}\right)\ast \mathcal{Y}=\mathcal{O},
	\end{equation*}
and $\mathcal{F}$ lies within the null space associated with $\mathcal{S}_w$.

We then demonstrate that tensor $\mathcal{F}$ possesses tubal rank equal to $c$.
Note that $\mathcal{F}={\left(\mathcal{X}^T\right)}^\dagger \ast \mathcal{Y}$, which implies $\mathcal{X}^T \ast \mathcal{F}= \mathcal{Y}$. Accordingly,
\begin{equation} \label{yier1}
  {\rm bcirc}\left( \mathcal{X}^T \right) {\rm unfold}\left( \mathcal{F}\right)={\rm unfold}\left( \mathcal{Y} \right),
\end{equation}
\eqref{yier1} may be recast into 
\begin{equation*} \label{yier2}
  \operatorname{unfold}\left( \mathcal{F} \right)
  = \left( \operatorname{bcirc}\left( \mathcal{X}^T \right) \right)^\dagger
    \operatorname{unfold}\left( \mathcal{Y} \right).
\end{equation*}
Given that $\left( \operatorname{bcirc}\left( \mathcal{X}^T \right) \right)^\dagger$ together with $\operatorname{unfold}\left( \mathcal{Y} \right)$ possess column-full-rank property, $\operatorname{unfold}\left( \mathcal{F} \right)$ inherits the column full-rank condition. Thus, it follows from Definition \ref{def8} that tensor $\mathcal{F}$ has tubal rank equal to $c$.
\end{proof}
Secondly, we eliminate components belonging to $\mathcal{N}\left(\mathcal{S}_t\right)$ from $span\left\{\mathcal{F}\right\}$ to get a closed form solution of \eqref{equa3.1}.
\begin{theorem} \label{Theo3.2}
	Let $\underline{\mathcal{Y}}=\left(\mathcal{I}_m-\mathcal{E}_m\ast\mathcal{E}_m^{T} /m \right)\ast\mathcal{Y}$. Assume that $\overline{\mathcal{Y}}\in \mathbb{R}^{m \times (c-1) \times n_3}$ provides an orthonormal basis spanning ${\rm span}\left\{\underline{\mathcal{Y}}\right\}$. Provided that ${\rm bcirc}\left(\mathcal{X}^T\right)$ has full-row rank, we have
	\begin{align}\label{equa Theo3.2}
		\overline{\mathcal{F}}={\left(\mathcal{X}^T\right)}^\dagger \ast \overline{\mathcal{Y}} \in \mathcal{N}\left(\mathcal{S}_w\right) \setminus \mathcal{N}\left(\mathcal{S}_t\right).
	\end{align}
	Specifically, if $\mathcal{F}_Q \in \mathbb{R}^{n_1 \times (c-1) \times n_3}$ stands for the $Q$-factor obtained from the economized QR factorization of $\overline{\mathcal{F}}$, $\mathcal{F}_Q $ constitutes one solution of \eqref{equa3.1}.
\end{theorem}
\begin{proof}
	 From  $\mathcal{X}^T\ast{\left(\mathcal{X}^T\right)}^\dagger=\mathcal{I}_m$ together with $\mathcal{W}\ast \mathcal{E}_m=\mathcal{E}_m$, we obtain
	\begin{align}\label{eq17}
		\mathcal{S}_w \ast\left({\left(\mathcal{X}^T\right)}^\dagger\ast\underline{\mathcal{Y}}\right)
		\!&=\!\mathcal{X} \ast\left(\mathcal{I}_m-\mathcal{W}\right) \ast\mathcal{X}^T \ast{\left(\mathcal{X}^T\right)}^\dagger\ast\underline{\mathcal{Y}}\nonumber\\
		\!&=\!\mathcal{X} \ast\left(\mathcal{I}_m-\mathcal{W}\right) \ast\mathcal{X}^T \ast{\left(\mathcal{X}^T\right)}^\dagger \ast\left(\mathcal{I}_m-\mathcal{E}_m\ast\mathcal{E}_m^{T}/m \right)\ast\mathcal{Y}\nonumber\\
		\!&=\!\mathcal{X} \ast\left(\mathcal{I}_m-\mathcal{W}\right)\ast\left(\mathcal{I}_m-\mathcal{E}_m\ast\mathcal{E}_m^{T}/m \right)\ast\mathcal{Y}\nonumber\\
		\!&=\!\mathcal{O},
	\end{align}
	and ${\left(\mathcal{X}^T\right)}^\dagger\ast\underline{\mathcal{Y}}$ belongs to the null space associated with $\mathcal{S}_w$.

On the other hand, we have from Lemma \ref{Lem1} that
	\begin{align*}
		\mathcal{S}_t \ast\left({\left(\mathcal{X}^T\right)}^\dagger\ast\underline{\mathcal{Y}}\right)
		\!&=\!\mathcal{X}\ast \mathcal{H} \ast \mathcal{X}^T\ast\left({\left(\mathcal{X}^T\right)}^\dagger\ast\underline{\mathcal{Y}}\right)\\
		\!&=\!\mathcal{X}\ast\mathcal{H} \ast\left(\mathcal{I}_m-\mathcal{E}_m\ast\mathcal{E}_m^{T}/m \right)\ast\mathcal{Y}\\
		\!&=\!\mathcal{X}\ast \mathcal{H}\ast\mathcal{Y},
	\end{align*}
in which we used $\mathcal{H}\ast\mathcal{E}_m=\mathcal{O}$.
	It follows from $\mathcal{N}\left(\mathcal{H}^{(1)}\right)={\rm span}\left\{\mathbf{1}_m\right\}$ that
\begin{align*}
\mathcal{H}\ast\overrightarrow{\mathcal{Y}_i}\!&=\! {\rm fold}\left({\rm bcirc}\left(\mathcal{H}\right){\rm  unfold}\left(\overrightarrow{\mathcal{Y}_i}\right)\right)\\
		\!&=\! {\rm fold}\left(
		\begin{bmatrix}
			\mathbf{I}_m-\mathbf{1}_m\mathbf{1}_m^T/m& \mathbf{O} & \cdots & \mathbf{O}\\
			\mathbf{O}& \mathbf{I}_m-\mathbf{1}_m\mathbf{1}_m^T/m & \cdots & \mathbf{O} \\
			\vdots& \vdots& \ddots & \vdots \\
			\mathbf{O} & \mathbf{O}& \cdots& \mathbf{I}_m-\mathbf{1}_m\mathbf{1}_m^T/m
		\end{bmatrix}
		\begin{bmatrix}
			{[\mathbf {0}^T,\ldots,\mathbf {1}_{m_i}^T,\ldots, \mathbf {0}^T]}^T \\
			\mathbf{0}\\
			\vdots \\
			\mathbf{0}
		\end{bmatrix}\right)\\
		\!&=\! {\rm fold}\left({\begin{bmatrix}
				 \left(\mathbf{I}_m-\mathbf{1}_m\mathbf{1}_m^T/m \right){[\mathbf {0}^T,\ldots,\mathbf {1}_{m_i}^T,\ldots, \mathbf {0}^T]}^T, &\mathbf{0}^T,&\cdots& ,\mathbf{0}^T \end{bmatrix}}^T\right)
		\ne \mathcal{O},~i=1,2,\ldots,c.
\end{align*}
	Consequently,
\begin{equation*}
  \mathcal{H}\ast\mathcal{Y}=[\mathcal{H}\ast\overrightarrow{\mathcal{Y}_1},\mathcal{H}\ast\overrightarrow{\mathcal{Y}_2},\ldots,\mathcal{H}\ast\overrightarrow{\mathcal{Y}_c}]\ne \mathcal{O},
\end{equation*}
and $\mathcal{X}\ast \mathcal{H}\ast\mathcal{Y}\ne \mathcal{O}$. As a result,
\begin{equation}\label{eq18}
  \mathcal{S}_t \ast\left({\left(\mathcal{X}^T\right)}^\dagger\ast\underline{\mathcal{Y}}\right)\ne \mathcal{O}.
\end{equation}
One may draw this conclusion based on formulas \eqref{eq17} and \eqref{eq18}:
	\begin{align*}
		{\left(\mathcal{X}^T\right)}^\dagger \ast \underline{\mathcal{Y}} \in \mathcal{N}\left(\mathcal{S}_w\right) \setminus \mathcal{N}\left
(\mathcal{S}_t\right).
	\end{align*}

Finally, given that $\overline{\mathcal{Y}}$ forms an orthonormal basis spanning ${\rm span}\left\{\underline{\mathcal{Y}}\right\}$, one obtains
$$
{\left(\mathcal{X}^T\right)}^\dagger \ast \overline{\mathcal{Y}} \in {\rm span}\left\{{\left(\mathcal{X}^T\right)}^\dagger \ast \underline{\mathcal{Y}}\right\},
$$
 and
 $$
 \overline{\mathcal{F}}={\left(\mathcal{X}^T\right)}^\dagger \ast \overline{\mathcal{Y}} \in \mathcal{N}\left(\mathcal{S}_w\right) \setminus \mathcal{N}\left(\mathcal{S}_t\right).
 $$
For the tensor $\mathcal{F}_Q \in \mathbb{R}^{n_1 \times (c-1) \times n_3}$,
which corresponds to the $Q$-component obtained upon performing reduced QR decomposition on $\overline{\mathcal{F}}$, \eqref{equa3.1} holds for $\mathcal{F}_{Q}$ as well.
\end{proof}

\begin{remark}
	{\it Theorem \ref{Theo3.2} indicates that if ${\rm bcirc}\left(\mathcal{X}^T\right)$ is of full row rank, \eqref{equa Theo3.2} is a solution to \eqref{equa3.1}. In real applications, this condition may not always be satisfied. However, numerical experiments show that \eqref{equa Theo3.2} still works very well; please see Section \ref{sec5}.
When the full row rank condition is violated, rigorous theoretical optimality guarantees may fail to hold. Under such circumstances, \eqref{equa Theo3.2} is merely treated as an empirical approximation.
}
\end{remark}

\subsection{A randomized algorithm for the tensor trace-ratio problem}\label{Sec3.2}
A direct calculation of the projection tensor in problem \eqref{equa3.1} for high-dimensional large-scale tensor data is computationally intractable. Recently, some inner-outer iteration methods are proposed for this type of problem. More precisely, inner iterations are used to solve the generalized eigenvalue problem for eigenslices, while outer iterations compute the trace-ratio objective as discussed in \cite{Tensor_product1, Tracevs, Bouallala, Zahir}. However, for large-scale tensor data, the workload of these methods will be heavy if outer-loop iterations grow excessively. It is therefore worthwhile to explore alternative approaches for large-scale tensor trace-ratio tasks without relying on the nested inner-outer iteration architecture.

Theorem \ref{Theo3.2} yields a closed-form solution for \eqref{equa3.1}, offering an alternative route toward the tensor trace-ratio problem. Indeed, tensor SVD under t-product (t-SVD) \cite{t-product1} is applicable for tackling the tensor least-squares formulation given in \eqref{equa Theo3.2}. However, it may suffer from the curse of dimensionality for large-scale and higher-order tensor data. To overcome this difficulty, we consider the randomized t-SVD instead, which exploits random sampling for approximating the low-rank tensor structure \cite{t-rsvd}. Algorithm \ref{Alg3} lists a randomized t-SVD accelerated by the power iteration.
\begin{algorithm}
	\caption{A randomized and truncated tensor SVD under t-product \cite{t-rsvd}}\label{Alg3}
	\begin{algorithmic}
		\item  {{\bf Input:} A tensor $\mathcal{X}\in \mathbb{R}^{n_1 \times m \times n_3}$, desired truncation rank $r$ satisfying $r\ll m$, oversampling parameter $p$ (typical choice $p\ge 4$), and power-iteration parameter $q$ (typical choice $q=2$); \\
			{\bf Output:} A rank-$r$ approximation $\widetilde{\mathcal{X}}^T $ to $\mathcal{X}^T$;} \\
		\item {{\bf 1.} Generate a Gaussian random tensor $\mathcal{W}\in\mathbb{R}^{m\times (r+p)\times n_3}$;}\\
		\item {{\bf 2.} Form a random projection of tensor $\mathcal{X}$ as $\mathcal{G}={(\mathcal{X}\ast\mathcal{X}^T)}^q\ast\mathcal{X}\ast\mathcal{W}$;}\\
		\item {{\bf 3.} Obtain tensor $\mathcal{Q}_1\in \mathbb{R}^{n_1 \times (r+p) \times n_3}$ via the t-QR factorization \cite{t-product2} applied to $\mathcal{G}$;} \\
		\item {{\bf 4.} Define tensor $\mathcal{B}=\mathcal{X}^T\ast\mathcal{Q}_1\in\mathbb{R}^{m \times (r+p) \times n_3}$, perform t-SVD decomposition $\mathcal{B}=\mathcal{U}_2\ast\mathcal{S}_1\ast\widetilde{\mathcal{V}}^T$, then assign $\mathcal{V}_1=\mathcal{Q}_1\ast\widetilde{\mathcal{V}}$;}\\
		\item {{\bf 5.} Set $\widetilde{\mathcal{U}}_r=\mathcal{U}_2(:,1:r,:)$, $\widetilde{\mathcal{V}}_r=\mathcal{V}_1(:,1:r,:)$, $\widetilde{\mathcal{S}}_r=\mathcal{S}_{\textcolor{red}{1}}(1:r,1:r,:)$, and assign $\widetilde{\mathcal{X}}^T=\widetilde{\mathcal{U}}_r\ast\widetilde{\mathcal{S}}_r\ast\widetilde{\mathcal{V}}_r^T$.}\\
	\end{algorithmic}
\end{algorithm}

With $\widetilde{\mathcal{X}}^T$ computed by using Algorithm \ref{Alg3}, then \eqref{equa Theo3.2} is transformed into the subsequent tensor least-squares problem:
\begin{equation}\label{equa3.2.1}
	\widetilde{\mathcal{F}}=\min_{\overline{\mathcal{F}}\in\mathbb{R}^{n_1\times {(c-1)}\times n_3}} {\left\|\widetilde{\mathcal{X}}^T\ast \overline{\mathcal{F}}-\overline{\mathcal{Y}}\right\|}_F.
\end{equation}
Algorithm \ref{Alg4} outlines the details of the main algorithm in this work.
\begin{algorithm}
	\caption{Randomized solver for \eqref{equa3.1}}\label{Alg4}
	\begin{algorithmic}
		\item  {{\bf Input:} A tensor $\mathcal{X}\in \mathbb{R}^{n_1 \times m \times n_3}$, desired truncation rank $r$ satisfying $r\ll m$, oversampling parameter $p$ (typical selection $p\ge 4$), together with power-iteration parameter $q$ (typical selection $q=2$); \\
			{\bf Output:} Approximate solution $\widetilde{\mathcal{F}}_Q $;} \\
		\item {{\bf 1.} Apply Algorithm \ref{Alg3} to obtain a rank-$r$ approximant for $\mathcal{X}^T$ $\colon$ $\widetilde{\mathcal{X}}^T=\widetilde{\mathcal{U}}_r\ast\widetilde{\mathcal{S}}_r\ast\widetilde{\mathcal{V}}_r^T$;}\\
		\item {{\bf 2.} Let $\widetilde{\mathcal{F}}={(\widetilde{\mathcal{X}}^T)}^\dagger\ast \overline{\mathcal{Y}}=\widetilde{\mathcal{V}}_r\ast\widetilde{\mathcal{S}}_r^\dagger\ast\widetilde{\mathcal{U}}_r^T\ast \overline{\mathcal{Y}}$;}\\
		\item {{\bf 3.} Compute $\widetilde{\mathcal{F}}_Q$ by computing the t-QR factorization \cite{t-product2} of the tensor $\widetilde{\mathcal{F}}$.}
	\end{algorithmic}
\end{algorithm}

\subsection{Theoretical analysis} \label{Sec3.3}
This subsection investigates the Frobenius-norm error ${\big\| \widetilde{\mathcal{F}}-\overline{\mathcal{F}}\big\|}_F$ between the closed-form solution $\overline{\mathcal{F}}$ derived in \eqref{equa Theo3.2} and the numerical solution $\widetilde{\mathcal{F}}$ yielded by solving \eqref{equa3.2.1}.
Let $\mathcal{X}_r^T$ denote the optimal tubal rank-$r$ approximation to $\mathcal{X}^T$ under the tensor Frobenius norm \cite{t-product1}, where $r$ corresponds to a user-defined input parameter.
Therefore, 
$\overline{\mathcal{F}}_r={\left(\mathcal{X}_r^T\right)}^\dagger\ast \overline{\mathcal{Y}}$ 
is an approximation to $\overline{\mathcal{F}}={\left(\mathcal{X}^T\right)}^\dagger\ast \overline{\mathcal{Y}}$. From Algorithm \ref{Alg3}, the rank-$r$ tensor $\widetilde{\mathcal{X}}^T$ is an approximation to $\mathcal{X}_r^T$. As a result, $\widetilde{\mathcal{F}}={\left(\widetilde{\mathcal{X}}^T\right)}^\dagger\ast \overline{\mathcal{Y}}$ approximates $\overline{\mathcal{F}}_r$.
We notice that
\begin{equation}
{\left\| \widetilde{\mathcal{F}}-\overline{\mathcal{F}}\right\|}_F \le
{\left\| \overline{\mathcal{F}}-\overline{\mathcal{F}}_r\right\|}_F+{\left\| \overline{\mathcal{F}}_r- \widetilde{\mathcal{F}}\right\|}_F,
\end{equation}
and it is sufficient to analyze the term ${\left\| \overline{\mathcal{F}}_r- \widetilde{\mathcal{F}}\right\|}_F$.
We turn our attention to discrepancy between the two subspaces $\widetilde{\mathbb{S}}=\operatorname{span}\left\{{\left(\widetilde{\mathcal{X}}^T\right)}^\dagger\ast \overline{\mathcal{Y}}\right\}$ and ${\mathbb{S}}=\operatorname{span}\left\{{\left(\mathcal{X}_r^T\right)}^\dagger\ast \overline{\mathcal{Y}}\right\}$.

We first present definitions on  tensor orthogonal projector and sine of the angles between two tensor subspaces.
\begin{definition}
  Let $\mathbb{S} \subseteq \mathbb{R}^{m \times 1 \times n}$ be a tensor subspace. Then $\mathcal{P} \in \mathbb{R}^{m \times m \times n}$ is the tensor orthogonal projector onto $\mathbb{S}$ if $\operatorname{span}(\mathcal{P}) = \mathbb{S}$, $\mathcal{P}^2 = \mathcal{P} \ast \mathcal{P} = \mathcal{P}$, and $\mathcal{P}^T = \mathcal{P}$.
\end{definition}
\begin{definition}\label{Def 3.1}
 Suppose $\mathbb{S}_1$ and $\mathbb{S}_2$ are two tensor subspaces of dimension $\mathbb{R}^{m \times 1 \times n}$ and $\dim(\mathbb{S}_1)=\dim(\mathbb{S}_2)$.
Then, define the sine of the angle between the two subspaces as
\begin{equation*}	
\sin\angle\left(\mathbb{S}_1,\mathbb{S}_2\right)={\left\|\mathcal{P}_1-\mathcal{P}_2\right\|}_2,
\end{equation*}
in which $\mathcal{P}_1$ and $\mathcal{P}_2$ serve as orthogonal projectors mapping onto tensor subspaces $\mathbb{S}_1$ and $\mathbb{S}_2$, respectively.
\end{definition}

To prepare for the subsequent main theorem, we start by providing a set of auxiliary lemmas for later derivations.
\begin{lemma}\label{lem 3.3}
	Suppose $\mathcal{P}\in \mathbb{R}^{n_1 \times n_2 \times n_3}$ represents an orthogonal projector tensor, while $\mathcal{M}\in \mathbb{R}^{n_2 \times l \times n_3}$ denotes another tensor. Then we have
	\begin{align} \label{eqle3}	
     {\left\|\mathcal{P}\ast\mathcal{M}\right\|}_2\leq{\left\|\mathcal{P}\ast{(\mathcal{M}\ast\mathcal{M}^H)}^q\ast\mathcal{M}\right\|}_2^{1/{2q+1}},
	\end{align}
where $q\geq 0$ is an integer.
\end{lemma}
\begin{proof}
	Suppose $\mathcal{M}=\mathcal{U}\ast\mathcal{S}\ast\mathcal{V}^H$ denotes its t-SVD decomposition. Within the Fourier transformed domain, one can derive the equality
	\begin{equation*}
		\widehat{\mathcal{P}}^{(i)}\widehat{\mathcal{M}}^{(i)}=\widehat{\mathcal{P}}^{(i)}\widehat{\mathcal{U}}^{(i)}\widehat{\mathcal{S}}^{(i)}\widehat{\mathcal{V}^T}^{(i)},~i=1,2,\ldots,n_3.
	\end{equation*}
	Furthermore, we obtain from \cite[Proposition 8.6]{Halko} that
	\begin{equation}\label{equa3.3.11}
		{\left\|\widehat{\mathcal{P}}^{(i)}\widehat{\mathcal{M}}^{(i)}\right\|}_2^{2(2q+1)}\leq{\left\|\widehat{\mathcal{P}}^{(i)}  {\left[\widehat{\mathcal{M}}^{(i)}{\left(\widehat{\mathcal{M}}^{(i)}\right)}^H\right]}^q  \widehat{\mathcal{M}}^{(i)}\right\|}_2^2, ~i=1,2,\ldots,n_3.
	\end{equation}
	Thus,
	\begin{equation*}
		{\left\|\mathcal{P}\ast\mathcal{M}\right\|}_2^{2(2q+1)}=\max_{i=1,2,\ldots,n_3}{\left\|\widehat{\mathcal{P}}^{(i)}\widehat{\mathcal{M}}^{(i)}\right\|}_2^{2(2q+1)},
	\end{equation*}
	and
	\begin{equation*}
		{\left\|\mathcal{P}\ast{\left(\mathcal{M}\ast\mathcal{M}^H\right)}^q\ast\mathcal{M}\right\|}_2^2=\max_{i=1,2,\ldots,n_3}{\left\|\widehat{\mathcal{P}}^{(i)}  {\left[\widehat{\mathcal{M}}^{(i)}{\left(\widehat{\mathcal{M}}^{(i)}\right)}^H\right]}^q  \widehat{\mathcal{M}}^{(i)}\right\|}_2^2.
	\end{equation*}
	So it follows from \eqref{equa3.3.11}  that
	\begin{equation*}
		\max_{i=1,2,\ldots,n_3}{\left\|\widehat{\mathcal{P}}^{(i)}\widehat{\mathcal{M}}^{(i)}\right\|}_2^{2(2q+1)}\leq\max_{i=1,2,\ldots,n_3}{\left\|\widehat{\mathcal{P}}^{(i)}  {\left[\widehat{\mathcal{M}}^{(i)}{\left(\widehat{\mathcal{M}}^{(i)}\right)}^H\right]}^q  \widehat{\mathcal{M}}^{(i)}\right\|}_2^2,
	\end{equation*}
	which gives \eqref{eqle3}.	
\end{proof}

\begin{lemma}\label{lem3.4}
Let $\mathcal{X}\in \mathbb{R}^{n_1 \times n_2 \times n_3}$, and consider an approximant $\widetilde{\mathcal{X}}$ of $\mathcal{X}$ such that ${\rm rank}(\widehat{\mathcal{X}}^{(i)})={\rm rank}(\widehat{\widetilde{\mathcal{X}}}^{(i)}),~i=1,2,\ldots,n_3$. It holds that
	\begin{equation} \label{eqle4}
		{\left\|\mathcal{X}^\dagger-\widetilde{\mathcal{X}}^\dagger\right\|}_2\leq\sqrt{2}{\left\|\mathcal{X}^\dagger\right\|}_2{\left\|\widetilde{\mathcal{X}}^\dagger\right\|}_2{\left\|\mathcal{X}-\widetilde{\mathcal{X}}\right\|}_2.
	\end{equation}
\end{lemma}
\begin{proof}
We notice that
	\begin{equation*}
		{\left\|\mathcal{X}^\dagger-\widetilde{\mathcal{X}}^\dagger\right\|}_2
=\max_{i=1,2,\ldots,n_3}{\left\|{{\left(\widehat{\mathcal{X}}^{(i)}\right)}^\dagger}-{{\left(\widehat{\widetilde{\mathcal{X}}}^{(i)}\right)}^\dagger}\right\|}_2.
	\end{equation*}
	From \cite[Theorem 4.1]{Wedin}, we have
	\begin{equation*}
		{\left\|{{\left(\widehat{\mathcal{X}}^{(i)}\right)}^\dagger}-{{\left(\widehat{\widetilde{\mathcal{X}}}^{(i)}\right)}^\dagger}\right\|}_2\leq\sqrt{2}{\left\|{\left(\widehat{\mathcal{X}}^{(i)}\right)}^\dagger\right\|}_2{\left\|{\left(\widehat{\widetilde{\mathcal{X}}}^{(i)}\right)}^\dagger\right\|}_2{\left\|\widehat{\mathcal{X}}^{(i)}-\widehat{\widetilde{\mathcal{X}}}^{(i)}\right\|}_2, i=1,2,\ldots,n_3.
	\end{equation*}

Therefore, it follows that
	\begin{align*}
		{\left\|\mathcal{X}^\dagger-\widetilde{\mathcal{X}}^\dagger\right\|}_2\!&=\!\max_{i=1,2,\ldots,n_3}{\left\|{{\left(\widehat{\mathcal{X}}^{(i)}\right)}^\dagger}-{{\left(\widehat{\widetilde{\mathcal{X}}}^{(i)}\right)}^\dagger}\right\|}_2  \\
\!&\leq\! \max_{i=1,2,\ldots,n_3}\left\{\sqrt{2}{\left\|{\left(\widehat{\mathcal{X}}^{(i)}\right)}^\dagger\right\|}_2{\left\|{\left(\widehat{\widetilde{\mathcal{X}}}^{(i)}\right)}^\dagger\right\|}_2{\left\|\widehat{\mathcal{X}}^{(i)}-\widehat{\widetilde{\mathcal{X}}}^{(i)}\right\|}_2\right\}\\
		\!&\leq\!
\sqrt{2}\max_{i=1,2,\ldots,n_3}\left\{{\left\|{\left(\widehat{\mathcal{X}}^{(i)}\right)}^\dagger\right\|}_2\right\} \max_{i=1,2,\ldots,n_3}\left\{{\left\|{\left(\widehat{\widetilde{\mathcal{X}}}^{(i)}\right)}^\dagger\right\|}_2\right\}
\max_{i=1,2,\ldots,n_3}\left\{{\left\|\widehat{\mathcal{X}}^{(i)}-\widehat{\widetilde{\mathcal{X}}}^{(i)}\right\|}_2\right\}\\
		\!&=\!\sqrt{2}{\left\|\mathcal{X}^\dagger\right\|}_2{\left\|\widetilde{\mathcal{X}}^\dagger\right\|}_2{\left\|\mathcal{X}-\widetilde{\mathcal{X}}\right\|}_2.
	\end{align*}
\end{proof}


We now proceed to prove the central theorem for this section; the proof consists of two stages.

{\tt (i)} To begin with, we derive an upper bound characterizing the discrepancy between $\mathcal{X}_r^T$ and our approximation $\widetilde{\mathcal{X}}^T$.
Denote by
\begin{align*}
	\mathcal{L}={\left(\mathcal{X}\ast\mathcal{X}^T\right)}^q\ast\mathcal{X},\quad{\rm and}\quad \mathcal{G}=	\mathcal{L}\ast	\mathcal{W},
\end{align*}
Algorithm \ref{Alg3} implies that $\mathcal{Q}_1\ast\mathcal{Q}_1^T$ corresponds to the orthogonal projection operator onto $\operatorname{span}\left\{\mathcal{G}\right\}$.

From \cite[Theorem 10.8]{Halko}, one obtains that
\begin{align*}	{\left\|(\mathcal{I}-\mathcal{Q}_1\ast\mathcal{Q}_1^T)\ast\mathcal{L}\right\|}_2\!&=\!\underset{i=1,2,\ldots,n_3}\max {\left\|\left(\mathbf{I}-\widehat{\mathcal{Q}}_1^{(i)}{\left(\widehat{\mathcal{Q}_1^T}\right)}^{(i)}\right)\widehat{\mathcal{L}}^{(i)}\right\|}_2 \\
\!&\equiv\!{\left\|\left(\mathbf{I}-\widehat{\mathcal{Q}}_1^{(l)}{\left(\widehat{\mathcal{Q}_1^T}\right)}^{(l)}\right)\widehat{\mathcal{L}}^{(l)}\right\|}_2 \\
\!&\leq\!\left(1+t\sqrt{\frac{3r}{p+1}}+u t\frac{{\rm e}\sqrt{r+p}}{p+1}\right)\widetilde{\sigma}_{r+1}^{(l)}+t\frac{{\rm e}\sqrt{r+p}}{p+1}{\left(\sum_{j>r}\left(\widetilde{\sigma}_j^{(l)}\right)^2\right)}^{1/2},
\end{align*}
where the failure probability is upper-bounded by $2t^{-p}+{\rm e}^{-u^2/2}$, in which $p\geq 4$. Observe that $\widetilde{\sigma}_1^{(l)}\geq \widetilde{\sigma}_2^{(l)}\geq\cdots\geq \widetilde{\sigma}_{\min\{n_1, m\}}^{(l)}\geq 0$ represent the singular values associated with $\widehat{\mathcal{L}}^{(l)}$.

For tensor $\mathcal{X}^T$, its economized singular value decomposition takes the form $\mathcal{X}^T=\mathcal{U}\ast\mathcal{S}\ast\mathcal{P}^T$. Thus, $\widehat{\mathcal{X}^T}^{(i)}=\widehat{\mathcal{U}}^{(i)}\widehat{\mathcal{S}}^{(i)}\widehat{\mathcal{P}^T}^{(i)}$, and
$$
\widehat{\mathcal{S}}^{(i)} = {\rm diag}\left(\sigma_1^{(i)}, \sigma_2^{(i)}, \ldots, \sigma_{\min\{n_1, m\}}^{(i)}\right) \in \mathbb{R}^{\min\{n_1, m\} \times \min\{n_1, m\}},
$$
with $\sigma_1^{(i)} \geq \sigma_2^{(i)} \geq \cdots \geq \sigma_{\min\{n_1, m\}}^{(i)} > 0$, $i=1,2,\ldots,n_3$. It holds that
$$
\mathcal{L}=\mathcal{P}\ast\mathcal{S}^{2q+1}\ast\mathcal{U}^T,
$$
and
$$
\widetilde{\sigma}_i^{(j)}={(\sigma_i^{(j)})}^{2q+1},~i=1,2,\ldots,\min\{n_1, m\},~j=1,2,\ldots,n_3.
$$
According to formula \eqref{eqle3}, one obtains
\begin{equation*}
	{\left\|\left(\mathcal{I}-\mathcal{Q}_1\ast\mathcal{Q}_1^T\right)\ast\mathcal{X}\right\|}_2\leq{\left\|\left(\mathcal{I}-\mathcal{Q}_1\ast\mathcal{Q}_1^T\right)\ast\mathcal{L}\right\|}_2^{1/{2q+1}}.
\end{equation*}
Consequently, one obtains
\begin{align}\label{equa3.3.12}
	&{\left\|\mathcal{X}^T \!-\!\mathcal{X}^T\ast\mathcal{Q}_1\ast\mathcal{Q}_1^T\right\|}_2
\!=\!	{\left\|\left(\mathcal{I}\!-\!\mathcal{Q}_1\ast\mathcal{Q}_1^T\right)\ast\mathcal{X}\right\|}_2 \nonumber\\
	\!&\leq\! {\left[\left(1+t\sqrt{\frac{3r}{p+1}}\!+\!ut\frac{{\rm e}\sqrt{r+p}}{p+1}\right){\left(\sigma_{r+1}^{(l)}\right)}^{2q+1}
		\!+\!t\frac{{\rm e}\sqrt{r+p}}{p+1}{\left(\sum_{j>r}{\left(\sigma_j^{(l)}\right)}^{2(2q+1)}\right)}^{\frac{1}{2}}\right]}^{\frac{1}{2q+1}},
\end{align}
where the failure probability is upper-bounded by $2t^{-p}+{\rm e}^{-u^2/2}$, with conditions $p\geq4$ and $\mu,~t\geq1$.

Moreover, we obtain from Algorithm \ref{Alg3} that
$$
\mathcal{X}^T\ast\mathcal{Q}_1\ast\mathcal{Q}_1^T=\mathcal{U}_2\ast\mathcal{S}_1\ast\mathcal{V}_1^T,
$$
where $\widetilde{\mathcal{X}}^T=\widetilde{\mathcal{U}}_r\ast\widetilde{\mathcal{S}}_r\ast\widetilde{\mathcal{V}}_r^T$ constitutes a rank-$r$ approximation for ${\mathcal{X}}^T$. Denote
$$
\sigma_{r+1}^{(j)}\left(\widehat{\mathcal{X}^T}^{(j)}\widehat{\mathcal{Q}_1}^{(j)}\widehat{\mathcal{Q}_1^T}^{(j)}\right),~j=1,2,\ldots,n_3,
$$
as the $(r+1)$-th dominant singular value belonging to the $j$-th frontal slice of
$$
{\rm fft}\left(\mathcal{X}^T\ast\mathcal{Q}_1\ast\mathcal{Q}_1^T,[~],3\right).
$$
As stated by interlacing theorem for singular values \cite{MatrixC}, we get
$$
\sigma_{r+1}^{(j)}\left(\widehat{\mathcal{X}^T}^{(j)}\widehat{\mathcal{Q}_1}^{(j)}\widehat{\mathcal{Q}_1^T}^{(j)}\right)\leq \sigma_{r+1}^{(j)},~j=1,2,\ldots,n_3,
$$
for which $\sigma_{r+1}^{(j)}$ denotes the $(r+1)$-th dominant singular value associated with the $j$-th frontal slice of $\widehat{\mathcal{X}^T}$.
Consequently,
\begin{align}\label{equa3.3.13}
	{\left\|\mathcal{U}_2\ast\mathcal{S}_1\ast\mathcal{V}_1^T-\widetilde{\mathcal{U}}_r\ast\widetilde{\mathcal{S}}_r\ast\widetilde{\mathcal{V}}_r^T\right\|}_2
\!&=\! \max_{j=1,2,\ldots,n_3}{\left\|\widehat{\mathcal{U}_2}^{(j)}\widehat{\mathcal{S}_1}^{(j)}\widehat{\mathcal{V}_1^T}^{(j)}-\widehat{\widetilde{\mathcal{U}}_r}^{(j)}\widehat{\widetilde{\mathcal{S}}_r}^{(j)}\widehat{\widetilde{\mathcal{V}}_r^T}^{(j)}\right\|}_2 \nonumber \\
\!&=\! \max_{j=1,2,\ldots,n_3} \sigma_{r+1}^{(j)}\left(\widehat{\mathcal{U}_2}^{(j)}\widehat{\mathcal{S}_1}^{(j)}\widehat{\mathcal{V}_1^T}^{(j)}\right) \nonumber \\
\!&\leq\! \max_{j=1,2,\ldots,n_3} \sigma_{r+1}^{(j)} \nonumber \\
\!&\equiv\! \sigma_{r+1}^{(l)}.
\end{align}
Combining \eqref{equa3.3.12} and \eqref{equa3.3.13} yields
\begin{align}\label{equa3.3.14}
\!&~\!	{\left\|\mathcal{X}^T\!-\!\widetilde{\mathcal{X}}^T\right\|}_2 \nonumber \\
\!&=\!{\left\|\mathcal{X}^T\!-\!\widetilde{\mathcal{U}}_r\ast\widetilde{\mathcal{S}}_r\ast\widetilde{\mathcal{V}}_r^T\right\|}_2
\leq {\left\|\mathcal{X}^T\!-\!\mathcal{X}^T\ast\mathcal{Q}_1\ast\mathcal{Q}_1^T\right\|}_2\!+\!{\left\|\mathcal{U}_2\ast\mathcal{S}_1\ast\mathcal{V}_1^T\!-\!\widetilde{\mathcal{U}}_r\ast\widetilde{\mathcal{S}}_r\ast\widetilde{\mathcal{V}}_r^T\right\|}_2 \nonumber \\
	\!&\leq\!{\left[\left(1+t\sqrt{\frac{3r}{p+1}}\!+\!u t\frac{{\rm e}\sqrt{r+p}}{p+1}\right){\left(\sigma_{r+1}^{(l)}\right)}^{2q+1}\!+\!t\frac{{\rm e}\sqrt{r+p}}{p+1}{\left(\sum_{j>r}{\left(\sigma_j^{(l)}\right)}^{2(2q+1)}\right)}^{\frac{1}{2}}\right]}^{\frac{1}{2q+1}} \nonumber \\
	\!&~+\!\sigma_{r+1}^{(l)},
\end{align}
where the failure probability is upper-bounded by $2t^{-p}+{\rm e}^{-u^2/2}$. Let
\begin{align}\label{equa3.3.15}
	\eta={\bigl[\bigl(1+{\rm e}\sqrt{\frac{3r}{p}}+{\rm e}^2\frac{\sqrt{2(r+p)}}{p}\bigr){\bigl(\sigma_{r+1}^{(l)}\bigr)}^{2q+1}+
		{\rm e}^2\frac{\sqrt{r+p}}{p}{\bigl(\sum_{j>r}{\bigl(\sigma_j^{(l)}\bigr)}^{2(2q+1)}\bigr)}^{\frac{1}{2}}\bigr]}^{\frac{1}{2q+1}},
\end{align}
suppose we set $t={\rm e},~u=\sqrt{2p}$, \eqref{equa3.3.14} reduces to
\begin{align}\label{equa3.3.16}
	{\left\|\mathcal{X}^T-\widetilde{\mathcal{X}}^T\right\|}_2\leq\eta+\sigma_{r+1}^{(l)}
=\left( \mu+1\right)\sigma_{r+1}^{(l)},
\end{align}
for which the failure probability is upper-bounded by $3{\rm e}^{-p}$, in which $\mu=\eta/{\sigma_{r+1}^{(l)}}$.
Moreover,
\begin{align}\label{equa3.3.17}
	{\left\|\mathcal{X}^T-\mathcal{X}_r^T\right\|}_2
=
\max_{j=1,2,\ldots,n_3}{\left\|{\left(\widehat{\mathcal{X}}^T\right)}^{(j)}-{\left(\widehat{\mathcal{X}}_r^T\right)}^{(j)}\right\|}_2
=
\max_{j=1,2,\ldots,n_3}\sigma_{r+1}^{(j)}=\sigma_{r+1}^{(l)}.
\end{align}
A combination of \eqref{equa3.3.16} and \eqref{equa3.3.17} yields
\begin{align}\label{equa3.3.18}
	{\left\|\mathcal{X}_r^T-\widetilde{\mathcal{X}}^T\right\|}_2\leq{\left\|\mathcal{X}^T-\mathcal{X}_r^T\right\|}_2+	{\left\|\mathcal{X}^T-\widetilde{\mathcal{X}}^T\right\|}_2
\leq \left( \mu+1\right)\sigma_{r+1}^{(l)}+\sigma_{r+1}^{(l)}
\leq \left(\mu+2\right)\sigma_{r+1}^{(l)}.
\end{align}

{\tt (ii)} In the second step, our analysis concerns the angle separating tensor subspaces $\widetilde{\mathbb{S}}={\rm span}\left\{{\left(\widetilde{\mathcal{X}}^T\right)}^\dagger\ast \overline{\mathcal{Y}}\right\}$ and $\mathbb{S}={\rm span}\left\{{\left(\mathcal{X}_r^T\right)}^\dagger\ast \overline{\mathcal{Y}}\right\}$.

Assume
\begin{equation}\label{eqq31}
{\rm rank}\left(\widehat{\left({\left(\mathcal{X}_r^T\right)}^\dagger\ast \overline{\mathcal{Y}}\right)}^{(i)}\right)={\rm rank}\left(\widehat{\left({\left(\widetilde{\mathcal{X}}^T\right)}^\dagger\ast \overline{\mathcal{Y}}\right)}^{(i)}\right),~i=1,2,\ldots,n_3.
\end{equation}
It holds that
\begin{align}\label{eqq32}
	{\left\|{(\mathcal{X}_r^T)}^\dagger\right\|}_2
\!&=\!
{\left\|\widehat{{\left(\mathcal{X}_r^T\right)}^\dagger}\right\|}_2
=\max_{i=1,2,\ldots,n_3}{ \left\|\widehat{{\left(\mathcal{X}_r^T\right)}^\dagger}^{(i)}\right\|}_2
=\max_{i=1,2,\ldots,n_3}{\left\| {\left(\widehat{\mathcal{X}_r^T}^{(i)}\right)} ^\dagger \right\|}_2 \nonumber\\
\!&=\!
\max_{i=1,2,\ldots,n_3}\frac{1}{\sigma_r^{(i)}}
\equiv\frac{1}{\sigma_r^{(s)}}.
\end{align}
Thus, it follows from \eqref{eqle4}, \eqref{equa3.3.18}and \eqref{eqq31} that
\begin{align}\label{eq32}
	{\left\|{\left(\mathcal{X}_r^T\right)}^\dagger\ast \overline{\mathcal{Y}}-{\left(\widetilde{\mathcal{X}}^T\right)}^\dagger\ast \overline{\mathcal{Y}}\right\|}_2\!&\leq\!{\left\|{\left(\mathcal{X}_r^T\right)}^\dagger-{\left(\widetilde{\mathcal{X}}^T\right)}^\dagger\right\|}_2{\left\| \overline{\mathcal{Y}}\right\|}_2  \nonumber\\
	\!&\leq\!\frac{{\left\|{\left(\mathcal{X}_r^T\right)}^\dagger-{\left(\widetilde{\mathcal{X}}^T\right)}^\dagger\right\|}_2}{{\left\|{\left(\mathcal{X}_r^T\right)}^\dagger\right\|}_2}{\left\|{\left(\mathcal{X}_r^T\right)}^\dagger\right\|}_2{\left\| \overline{\mathcal{Y}}\right\|}_2     \nonumber\\
	\!&\leq\!\sqrt{2}{\left\|{\left(\widetilde{\mathcal{X}}^T\right)}^\dagger\right\|}_2{\left\|\mathcal{X}_r^T-\widetilde{\mathcal{X}}^T\right\|}_2{\left\|{\left(\mathcal{X}_r^T\right)}^\dagger\right\|}_2{\left\| \overline{\mathcal{Y}}\right\|}_2    \nonumber\\
	\!&\leq\!\sqrt{2}\left(\mu+2\right){\left\|{\left(\widetilde{\mathcal{X}}^T\right)}^\dagger\right\|}_2\frac{\sigma_{r+1}^{(l)}}{\sigma_r^{(s)}},
\end{align}
where we used the fact that ${\left\|\overline{\mathcal{Y}}\right\|}_2=1$.

Moreover,
\begin{align}\label{equa3.3.21}
{\left\|{\left[{\left(\mathcal{X}_r^T\right)}^\dagger\ast \overline{\mathcal{Y}}\right]}^\dagger\right\|}_2\!&=\!
\max_{i=1,2,\ldots,n_3}{\left\|\widehat{{\left[{\left(\mathcal{X}_r^T\right)}^\dagger\ast \overline{\mathcal{Y}}\right]}^\dagger}^{(i)}\right\|}_2 \nonumber
= \max_{i=1,2,\ldots,n_3}{\left\|\widehat{{\left[{\left(\mathcal{X}_r^T\right)}^\dagger\ast \overline{\mathcal{Y}}\right
]}^{(i)}}^\dagger\right\|}_2 \nonumber \\
\!&=\!
\max_{i=1,2,\ldots,n_3} \frac{1}{\sigma_{\min}^{(i)}\left( {\widehat{{\left(\mathcal{X}_r^T\right)}^\dagger}}^{(i)} {\widehat{\overline{\mathcal{Y}}}}^{(i)}\right)} \nonumber
\leq
\max_{i=1,2,\ldots,n_3} \frac{1}{\sigma_{\min}^{(i)}\left( {\widehat{{\left(\mathcal{X}_r^T\right)}^\dagger}}^{(i)}\right)} \nonumber \\
\!&=\!
\max_{i=1,2,\ldots,n_3} \sigma_{\max}^{(i)}\left({\widehat{{\left(\mathcal{X}_r^T\right)}}}^{(i)}\right) \nonumber \\
\!&\equiv\! \sigma_1^{(w)}.
\end{align}

Denote $\mathcal{P}_{{\left(\mathcal{X}_r^T\right)}^\dagger\ast \overline{\mathcal{Y}}}$ and $\mathcal{P}_{{\left(\widetilde{\mathcal{X}}^T\right)}^\dagger\ast \overline{\mathcal{Y}}}$ as the orthogonal projectors corresponding to subspaces ${\rm span}\left\{{\left(\mathcal{X}_r^T\right)}^\dagger\ast \overline{\mathcal{Y}}\right\}$ and ${\rm span}\left\{{\left(\widetilde{\mathcal{X}}^T\right)}^\dagger\ast \overline{\mathcal{Y}}\right\}$, respectively. Thus, it follows from {\cite[125--127]{Sun}, \eqref{eq32} and \eqref{equa3.3.21} that
\begin{align*}
	\sin\angle\left(\mathbb{S},\widetilde{\mathbb{S}}\right)\!&=\!{\left\|\mathcal{P}_{{\left(\mathcal{X}_r^T\right)}^\dagger\ast \overline{\mathcal{Y}}}-\mathcal{P}_{{\left(\widetilde{\mathcal{X}}^T\right)}^\dagger\ast \overline{\mathcal{Y}}}\right\|}_2
= \max_{i=1,2,\ldots,n_3} {\left\|\widehat{\mathcal{P}}_{{\left(\mathcal{X}_r^T\right)}^\dagger\ast \overline{\mathcal{Y}}}^{(i)}-\widehat{\mathcal{P}}_{{\left(\widetilde{\mathcal{X}}^T\right)}^\dagger\ast \overline{\mathcal{Y}}}^{(i)}\right\|}_2 \\
\!&\leq\! \max_{i=1,2,\ldots,n_3}\left\{
    \min\left\{
        \left\| \left[ \widehat{ \left( \left( \mathcal{X}_r^T \right)^\dagger \ast \overline{\mathcal{Y}} \right) }^{(i)} \right]^\dagger \right\|_2,
        \left\| \left[ \widehat{ \left( \left( \widetilde{\mathcal{X}}^T \right)^\dagger \ast \overline{\mathcal{Y}} \right) }^{(i)} \right]^\dagger \right\|_2
    \right\}\times\right.\\
    \!&~\!\left.\left\| \left( \widehat{ \left( \mathcal{X}_r^T \right)^\dagger \ast \overline{\mathcal{Y}} } - \widehat{ \left( \widetilde{\mathcal{X}}^T \right)^\dagger \ast \overline{\mathcal{Y}} } \right)^{(i)} \right\|_2
\right\}\\
\!&\leq\! \max_{i=1,2,\ldots,n_3} \left\{
    \left\| \left[ \widehat{ \left( \left( \mathcal{X}_r^T \right)^\dagger \ast \overline{\mathcal{Y}} \right) }^{(i)} \right]^\dagger \right\|_2
    \left\| \left( \widehat{ \left( \mathcal{X}_r^T \right)^\dagger \ast \overline{\mathcal{Y}} } - \widehat{ \left( \widetilde{\mathcal{X}}^T \right)^\dagger \ast \overline{\mathcal{Y}} } \right)^{(i)} \right\|_2
\right\}\\
\!&\leq\!\sqrt{2}\left(\mu+2\right){\left\|{\Big(\widetilde{\mathcal{X}}^T\Big)}^\dagger\right\|}_2 \sigma_1^{(w)}\sigma_{r+1}^{(l)}/{\sigma_r^{(s)}}.
\end{align*}

In conclusion, let $\mathcal{F}_Q$ and $\widetilde{\mathcal{F}}_Q$ stand for the $Q$-factors yielded by economical $QR$ decompositions for ${\left(\mathcal{X}_r^T\right)}^\dagger\ast \overline{\mathcal{Y}}$ and ${\left(\widetilde{\mathcal{X}}^T\right)}^\dagger\ast \overline{\mathcal{Y}}$, correspondingly.
Theorem \ref{theo3.3} characterizes the distance separating the truncated solution $\mathcal{F}_Q$ from the approximate solution $\widetilde{\mathcal{F}}_Q$.
\begin{theorem}\label{theo3.3}
	Define two tensor subspaces as $\mathbb{S}={\rm span}\left\{{\left(\mathcal{X}_r^T\right)}^\dagger\ast \overline{\mathcal{Y}}\right\}={\rm span}\left\{ \overline{\mathcal{F}}_Q\right\}$
and
$\widetilde{\mathbb{S}}={\rm span}\left\{{\left(\widetilde{\mathcal{X}}^T\right)}^\dagger\ast \overline{\mathcal{Y}}\right\}={ \rm span}\left\{\widetilde{\mathcal{F}}_Q\right\}$.
Under the rank-preserving assumption
 $$
 {\rm rank}\left(\widehat{\left({\left(\mathcal{X}_r^T\right)}^\dagger\ast \overline{\mathcal{Y}}\right)}^{(i)}\right)={\rm rank}\left(\widehat{\left({\left(\widetilde{\mathcal{X}}^T\right)}^\dagger\ast \overline{\mathcal{Y}}\right)}^{(i)}\right),~i=1,2,\ldots,n_3,
 $$
the following subspace-angle bound holds:
	\begin{equation*}
\sin\angle\left(\mathbb{S},\widetilde{\mathbb{S}}\right)\leq\sqrt{2}\left(\mu+2\right){\left\|{(\widetilde{\mathcal{X}}^T)}^\dagger\right\|}_2\sigma_1^{(w)}/{\sigma_r^{(s)}}\cdot\sigma_{r+1}^{(l)},
	\end{equation*}
	The above estimate holds with failure probability no larger than $3 {\rm e}^{-p}$ with $p\ge4$, where $r$ is the tubal rank.
$\sigma_1^{(w)},\sigma_{r+1}^{(l)}$ and $\sigma_r^{(s)}$
are specified by \eqref{equa3.3.21}, \eqref{equa3.3.13} and \eqref{eqq32}, correspondingly.
\end{theorem}

\begin{remark}
{\it Results derived from Theorem \ref{theo3.3} reveal how the separation measure between the two tensor subspaces depends on several key contributing quantities. These terms include $\sigma_{r+1}^{(l)}$, which represents the $(r+1)$-th largest singular value computed from the $l$-th frontal slice of $\widehat{\mathcal{X}^T}$; $\sigma_1^{(w)}$, corresponding to the maximum singular value associated with the $w$-th frontal slice of $\widehat{{\mathcal{X}_r^T}}$; $\sigma_r^{(s)}$, standing for the $r$-th dominant singular value obtained from the $s$-th frontal slice of $\widehat{{\mathcal{X}_r^T}}$, together with the matrix norm of the MP-inverse ${(\widetilde{\mathcal{X}}^T)}^\dagger$.}
\end{remark}

\section{A sparse multilinear discriminant analysis method for tensor feature selection}\label{sec4}

Traditional multilinear discriminant analysis methods often face with limitations including the curse of dimensionality and interference from redundant features, which may lead to suboptimal performances. Therefore, it is desirable to seek new technologies capable of eliminating redundant information or filtering out irrelevant features during the feature extraction step.

Within this section, we develop a sparse multilinear discriminant analysis approach based on tensor $t$-product for feature selection and dimensionality reduction. Compared with prior approaches, the presented method addresses the tensor trace-ratio problem in a direct manner. We first construct a sparse model utilizing the closed-form solution associated with the tensor trace-ratio problem together with $\ell_{2,0}$-norm regularization, and then consider how to solve it efficiently. Our strategy not only reduces computational complexities, but also enhances the model's generalization ability.

\subsection{A sparse multilinear discriminant analysis model}
In the present subsection, a sparse multilinear discriminant analysis model (SMDA) is constructed relying on the tensor least-squares problem \eqref{equa3.2.1} as well as the tensor $\ell_{2,0}$-norm. To this end, we first provide the definition associated with the tensor $\ell_{2,0}$-norm.

\begin{definition} \label{defm0}
	\cite{Nie2020Subspace} The $\ell_{2,0}$-norm of a matrix $\mathbf{A}\in \mathbb{R}^{n\times m}$ is defined as
	\begin{equation*}
		{\|\mathbf{A}\|}_{2,0} = \sum_{i=1}^{n} {\left \| \sum_{j=1}^{m} {\left(\mathbf{A}(i,j)\right)}^{2} \right \|}_{0},
	\end{equation*}
where for a scalar $a$, $\|a\|_0 = 1$ if $a \neq 0$, $\|a\|_0 = 0$ if $a = 0$.
\end{definition}

By analogy, the tensor $\ell_{2,0}$-norm can be formulated as shown below.
\begin{definition} \label{def0norm}
    For a tensor $\mathcal{M}\in\mathbb{R}^{n_1\times n_2 \times n_3}$, its $\ell_{2,0}$-norm reads
	\begin{equation*}
	{|||\mathcal{M}|||}_{2,0}=\sum_{i=1}^{n_1}{\rm sn}({\left \|\mathcal{M}(i,:,:)\right\|}_{2,0}),
	\end{equation*}
in which $\text{sn}(\cdot)$ denotes an indicator function, with ${\rm sn}({\left \|\mathcal{M}(i,:,:)\right\|}_{2,0})$=1 if ${\left \|\mathcal{M}(i,:,:)\right\|}_{2,0}\ne0$, and ${\rm sn}({\left \|\mathcal{M}(i,:,:)\right\|}_{2,0})$=0 otherwise.
\end{definition}
We are ready to present a sparse multilinear discriminant analysis model for multilinear discriminant analysis feature selection
\begin{align} \label{eqtr35}
  \max_{\mathcal{V} \in \mathbb{R}^{n_1\times d\times n_3}} \frac{{\rm tr}\left(\mathcal{V}^T\ast\mathcal{S}_b\ast\mathcal{V}\right)}{{\rm tr}\left(\mathcal{V}^T\ast\mathcal{S}_t\ast\mathcal{V}\right)} \\
 \operatorname{s.t.} \mathcal{V}^T\ast\mathcal{V}=\mathcal{I}_d,~{|||\mathcal{V}|||}_{2,0}=k. \nonumber
\end{align}
Motivated by \eqref{equa3.2.1}, we reformulate \eqref{eqtr35} as
\begin{align}\label{equa4.1}
	\min_{\mathcal{F}\in\mathbb{R}^{n_1\times d\times n_3}} {\left\|{\mathcal{X}}^T\ast\mathcal{F}-\overline{\mathcal{Y}}\right\|}_F \\
	\operatorname{s.t.} \mathcal{F}^T\ast\mathcal{F}=\mathcal{I}_d,~{|||\mathcal{F}|||}_{2,0}=k, \nonumber
\end{align}
in which $d=c-1$ denotes the reduced dimension, while $k$ counts the selected features, i.e., the nonzero horizontal slices in $\mathcal{F}$, satisfying $d\leq k\leq n_1$.

\subsection{A sparse multilinear discriminant analysis method}
Within this subsection, we investigate approaches for solving the optimization model \eqref{equa4.1}. The subsequent equivalent-form theorem is required.
\begin{theorem}\label{Theo4.1}
Let $\widetilde{\mathcal{A}}=-{\mathcal{X}}\ast{\mathcal{X}}^T+{\mathcal{X}}\ast\overline{\mathcal{Y}}\ast\mathcal{F}^T+\mathcal{F}\ast \overline{\mathcal{Y}}^T\ast{\mathcal{X}}^T$, then \eqref{equa4.1}
is equivalent to
\begin{align} \label{eqstr}
\max_{\mathcal{F}\in\mathbb{R}^{n_1\times d\times n_3}} {\rm tr}\left(\mathcal{F}^T\ast\widetilde{\mathcal{A}}\ast\mathcal{F}\right) \nonumber\\
\operatorname{s.t.} \mathcal{F}^T\ast\mathcal{F}=\mathcal{I}_d,~{|||\mathcal{F}|||}_{2,0}=k.
\end{align}
\end{theorem}
\begin{proof}
Notice that
\begin{align*}
{\left\|{\mathcal{X}}^T\ast\mathcal{F}-\overline{\mathcal{Y}}\right\|}_F
\!&=\!
{\rm tr}\left({\left({\mathcal{X}}^T\ast\mathcal{F}-\overline{\mathcal{Y}}\right)}^T\ast\left({\mathcal{X}}^T\ast\mathcal{F}-\overline{\mathcal{Y}}\right)\right) \\
\!&=\!
{\rm tr}\left(\mathcal{F}^T\ast{\mathcal{X}}\ast{\mathcal{X}}^T\ast\mathcal{F}-\mathcal{F}^T\ast{\mathcal{X}}\ast\overline{\mathcal{Y}}-\overline{\mathcal{Y}}^T\ast{\mathcal{X}}^T\ast\mathcal{F}+\overline{\mathcal{Y}}^T\ast \overline{\mathcal{Y}}\right).
\end{align*}
As ${\rm tr}(\overline{\mathcal{Y}}^T\ast\overline{\mathcal{Y}})$ is a constant, let $\widetilde{\mathcal{A}}=-{\mathcal{X}}\ast{\mathcal{X}}^T+{\mathcal{X}}\ast\overline{\mathcal{Y}}\ast\mathcal{F}^T+\mathcal{F}\ast \overline{\mathcal{Y}}^T\ast{\mathcal{X}}^T$, then we obtain
\begin{align*}
{\rm tr}\left(\mathcal{F}^T\ast{\mathcal{X}}\ast{\mathcal{X}}^T\ast\mathcal{F}-\mathcal{F}^T\ast{\mathcal{X}}\ast\overline{\mathcal{Y}}-\overline{\mathcal{Y}}^T\ast{\mathcal{X}}^T\ast\mathcal{F}\right) ={\rm tr}\left(\mathcal{F}^T\ast\widetilde{\mathcal{A}}\ast\mathcal{F}\right),
\end{align*}
which completes the proof.
\end{proof}
However, solving \eqref{eqstr} directly is NP-hard, as the two nonlinear constraints are difficult to satisfy. Motivated by the sparsity strategy proposed in \cite{Nie2020Subspace}, we rewrite \eqref{eqstr} as the following optimization problem
	\begin{align} \label{eqstrtwo}
		\max_{\mathcal{F}\in\mathbb{R}^{n_1\times d\times n_3}} {\rm tr}\left(\mathcal{F}^T\ast\mathcal{A}\ast\mathcal{F}\right) \\
		\operatorname{s.t.} \mathcal{F}^T\ast\mathcal{F}=\mathcal{I}_d,~{|||\mathcal{F}|||}_{2,0}=k, \nonumber
	\end{align}
where
$$
\mathcal{A}=-{\mathcal{X}}\ast{\mathcal{X}}^T+{\mathcal{X}}\ast\overline{\mathcal{Y}}\ast\mathcal{F}^T+\mathcal{F}\ast \overline{\mathcal{Y}}^T\ast{\mathcal{X}}^T+ \gamma \mathcal{I}
$$
is positive semi-definite as a tensor when $\gamma>0$ is large enough.
We next consider how to choose $\gamma$ in practice. Recall that $\mathcal{F}^T \ast \mathcal{F}=\mathcal{I}$ and $\overline{\mathcal{Y}}^T \ast \overline{\mathcal{Y}}=\mathcal{I}$. From
  $$
  \widehat{\mathcal{A}}^{(i)}=
  -\widehat{{\mathcal{X}}}^{(i)} \widehat{{\mathcal{X}}^T}^{(i)} +\widehat{{\mathcal{X}}}^{(i)} \widehat{\overline{\mathcal{Y}}}^{(i)} \widehat{\mathcal{F}^T}^{(i)} +\widehat{\mathcal{F}}^{(i)} {\left(\widehat{\overline{\mathcal{Y}}^T}\right)}^{(i)} \widehat{{\mathcal{X}}^T}^{(i)} + \gamma \widehat{\mathcal{I}}^{(i)},\quad i=1,2,\ldots,n_3,
  $$
we obtain 
\begin{align*}
\!&~\!\left\| -\widehat{{\mathcal{X}}}^{(i)} \widehat{{\mathcal{X}}^T}^{(i)} +\widehat{{\mathcal{X}}}^{(i)} \widehat{\overline{\mathcal{Y}}}^{(i)} \widehat{\mathcal{F}^T}^{(i)} +\widehat{\mathcal{F}}^{(i)} {\left(\widehat{\overline{\mathcal{Y}}^T}\right)}^{(i)}
\widehat{{\mathcal{X}}^T}^{(i)} \right\|_2 \\
\!&\leq\!
\left\| \widehat{{\mathcal{X}}}^{(i)} \widehat{{\mathcal{X}}^T}^{(i)} \right\|_2 + 2\left\| \widehat{{\mathcal{X}}}^{(i)} \widehat{\overline{\mathcal{Y}}}^{(i)} \widehat{\mathcal{F}^T}^{(i)} \right\|_2  \\
\!&=\! \left\|\widehat{{\mathcal{X}}}^{(i)} \right\|_2^2 + 2\left\| \widehat{{\mathcal{X}}}^{(i)} \right\|_2 \\
\!&\leq\! \left\|\widehat{{\mathcal{X}}}^{(i)} \right\|_1^2 + 2\left\| \widehat{{\mathcal{X}}}^{(i)} \right\|_1,~ i=1,2,\ldots,n_3.
\end{align*}
Thus, we can use
 \begin{equation}\label{eqq39}
   \gamma={\left\{\max_{j=1,2,\ldots,m}{|||\widehat{{\mathcal{X}}}(:,j,:)|||}_1\right\}}^2+2\left\{\max_{j=1,2,\ldots,m}{|||\widehat{{\mathcal{X}}}(:,j,:)|||}_1\right\},
\end{equation}
where $|||\cdot|||_1$ is the $\ell_1$-norm of a tensor \cite{t-svd}.


In summary, \eqref{eqstrtwo} is mathematically equivalent to
\begin{equation}\label{eqs40}
  \max_{\mathcal{F}\in\mathbb{K}} \quad \operatorname{tr}\bigl(\mathcal{F}^T \ast \mathcal{T} \ast \mathcal{F}\bigr),
\end{equation}
in which
\begin{equation}\label{eqq40}
\mathcal{T}=\mathcal{A}\ast \mathcal{F} \ast {\left( \mathcal{F}^T \ast \mathcal{A} \ast \mathcal{F}\right)}^\dagger \ast \mathcal{F}^T \ast \mathcal{A},
\end{equation}
and
$$
\mathbb{K}=\{ \mathcal{F}\in \mathbb{R}^{n_1 \times d \times n_3} \mid \mathcal{F}^T \ast \mathcal{F}=\mathcal{I}_d,~{|||\mathcal{F}|||}_{2,0}=k\}.
$$
We point out that the final model \eqref{eqs40} is essentially different from the matrix model proposed in \cite{Nie2020Subspace}, and the former is not a trivial extension of the latter one.
Instead, we reformulate \eqref{equa4.1} as a tensor trace-optimization problem with constraint rather than a trace-ratio or ratio-trace formulation.

Inspired by \cite[Algorithm 3]{Nie2020Subspace}, we are in a position to consider how to construct an index set $\mathbf{J}$ for feature selection based on \eqref{eqs40}.
The difference between our work and that in \cite{Nie2020Subspace} is that, for a given index set $\mathbf{J}$, we need to define horizontal slice extraction tensors used for feature selection. We first need the following definition.

\begin{definition} \label{des18}
	Define the horizontal slice extraction tensor $\mathcal{Q}\in \mathbb{R}^{n_1\times k \times n_3}$ by
	\[ \mathcal{Q}(i,j,v)=
	\begin{cases}
		1, & \text{if } i = \mathbf{J} (j) \text{ and } v=1,\\
		0, & \text{otherwise},
	\end{cases}\]
	where $\mathbf{J}$ is a given index set with elements ranging from 1 to $n_1$.
\end{definition}

We first analyze some properties of $\mathcal{T}$, which is used to seek the index set $\mathbf{J}$.
\begin{theorem} \label{Theo8}
 The tensor $\mathcal{T}$ possesses the property of positive semi-definiteness, while its tensor-rank satisfies ${\rm rank_t}\left(\mathcal{T}\right)\leq d $.
\end{theorem}
\begin{proof}
  It can be noted that $\mathcal{A}$ has  the property of symmetric and positive semi-definiteness. Set
  $$
  \widehat{\mathcal{M}}^{(i)} = \widehat{\mathcal{A}}^{(i)} \widehat{\mathcal{F}}^{(i)} {\left(\widehat{\left( \mathcal{F}^T \ast \mathcal{A} \ast \mathcal{F} \right)^{\dagger}}\right)}^{(i)}\widehat{\mathcal{F}^T}^{(i)}{\left(\widehat{\mathcal{A}}^{(i)}\right)}^{1/2},\quad i=1,2,\ldots,n_3.
  $$
  Accordingly, one obtains
  \begin{align*}
    \widehat{\mathcal{M}}^{(i)} {\left(\widehat{\mathcal{M}}^{(i)}\right)}^H
    \!&=\!\widehat{\mathcal{A}}^{(i)} \widehat{\mathcal{F}}^{(i)} {\left(\widehat{\left( \mathcal{F}^T \ast \mathcal{A} \ast \mathcal{F} \right)^{\dagger}}\right)}^{(i)}\widehat{\mathcal{F}^T}^{(i)}\widehat{\mathcal{A}}^{(i)} \widehat{\mathcal{F}}^{(i)} {\left(\widehat{\left( \mathcal{F}^T \ast \mathcal{A} \ast \mathcal{F} \right)^{\dagger}}\right)}^{(i)}\widehat{\mathcal{F}^T}^{(i)}\widehat{\mathcal{A}}^{(i)} \\
    \!&=\!\widehat{\mathcal{T}}^{(i)},\quad i=1,2,\ldots,n_3.
  \end{align*}
  This demonstrates that $\mathcal{T}$ satisfies the positive semi-definite property. Moreover, ${\rm rank_t}\left(\mathcal{T}\right)\leq d$ can be derived by using ${\rm rank}\big(\widehat{\mathcal{T}}^{(i)}\big) \leq d,~i=1,2,\ldots,n_3$.
\end{proof}
In view of Theorem \ref{Theo8}, it suffices to consider the two cases as follows for constructing the index set $\mathbf{J}$.

{\tt (i)} The case of $k=d$. From \eqref{eqs40}, $\mathcal{F}$ can be rewritten as $\mathcal{F}=\mathcal{Q}\ast \mathcal{M}$, with $\mathcal{M}\in \mathbb{R}^{d \times d \times n_3}$ being an orthogonal tensor, and $\mathcal{Q}\in \mathbb{R}^{n_1 \times d \times n_3}$. Let $\widetilde{\mathcal{T}}=\mathcal{Q}^T\ast \mathcal{T} \ast \mathcal{Q}$, then \eqref{eqs40} can be reformulated as
\begin{equation}\label{eqr41}
  \max_{\mathcal{M}\in \mathbb{R}^{d \times d \times n_3},~\mathcal{M}\ast \mathcal{M}^T=\mathcal{M}^T \ast \mathcal{M}=\mathcal{I}} \quad \operatorname{tr}\bigl(\mathcal{M}^T \ast \widetilde{\mathcal{T}} \ast \mathcal{M}\bigr).
\end{equation}

From the definition of the tensor trace, we have that
\begin{equation*}
  {\rm tr}\bigl(\mathcal{M}^T \ast \widetilde{\mathcal{T}} \ast \mathcal{M}\bigr)
  ={\rm tr}\bigl(\widetilde{\mathcal{T}} \ast \mathcal{M} \ast \mathcal{M}^T\bigr)
  ={\rm tr}\bigl(\widetilde{\mathcal{T}}\bigr)
  ={\rm tr}\bigl(\widetilde{\mathcal{T}}^{(1)}\bigr).
\end{equation*}
Hence, \eqref{eqr41} is equivalent to
\begin{equation}\label{eqr42}
  \max_{\mathcal{Q}\in \mathbb{R}^{n_1\times d \times n_3},~\mathcal{Q}^T \ast \mathcal{Q}=\mathcal{I}} \quad {\rm tr}\bigl(\widetilde{\mathcal{T}}^{(1)}\bigr)=\max_{\mathcal{Q}^{(1)}\in \mathbb{R}^{n_1\times d \times 1}} \quad {\rm tr}\left({\bigl(\mathcal{Q}^{(1)}\bigr)}^T \mathcal{T}^{(1)} \mathcal{Q}^{(1)}\right).
\end{equation}
To maximize the optimization problem \eqref{eqr42}, one may select the rows together with columns associated with the top-$d$ largest diagonal entries from $\mathcal{T}^{(1)}$.
Thus, the index set $\mathbf{J}$ for constructing $\mathcal{Q}$ is formed by the corresponding indices of the first $d$ largest diagonal elements of $\mathcal{T}^{(1)}$.

{\tt (ii)} The case of $k\neq d$. Similarly,
it follows from \eqref{eqs40} that $\mathcal{F}$ can be rewritten as $\mathcal{F}=\mathcal{Q}\ast \mathcal{M}$, with $\mathcal{M}\in \mathbb{R}^{k \times d \times n_3}$ being an orthogonal tensor, and $\mathcal{Q}\in \mathbb{R}^{n_1 \times k \times n_3}$. Let $\widetilde{\mathcal{T}}=\mathcal{Q}^T\ast \mathcal{T} \ast \mathcal{Q}$, then \eqref{eqs40} can be reformulated as
\begin{equation}\label{eqr441}
  \max_{\mathcal{M}\in \mathbb{R}^{k \times d \times n_3},~\mathcal{M}^T \ast \mathcal{M}=\mathcal{I}} \quad \operatorname{tr}\bigl(\mathcal{M}^T \ast \widetilde{\mathcal{T}} \ast \mathcal{M}\bigr).
\end{equation}

From the definition of the tensor trace, it follows that
\begin{align*}
  {\rm tr}\bigl(\mathcal{M}^T \ast \widetilde{\mathcal{T}} \ast \mathcal{M}\bigr)
  &=
  {\rm tr}\bigl({(\mathcal{M}^T \ast \widetilde{\mathcal{T}} \ast \mathcal{M})}^{(1)}\bigr)
  ={\rm tr}\bigl(
  \frac{1}{n_3} \sum_{i=1}^{n_3}
  {(\widehat{\mathcal{M}}^T)}^{(i)}{\widehat{\widetilde{\mathcal{T}}}}^{(i)}{\widehat{\mathcal{M}}}^{(i)}
  \bigr)\\
  &=\frac{1}{n_3} \sum_{i=1}^{n_3}{\rm tr}
  \bigl(
  {(\widehat{\mathcal{M}}^{(i)})}^T{\widehat{\widetilde{\mathcal{T}}}}^{(i)}{\widehat{\mathcal{M}}}^{(i)}
  \bigr).
\end{align*}

From the property that ${(\widehat{\mathcal{M}}^{(i)})}^T \widehat{\mathcal{M}}^{(i)}=\mathbf{I}_d$, holding for every index $i=1,2,\ldots,n_3$, \eqref{eqr441} admits an alternative formulation given by

\begin{equation*}
\max_{\mathcal{Q}\in \mathbb{R}^{n_1\times k \times n_3},~\mathcal{Q}^T \ast \mathcal{Q}=\mathcal{I}}
\frac{1}{n_3}\sum_{i=1}^{n_3}
{\rm tr}\bigl({\widehat{\widetilde{\mathcal{T}}}}^{(i)}\bigr)
=\max_{\mathcal{Q}\in \mathbb{R}^{n_1\times k \times n_3},~\mathcal{Q}^T \ast \mathcal{Q}=\mathcal{I}}
{\rm tr}\bigl(\widetilde{\mathcal{T}}^{(1)}\bigr).
\end{equation*}

As a result, \eqref{eqs40} can be rewriten as
\begin{equation}\label{eqrr43}
  \max_{\mathcal{Q}\in \mathbb{R}^{n_1\times k \times n_3},~\mathcal{Q}^T \ast \mathcal{Q}=\mathcal{I}} \quad {\rm tr}\bigl(\widetilde{\mathcal{T}}^{(1)}\bigr)=\max_{\mathcal{Q}^{(1)}\in \mathbb{R}^{n_1\times k \times 1}} \quad {\rm tr}\left({\bigl(\mathcal{Q}^{(1)}\bigr)}^T \mathcal{T}^{(1)} \mathcal{Q}^{(1)}\right).
\end{equation}
According to \cite[Corollary  4.3.18]{Kyfan}, the optimum value is attained whenever the index set $\mathbf{J}$ associated with tensor $\mathcal{Q}$ is assembled from indices corresponding to the top-$k$ diagonal entries of $\mathcal{T}^{(1)}$. To sum up, the index set $\mathbf{J}$ is constructed from the top-$k$ diagonal entries belonging to $\mathcal{T}^{(1)}$ for both considered scenarios.

Thus, for the given tensor $\mathcal{Q}$, the optimization problem \eqref{eqstrtwo} reduces to
 \begin{equation}\label{eqr38}
  \max_{\mathcal{M}\in \mathbb{R}^{k \times d \times n_3},~\mathcal{M}^T \ast \mathcal{M}=\mathcal{I}_d} \quad \operatorname{tr}\bigl(\mathcal{M}^T \ast \mathcal{Q}^T\ast \mathcal{A} \ast \mathcal{Q} \ast \mathcal{M} \bigr),
\end{equation}
which attains the maximum value as the lateral slices belonging to $\mathcal{M}$  coincide with the eigenslices associated with the leading-$k$ dominant eigentubes computed from $\mathcal{Q}^T\ast \mathcal{A} \ast \mathcal{Q}$.

We proceed to introduce Algorithm \ref{Alg5} for the sparse multilinear discriminant analysis model \eqref{equa4.1}, which is implemented in the Fourier domain.
For this algorithm, $d$ is assigned as $c-1$. Furthermore, it should be emphasized that explicitly constructing or saving tensor $\mathcal{A}$ is unnecessary, and only the calculation of $\mathcal{A} \ast\mathcal{F}$ is needed; see Step 2 of Algorithm \ref{Alg5}. Moreover, suppose the t-SVD decomposition satisfies $ \mathcal{F}^T \ast (\mathcal{A} \ast\mathcal{F})=\mathcal{U}\ast\mathcal{S}\ast\mathcal{U}^T$. One then obtains
\begin{equation*}
  \mathcal{T}=\mathcal{A}\ast \mathcal{F} \ast \mathcal{U}\ast \mathcal{S}^\dagger\ast \mathcal{U}^T \ast \mathcal{F}^T \ast \mathcal{A},
\end{equation*}
refer to Step 4 of Algorithm \ref{Alg5}.

\begin{algorithm} 
	\caption{A sparse multilinear discriminant analysis method with tensor $\ell_{2,0}$-norm regularization} \label{Alg5}
\begin{algorithmic}
\item  {{\bf Input:} Given the tensor dataset $\mathcal{X}\in \mathbb{R}^{n_1\times m\times n_3}$, $\overline{\mathcal{Y}}\in \mathbb{R}^{m\times (c-1)\times n_3}$, the initial guess $\mathcal{F}_0\in \mathbb{R}^{n_1\times d \times n_3}$, the dimensionality reduction number $d=c-1$, the feature selection number $k$, the regularization coefficient $\gamma$ defined in \eqref{eqq39}, the upper bound for iteration count maxiter, along with the convergence tolerance {\it tol}; \\
    {\bf Output:} $\mathcal{F}\in \mathbb{R}^{n_1\times d\times n_3}$;} \\
		\item  {{\bf Repeat $\colon$}}
	    \item  {{\bf 1.} Let $iter=0,1,\dots,maxiter-1$;}\\
		\item  {{\bf 2.} Compute $\mathcal{A}\ast\mathcal{F}_{iter}=-\mathcal{X}\ast \left(\mathcal{X}^T\ast\mathcal{F}_{iter}\right)+\mathcal{X}\ast\overline{\mathcal{Y}}+\mathcal{F}_{iter}\ast \overline{\mathcal{Y}}^T\ast\left(\mathcal{X}^T \ast \mathcal{F}_{iter}\right)+\gamma \mathcal{F}_{iter}$ ;}\\
		\item  {{\bf 3.} Perform $t$-SVD: $\mathcal{F}_{iter}^T\ast(\mathcal{A}\ast\mathcal{F}_{iter})=\mathcal{U}_{iter}\ast \mathcal{S}_{iter}\ast \mathcal{U}_{iter}^T$;} \\
		\item  {{\bf 4.} Compute $\mathcal{T}_{iter}=\left(\mathcal{A}\ast \mathcal{F}_{iter}\right) \ast \mathcal{U}_{iter}\ast \mathcal{S}_{iter}^\dagger\ast {\mathcal{U}_{iter}}^T \ast {\left(\mathcal{A}\ast \mathcal{F}_{iter}\right)}^T $;} \\
		\item  {{\bf 5.} Sort the diagonal elements of $\mathcal{T}_{iter}^{(1)}$, and construct the indices set $\mathbf{J}_{iter}$ by using the largest $k$ diagonal elements of $\mathcal{T}^{(1)}_{iter}$;} \\
		\item  {{\bf 6.} Compute the horizontal slices extraction tensor $\mathcal{Q}_{iter}$ by using the indices $\mathbf{J}_{iter}$ and Definition \ref{des18};} \\
		\item {{\bf 7.} Let $\mathcal{V}_{\textit{iter}}$ consist of the eigenslices associated with the top-$d$ eigentubes of $\mathcal{Q}_{\textit{iter}}^T \ast (\mathcal{A} \ast \mathcal{Q}_{\textit{iter}})$;}\\
		\item  {{\bf 8.}  Let  $\mathcal{F}_{iter+1}=\mathcal{Q}_{iter}\ast\mathcal{V}_{iter}$. If
$$
\left|
\frac{{\rm tr}\left(\mathcal{F}_{iter+1}^T \ast \mathcal{A} \ast \mathcal{F}_{iter+1}\right)-{\rm tr}\left(\mathcal{F}_{iter}^T \ast \mathcal{A} \ast \mathcal{F}_{iter}\right)}{{\rm tr}\left(\mathcal{F}_{iter}^T \ast \mathcal{A} \ast \mathcal{F}_{iter}\right)}\right| \leq tol
$$
{\rm then break}; ~{\rm else}~{\rm continue};
}\\
		\item  {{\bf Until convergence}}
\end{algorithmic}
\end{algorithm}

\section{Numerical experiments}\label{sec5}
This part investigates the numerical performance achieved by our proposed algorithm. Our presented methods are tested on multiple real-world high-dimensional large-scale tensor datasets to make comparisons with currently available state-of-the-art approaches. Among our numerical experiments, Examples \ref{ex1}, \ref{ex2} and \ref{ex4} are executed on a Dell workstation equipped with dual 16-core Intel(R) Xeon(R) E5-2620 v4 processors operating at 2.10\,GHz and 64\,GB RAM. The workstation runs 64-bit Windows 10, and all numerical outputs are generated via {\sc Matlab} R2024a.
For the large-scale problem in Example \ref{ex3}, simulations are carried out on a Sugon workstation fitted with dual 64-core Intel(R) Xeon(R) Platinum E5-2637 processors clocked at 3.50\,GHz together with 256\,GB RAM. This machine also employs 64-bit Windows 10, and corresponding numerical outputs are produced using {\sc Matlab} R2018b.


\subsection{Datasets}\label{sec5.1}
Table \ref{Tal2data} summarizes datasets adopted in this section, covering face data, handwritten data and video data. Figure \ref{picvis1} visualize some samples from five datasets.
\begin{table}[htbp]
	\centering
	\caption{Summary description for the utilized datasets}
	\label{Tal2data}
	\begin{tabular}{ccccc}
		\toprule
		Dataset & Samples($n$) & Dimension($d$) & Classes($c$)& Background \\
		\midrule
		COIL20 & 1440 & $128 \times 128$ & 20 & Object images\\
		COIL100 & 7200 & $128 \times 128$ & 100 & Object images\\
		ORL & 400 & $112 \times 92$ & 40 & Face images \\
        GTDB & 750 & $60 \times 80$ & 50 & Face images \\
		Extended YaleB & 2432 & $100 \times 100$ & 38 & Face images \\
		AR2600 & 2600 & $120 \times 165$ & 100 & Face images \\
		AR1680 & 1680 & $50 \times 40$ & 120 & Face images \\
		CFERET & 3528 & $192 \times 128$ & 268 & Face images \\
		MNIST & 4000 & $28 \times 28$ & 10 & Handwritten digits \\
        FEI & 2800 & $480 \times 640$ & 200 & Face images \\
		CAS PEAL & 21840 & $360 \times 480$ & 1040 & Face images \\
		Youtube & 124819 & $320 \times 320$ & 1595 & Video data \\
		\bottomrule
	\end{tabular}
\end{table}
\begin{figure}[htbp]
    \centering
    \begin{subfigure}[b]{0.15\textwidth}
        \centering
        \includegraphics[width=\linewidth]{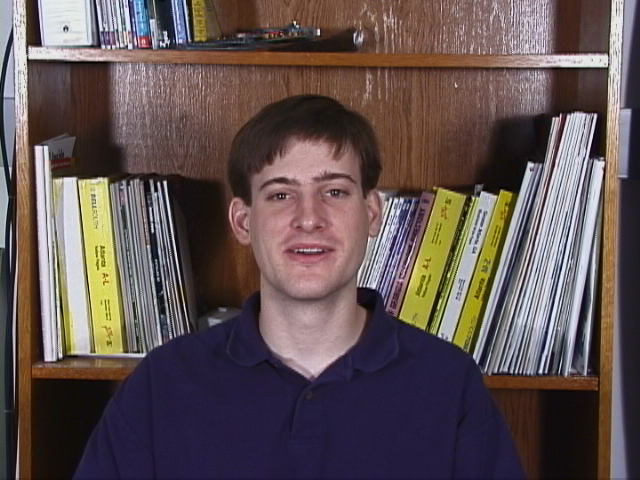}
    \end{subfigure}
    \begin{subfigure}[b]{0.15\textwidth}
        \centering
        \includegraphics[width=\linewidth]{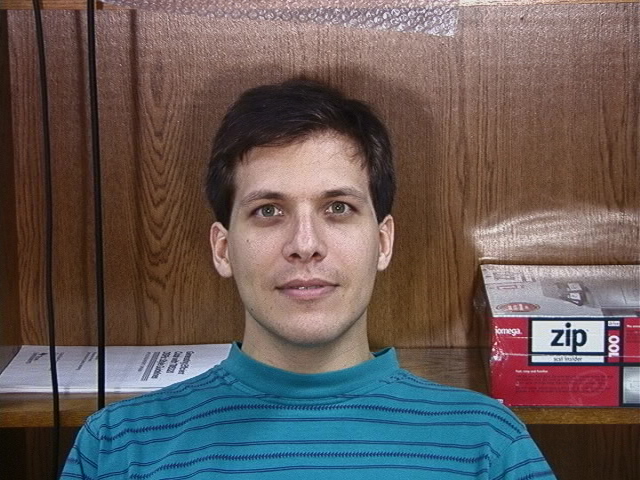}
    \end{subfigure}
    \begin{subfigure}[b]{0.15\textwidth}
        \centering
        \includegraphics[width=\linewidth]{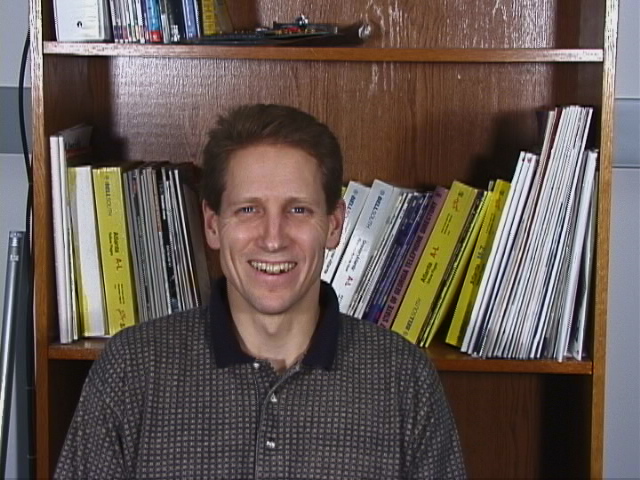}
    \end{subfigure}
    \begin{subfigure}[b]{0.15\textwidth}
        \centering
        \includegraphics[width=\linewidth]{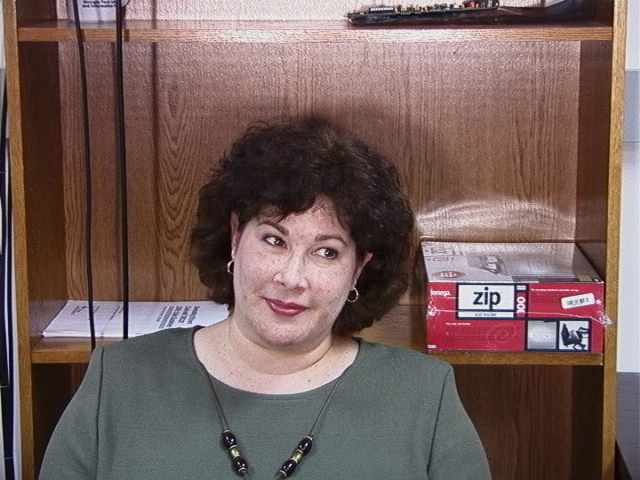}
    \end{subfigure}
    \begin{subfigure}[b]{0.15\textwidth}
        \centering
        \includegraphics[width=\linewidth]{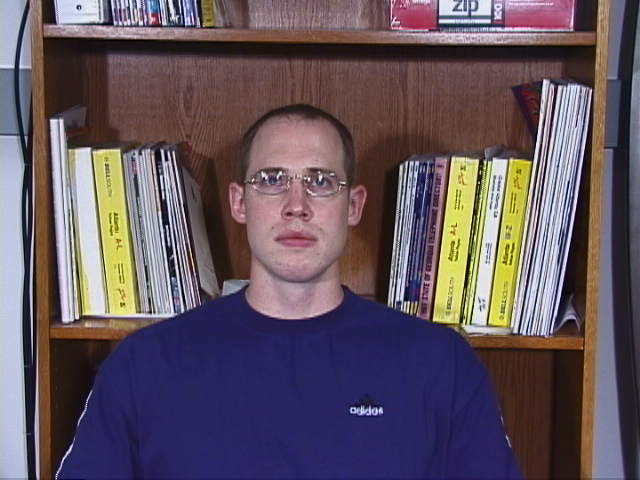}
    \end{subfigure}

    \begin{subfigure}[b]{0.15\textwidth}
        \centering
        \includegraphics[width=\linewidth]{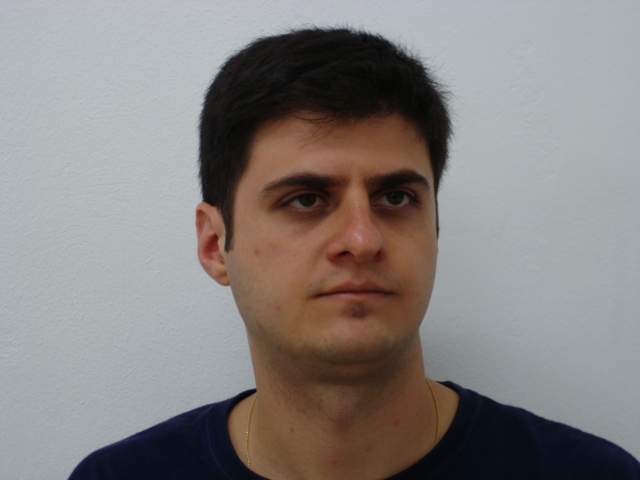}
    \end{subfigure}
    \begin{subfigure}[b]{0.15\textwidth}
        \centering
        \includegraphics[width=\linewidth]{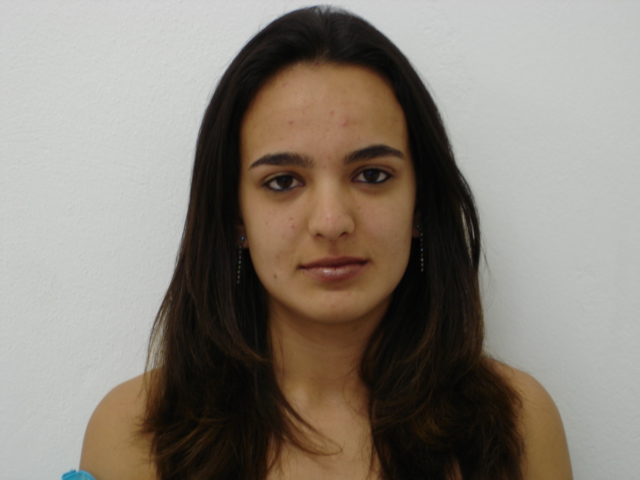}
    \end{subfigure}
    \begin{subfigure}[b]{0.15\textwidth}
        \centering
        \includegraphics[width=\linewidth]{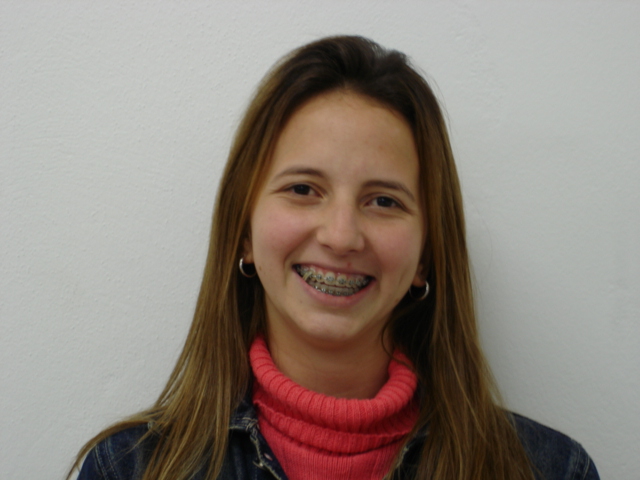}
    \end{subfigure}
    \begin{subfigure}[b]{0.15\textwidth}
        \centering
        \includegraphics[width=\linewidth]{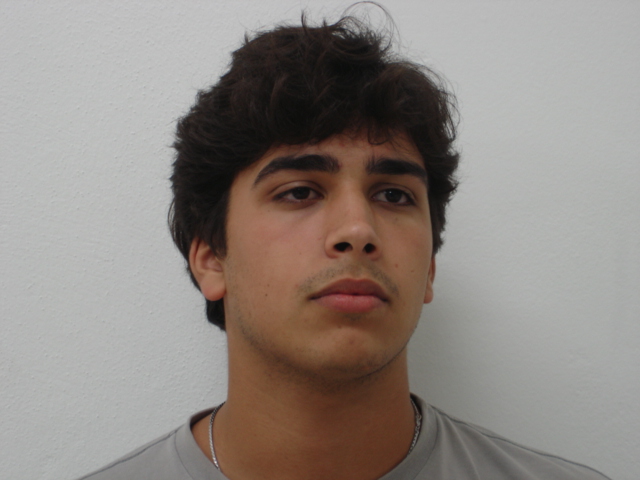}
    \end{subfigure}
    \begin{subfigure}[b]{0.15\textwidth}
        \centering
        \includegraphics[width=\linewidth]{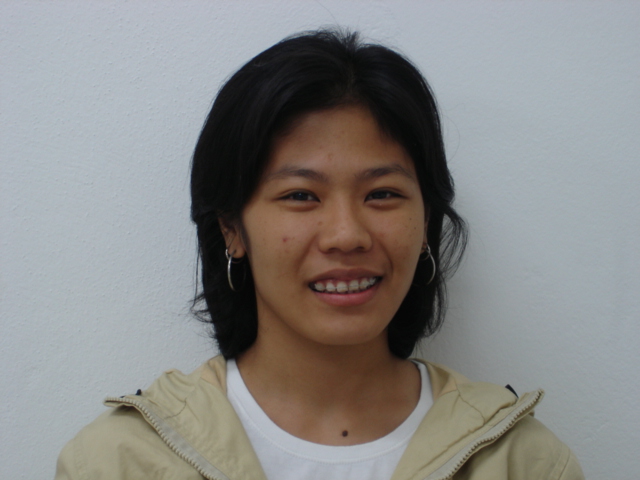}
    \end{subfigure}

    \begin{subfigure}[b]{0.15\textwidth}
        \centering
        \includegraphics[width=\linewidth]{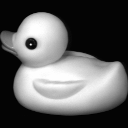}
    \end{subfigure}
    \begin{subfigure}[b]{0.15\textwidth}
        \centering
        \includegraphics[width=\linewidth]{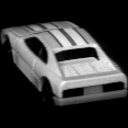}
    \end{subfigure}
    \begin{subfigure}[b]{0.15\textwidth}
        \centering
        \includegraphics[width=\linewidth]{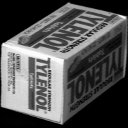}
    \end{subfigure}
    \begin{subfigure}[b]{0.15\textwidth}
        \centering
        \includegraphics[width=\linewidth]{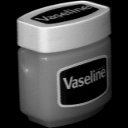}
    \end{subfigure}
    \begin{subfigure}[b]{0.15\textwidth}
        \centering
        \includegraphics[width=\linewidth]{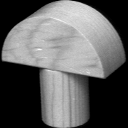}
    \end{subfigure}

    \begin{subfigure}[b]{0.15\textwidth}
        \centering
        \includegraphics[width=\linewidth]{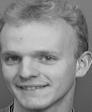}
    \end{subfigure}
    \begin{subfigure}[b]{0.15\textwidth}
        \centering
        \includegraphics[width=\linewidth]{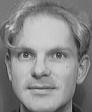}
    \end{subfigure}
    \begin{subfigure}[b]{0.15\textwidth}
        \centering
        \includegraphics[width=\linewidth]{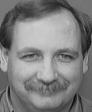}
    \end{subfigure}
    \begin{subfigure}[b]{0.15\textwidth}
        \centering
        \includegraphics[width=\linewidth]{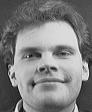}
    \end{subfigure}
    \begin{subfigure}[b]{0.15\textwidth}
        \centering
        \includegraphics[width=\linewidth]{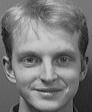}
    \end{subfigure}

    \begin{subfigure}[b]{0.15\textwidth}
        \centering
        \includegraphics[width=\linewidth]{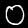}
    \end{subfigure}
    \begin{subfigure}[b]{0.15\textwidth}
        \centering
        \includegraphics[width=\linewidth]{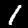}
    \end{subfigure}
    \begin{subfigure}[b]{0.15\textwidth}
        \centering
        \includegraphics[width=\linewidth]{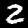}
    \end{subfigure}
    \begin{subfigure}[b]{0.15\textwidth}
        \centering
        \includegraphics[width=\linewidth]{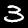}
    \end{subfigure}
    \begin{subfigure}[b]{0.15\textwidth}
        \centering
        \includegraphics[width=\linewidth]{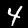}
    \end{subfigure}
    \caption{Representative sample images from {\tt GTDB} (first row), {\tt FEI} (second row), {\tt COIL20} (third row), {\tt ORL} (fourth row), and {\tt MNIST} (final row).} \label{picvis1}
\end{figure}

\begin{itemize}
	\item {\bf ORL} dataset \footnote{\url{https://cam-orl.co.uk/facedatabase.html}}
	comprises images of $40$ individuals captured between April $1992$ and April $1994$. Every subject within this dataset has $10$ corresponding $8$-bit grayscale PGM images at a resolution of $112\times 92$ pixels. These images were captured on separate occasions, exhibiting diversities in illumination, facial expressions and personal facial attributes (e.g., whether subjects wear glasses or not). Every image was collected against a uniformly dark background, where participants kept an upright frontal pose while permitting moderate lateral positional offsets. The dimension of dataset is $112 \times 400 \times 92$.
\end{itemize}

\begin{itemize}
	\item  {\bf Extended YaleB}  dataset
	\footnote{\url{http://cvc.yale.edu/projects/yalefacesB/yalefacesB.html}}
	contains $5760$ single-light-source images across 10 subjects. Each subject is recorded under $576$ viewing configurations, covering $9$ distinct poses and $64$ illumination conditions. The captured samples exhibit normal, sleepy, sad, and surprised facial expressions. For the experiments in this paper, a subset of the dataset is used, consisting of $2432$ images from $38$ individuals ($64$ images per person).  In the experiment, all images are cropped to $100\times 100$ pixels. The dimension of dataset is $100 \times 2432 \times 100$.
\end{itemize}

\begin{itemize}
	\item {\bf AR} dataset
	\footnote{\url{http://www2.ece.ohio-state.edu/~ aleix/ARdatabase.html}}
	dataset holds more than $4000$ colored facial samples from $126$ participants, among whom 70 are male and 56 are female. We select a subset of 100 subjects (50 males and 50 females) for evaluation. This subset contains frontal-view faces with varied expressions, illumination variations, and occlusions including sunglasses and scarves. Recordings were acquired under well-controlled environments, without constraints on participants' clothing, eyewear, makeup, or hairstyle. Each subject was photographed across two separate sessions spaced two weeks apart. Two AR-database subsets are adopted in our numerical tests. The AR2600 subset comprises $2600$ samples from $100$ subjects ($26$ samples for each subject), cropped to $120\times 165$ pixels and yielding a tensor size of $120 \times 2600 \times 165$. The other is the subset AR1680, which includes 1680 facial images of 120 individuals (14 images per person), with all samples resized to $50\times 40$ pixels and a tensor dimension of $50 \times 1680 \times 40$.
\end{itemize}

\begin{itemize}
	\item {\bf Color FERET} dataset
	\footnote{\url{https://www.face-rec.org/databases/}} was collected across 15 recording sessions spanning August 1993 to July 1996. This database comprises 1564 image sets totaling 14126 images, corresponding to 1199 distinct individuals along with 365 repeated image sets. The collected samples cover viewing angles ranging from frontal faces toward left-and right-profile poses. A subset taken from this database serves as our experimental data, containing $3528$ images belonging to $269$ different subjects. The images are cropped to $192\times 128$ pixels. The dimension of dataset is $192 \times 3528 \times 128$.
\end{itemize}

\begin{itemize}
	\item {\bf FEI} dataset
	\footnote{\url{http://www.fei.edu.br/~cet/facedatabase.html/}} consists of RGB facial samples collected within the Artificial Intelligence Laboratory at FEI. Each among the $200$ subjects is represented by 14 photographs, yielding $2800$ samples in total. All shots are taken against a solid white background with subjects facing straight forward, and the profile rotation angles may go up to around $180$ degrees. The size of each image is $480 \times 640$. The dimension of dataset is $480 \times 2800 \times 640$.
\end{itemize}

\begin{itemize}
	\item {\bf Youtube} dataset \cite{Wolf2011Face} is a sizable benchmark built for video classification tasks. This resource contains 80 million YouTube video links labelled with 4800 unique entities derived from knowledge graphs. In the experiment, we use a dataset of $320\times 320$ feature images from 1595 categories, where each category contains at most 90 sample classes. In the experiments, we randomly select 60000 data points as a subset for computation. The dimension of dataset is $320 \times 60000 \times 320$.
\end{itemize}

\begin{itemize}
	\item {\bf CAS-PEAL} face dataset \cite{Gao2008CASPEAL} was built by the Face Recognition Group from JDL, ICT and CAS, with support from the National Hi-Tech Program and ISVISION. This dataset holds 99594 photographs collected from 1040 participants (595 male and 445 female), covering diversified poses, facial expressions, accessories and illumination conditions (PEAL). For every participant, nine evenly-distributed cameras mounted on a horizontal semicircular shelf are deployed to record multi-pose images within a single acquisition. Participants are further required to perform upward and downward gaze movements, yielding 18 additional images over two separate captures. Five expression patterns, six accessory types (three pairs of glasses and three caps), together with fifteen illumination orientations are taken into account. At present, only partial access to this face database is open to researchers. One CAS-PEAL subset provides 21840 samples corresponding to 1040 subjects. This subset covers 21 distinct poses per subject while excluding other interference factors for research usage. Within our experiments, raw samples are cropped into 360 $\times$ 480-pixel resolution. The dimension of the dataset is 360 $\times$ 21840 $\times$ 480.
\end{itemize}

\begin{itemize}
	\item {\bf MNIST} dataset
	\footnote{\url{http://yann.lecun.com/exdb/mnist/}}
comprises 60000 samples for training together with 10000 samples for testing, which contain annotated hand-written digit instances. Each image is formatted as a $28\times 28$ pixel grayscale image. We use a subset of this database as experimental data, which contains 4000 images of 10 different categories of objects. The dimension of dataset is $28 \times 4000 \times 28$.
\end{itemize}

\begin{itemize}
	\item {\bf COIL100} dataset \footnote{\url{https://www.cs.columbia.edu/CAVE/software/softlib/coil-100.php}}
 holds 7200 color samples from 100 distinct objects, yielding 72 snapshots for each object. These items possess varied and intricate geometric features as well as reflective properties. Target objects are positioned upon an electrically-driven turntable with a black backdrop. Full 360-degree rotation of the turntable enables pose variation with respect to a stationary color camera. Image acquisition proceeds at 5-degree angular steps, which generates 72 different viewpoints for every object. Every sample has a resolution of 128 $\times$ 128 pixels. The dimension of dataset is $128 \times 7200 \times 128$.
\end{itemize}

\begin{itemize}
	\item {\bf COIL20} dataset \footnote{\url{http://www.cs.columbia.edu/CAVE/software/softlib/coil-20.php}}
	is a color image collection featuring 20 objects photographed at 5-degree intervals around a rotating turntable, resulting in 72 images per object with uniform $128\times 128$ pixel dimensions. This dataset consists of two subsets: the first one holds 720 unprocessed snapshots for 10 objects, whereas the second part contains 1440 pre-processed samples covering all 20 objects, systematically capturing pose variations across 360 degrees against a uniform background to facilitate evaluations of algorithms sensitive to geometric transformations and visual feature invariance. The dimension of dataset is $128 \times 1440 \times 128$.
\end{itemize}

\begin{itemize}
	\item {\bf GTDB} dataset
	\footnote{\url{https://www.anefian.com/research/face_reco.htm}} comprises 750 colour JPEG samples from 50 subjects, where every individual has precisely 15 snapshots exhibiting varied facial expressions, scaling factors and illumination setups. The original resolution reaches 480 $\times$ 640. For experimental evaluation, these samples are cropped down to 60 $\times$ 80 pixels, and the tensor dimension of this dataset equals 60 $\times$ 750 $\times$ 80.
\end{itemize}

\subsection{The algorithms for comparison and experiment settings}\label{sec5.2}
To demonstrate the numerical behavior of Algorithm \ref{Alg4} and Algorithm \ref{Alg5}, we carry out comparisons against several well-established multilinear as well as linear discriminant approaches. Detailed descriptions of these methods are described in Table \ref{Tal3methods}.
\begin{table}[htbp]
	\centering
	\caption{Details on the compared algorithms}
	\label{Tal3methods}
	\begin{tabular}{lcc}
		\toprule
		\text{Acronym} & \text{Full name} & \text{Linear/Multilinea}r \\
		\midrule
		CMDA \cite{singp2}  & \makecell{Multilinear discriminant analysis for higher-order \\tensor data classification} & Multilinear \\
		DATER \cite{Yan} & \makecell{Discriminant analysis with tensor representation} & Multilinear \\
		UMLDA \cite{Lu} & \makecell{Uncorrelated multilinear discriminant analysis \\ with regularization} & Multilinear \\
		t-MLDE \cite{Bouallala} & \makecell{Trace ratio based manifold learning \\with tensor data} & Multilinear \\
		CF-LDA  \cite{Shi}  & \makecell{A new algorithm for trace-ratio problem with \\ application to high-dimension and large-sample \\ data dimensionality reduction} & Linear \\
		2D-LDA \cite{Ye} & \makecell{Two-dimensional linear discriminant analysis} & Linear \\
		2D-QLDA \cite{Ling2022Joint} & \makecell{Joint diagonalization for a pair of Hermitian \\ quaternion matrices and applications to color \\ face recognition} & Linear \\
		STDA \cite{Lai_Sparse} & \makecell{Sparse tensor discriminant analysis} & Multilinear \\
		SLDA \cite{Nie2020Subspace} & \makecell{Subspace sparse discriminative feature selection} & Linear \\
		\bottomrule
	\end{tabular}
\end{table}

$\bullet$ CMDA\footnote{We thank Dr. Qun Li for providing us the {\sc Matlab} codes of this algorithm.}: Constrained multilinear discriminant analysis \cite{singp2}. This constrained multilinear discriminant analysis technique aims to obtain a collection of projection matrices possessing orthonormal columns through iterative maximization of the scatter ratio criterion.

$\bullet$ DATER: Discriminant analysis with tensor representation \cite{Yan}. It is a discriminant analysis method  with tensor representation in terms of tensor-to-tensor projection.

$\bullet$ UMLDA\footnote{We thank Dr. Haiping Lu for providing us the {\sc Matlab} codes of this algorithm.}: Uncorrelated multilinear discriminant analysis with regularization \cite{Lu}. This multilinear discriminant analysis technique built upon tensor-to-vector projection (TVP) yields uncorrelated discriminative features by means of TVP that optimizes the scalar scatter ratio criterion.

$\bullet$ t-MLDE: Trace-ratio based manifold learning with tensor data \cite{Bouallala}. This method is a generalization of the matrix trace-ratio criterion by using the properties of the t-product.

$\bullet$ CF-LDA:
A randomized algorithm for the trace-ratio problem on high-dimension
and large-sample dense data sets \cite{Shi}. The approach builds on the closed-form solution to the matrix trace-ratio problem together with randomized singular value decomposition.

$\bullet$ 2D-LDA: The two-dimensional linear discriminant analysis \cite{Ye}. As a 2-D linear discriminant analysis technique, it processes data under matrix form to mitigate the high computational overhead arising from eigen-decomposition.

$\bullet$ 2D-QLDA: A structure-preserving algorithm for joint diagonalization of quaternion hermitian matrices \cite{Ling2022Joint}. Designed for color face recognition tasks, this method puts forward an original joint diagonalization framework aimed at Hermitian quaternion matrix pairs, which is derived using a real-valued structure-preserving manipulation.

$\bullet$ STDA: Sparse tensor discriminant analysis \cite{Lai_Sparse}. As a sparse discriminant subspace learning approach, it picks every nonzero entry within each subspace among the most informative variables or factors so as to boost classification performance for tensor-oriented approaches.

$\bullet$ SLDA: Subspace sparse discriminative feature selection \cite{Nie2020Subspace}. A feature selection method that proposes a novel approach by explicitly tackling the long-standing subspace sparsity issue.

For every experimental configuration, we sample $m=80\%N$, $m=70\%N$ and $m=60\%N$ randomly out of the whole sample set for constructing the training subset, while the remaining samples are used for testing. In this notation, $N$ denotes the total sample count. Unless otherwise stated, $d=c-1$ is adopted, in which $c$ represents the class quantity. Besides, the nearest-neighbor (NN) classifier \cite{Cover1967NN} is utilized to carry out classification. Each trial runs for five separate runs, and all reported measurements such as CPU consumption (in seconds), recognition rate and standard deviation (Std-Dev) are computed as the average obtained across these five independent replicates.

\subsection{Numerical examples}\label{sec5.3}

\begin{example} \label{ex1}
This illustrative test contrasts Algorithm \ref{Alg4} against four multilinear discriminant approaches: CMDA \cite{singp2}, DATER \cite{Yan}, UMLDA \cite{Lu} and t-MLDE \cite{Bouallala}, and show the efficiency of our algorithm. The datasets are {\tt ORL}, {\tt Extended YaleB}, {\tt AR1680} and {\tt FEI}. Corresponding numerical outcomes are summarized in Table \ref{table1.1}, Table \ref{table1.2}, Table \ref{table1.3} and Table \ref{table1.4}.

Some remarks are given. First, one can readily notice from these tables that Algorithm \ref{Alg4} achieves considerably higher computational speed compared with the remaining four competing approaches.
We observe from Table \ref{table1.3} that, for the {\tt AR1680} dataset, when $m=80\% N$, Algorithm \ref{Alg4} consumes merely $0.32$ seconds. This indicates roughly $100-200$-fold speedup over t-MLDE, UMLDA, and CMDA. Notably, Algorithm \ref{Alg4} delivers around $1000$-fold acceleration relative to the DATER approach. Second, our method may outperform or underperform several competing techniques in terms of recognition rate; such outcomes are frequently data-dependent. Take the {\tt ORL} dataset as an example, it is seen from Table \ref{table1.1} that, when $m=60\%N$, the accuracy obtained from Algorithm \ref{Alg4} is approximately $91.87\%$, while those of t-MLDE, CMDA and DATER are about $94.63\%$, $94.75\%$ and $93.75\%$, respectively, and the accuracy of UMLDA is around $90.12\%$. In addition, the standard errors of Algorithm \ref{Alg4} are comparable to those of the other four algorithms, indicating that they have similar stability.

Third, for the large-scale problem {\tt FEI} in Table \ref{table1.4}, whose size is $480 \times 2800\times 640$, all the compared algorithms fail to work. Specifically, they encounter the well-known ``out of memory'' issue, or do not work within 1800 seconds. As a comparison, Algorithm \ref{Alg4} provides reliable results within 30 seconds. Our algorithm benefits from the closed-form solution, in contrast to iterative model solvers for competing approaches, which demonstrates its efficiency on large-scale tensor tasks.
\begin{table}[htb]
	\centering
	\caption{Example \ref{ex1}: Experimental outcomes obtained over the {\tt ORL} database, the dataset size is $112 \times 400\times 92$, with the reducing dimension $d = 39$.}\label{table1.1}
	\begin{tabular}{lccc}
		\toprule
		\multirow{2}{*}{\makecell[c]{Algorithms\\(Methods)}} & \multicolumn{3}{c}{CPU Time / s ( Recognition rate $\pm$ Std-Dev\% )} \\
		\cmidrule(lr){2-4}
		& $m = 60\%N$ & $m = 70\%N$ & $m = 80\%N$ \\
		\midrule
		\bf{Algorithm \ref{Alg4}} & 0.31 (91.87$\pm$2.90\%) & 0.34 (92.83$\pm$1.63\%) & 0.37 (93.50$\pm$3.10\%) \\
		t-MLDE & 33.70 (94.63$\pm$2.70\%) & 33.73 (94.83$\pm$1.33\%) & 34.68 (96.25$\pm$1.77\%) \\
		UMLDA & 65.27 (90.12$\pm$1.55\%) & 70.86 (92.00$\pm$3.23\%) & 75.99 (95.00$\pm$3.62\%) \\
		CMDA & 691.97(94.75$\pm$0.64\%) & 751.83 (95.17$\pm$1.43\%) & 838.58 (96.25$\pm$1.58\%) \\
		DATER & 696.19 (93.75$\pm$2.68\%) & 778.46 (94.67$\pm$2.10\%) & 856.46 (96.00$\pm$2.67\%) \\
		\bottomrule
	\end{tabular}
\end{table}
\begin{table}[htb]
	\centering
	\caption{Example \ref{ex1}: Experimental outcomes obtained over the {\tt Extended YaleB} database, the dataset size is $100 \times 2432\times 100$,  with the reducing dimension $d = 37$.}\label{table1.2}
	\begin{tabular}{lccc}
		\toprule
		\multirow{2}{*}{\makecell[c]{Algorithms\\(Methods)}} & \multicolumn{3}{c}{CPU Time / s ( Recognition rate $\pm$ Std-Dev\% )} \\
		\cmidrule(lr){2-4}
		& $m = 60\%N$ & $m = 70\%N$ & $m = 80\%N$ \\
		\midrule
		\bf{Algorithm \ref{Alg4}} & 0.92 (89.91$\pm$0.30\%) & 1.17 (90.74$\pm$0.80\%) & 1.36 (91.13$\pm$0.60\%) \\
		t-MLDE & 46.54 (70.63$\pm$0.62\%) & 59.42 (72.77$\pm$1.25\%) & 67.34 (74.13$\pm$1.89\%) \\
		UMLDA & 232.16 (98.64$\pm$0.18\%) & 262.91 (99.01$\pm$0.34\%) & 302.9 (99.10$\pm$0.56\%) \\
		CMDA & -- & -- & -- \\
		DATER & -- & -- & -- \\
		\bottomrule
	\end{tabular}
	\vspace{0.5em}
	\begin{tablenotes}
		\footnotesize
		\item[1] Here, ``$-$'' refers to CPU computation time exceeding 1800 seconds.
	\end{tablenotes}
\end{table}

\begin{table}[htb]
	\centering
	\caption{Example \ref{ex1}: Experimental outcomes obtained over the {\tt AR1680} database, the dataset size is $50 \times 1680\times 40$,  with the reducing dimension $d = 25$.}\label{table1.3}
	\begin{tabular}{lccc}
		\toprule
		\multirow{2}{*}{\makecell[c]{Algorithms\\(Methods)}} & \multicolumn{3}{c}{CPU Time / s ( Recognition rate $\pm$ Std-Dev\% )} \\
		\cmidrule(lr){2-4}
		& $m = 60\%N$ & $m = 70\%N$ & $m = 80\%N$ \\
		\midrule
		\bf{Algorithm \ref{Alg4}} & 0.27 (92.41$\pm$0.60\%) & 0.30 (94.68$\pm$0.64\%)& 0.32 (96.13$\pm$1.22\%) \\
		t-MLDE & 39.75 (90.62$\pm$0.45\%) & 41.19 (92.86$\pm$1.43\%) & 54.02 (94.94$\pm$1.63\%) \\
		UMLDA & 36.33 (97.56$\pm$0.67\%) & 39.63 (98.17$\pm$0.46\%) & 43.99 (98.87$\pm$0.48\%) \\
		CMDA & 61.58 (96.03$\pm$0.66\%) & 64.60 (97.46$\pm$0.83\%) & 69.30 (97.67$\pm$0.48\%) \\
		DATER & 312.30 (97.78$\pm$0.73\%) & 338.68 (98.5$\pm$0.31\%) & 345.36 (98.61$\pm$0.39\%) \\
		\bottomrule
	\end{tabular}
\end{table}
\begin{table}[htbp]
	\begin{center}
		\centering
		\caption{Example \ref{ex1}: Experimental outcomes obtained over the {\tt FEI}, the dataset size is $480 \times 2800\times 640$,  with the reducing dimension $d = 199$.}\label{table1.4}
		\begin{tabular}{lccc}
			\toprule
			\multirow{2}{*}{\makecell[c]{Algorithms\\(Methods)}} & \multicolumn{3}{c}{CPU Time / s ( Recognition rate $\pm$ Std-Dev\% )} \\
			\cmidrule(lr){2-4}
			& $m = 60\%N$ & $m = 70\%N$ & $m = 80\%N$ \\
			\midrule
			{\bf Algorithm \ref{Alg4}} & 22.74 (77.45$\pm$1.45\%) & 24.91 (78.81$\pm$0.86\%) & 26.65 (81.39$\pm$1.37\%) \\
			t-MLDE & -- & -- & -- \\
			UMLDA & -- & -- & -- \\
			CMDA & O.M. & O.M. & O.M. \\
			DATER & O.M. & O.M. & O.M. \\
			\bottomrule
		\end{tabular}
		\vspace{0.5em}
		\begin{tablenotes}
			\footnotesize
			\item[1] Here, ``O.M.'' indicates that the algorithm runs into an out-of-memory failure, while ``$-$'' refers to CPU computation time exceeding 1800 seconds.
		\end{tablenotes}
	\end{center}
\end{table}
\end{example}

\begin{example}\label{ex2}
	For the present illustration, Algorithm \ref{Alg4} is contrasted against several linear discriminant approaches including CF-LDA \cite{Shi}, 2D-LDA \cite{Ye} and 2D-QLDA \cite{Ling2022Joint}. The datasets used in this example are {\tt COIL100} and {\tt CFERET}. The corresponding numerical outcomes are summarized in Table \ref{teble2.1}--Table \ref{table2.2}.

Several observations can be drawn from the tabulated results. As observed in the two tables, Algorithm \ref{Alg4} runs marginally slower than the CF-LDA approach developed in \cite{Shi}. Even so, our recognition accuracies exceed or match those of the reference method. Second, compared with the 2D-LDA method, Algorithm \ref{Alg4} is much faster. As illustrated by Table \ref{table2.2}, when $m=80\%N$, Algorithm \ref{Alg4} takes $6.12$ seconds, while 2D-LDA takes $32.32$ seconds. In terms of recognition rate, our approach can outperform 2D-LDA on certain datasets, depending on the datasets utilized. For example, as is shown in Table \ref{teble2.1} on the {\tt COIL100} dataset, when $m=70\%N$, the accuracy of Algorithm \ref{Alg4} reaches approximately $97.21\%$, whereas that of 2D-LDA is $90.34\%$. That is, our method outperforms 2D-LDA by around $7\%$ in accuracy. However, on the {\tt CFERET} dataset, when $m=60\%N$, Table \ref{table2.2} shows Algorithm \ref{Alg4} attains an accuracy of roughly 34.89\%, versus 41.36\% obtained by 2D-LDA. Meanwhile, Table \ref{table2.2} also reveals that the standard errors for 2D-LDA is higher than those of Algorithm \ref{Alg4}, which means that Algorithm \ref{Alg4} is relatively more stable.

Third, it is observed that the 2D-QLDA method fails to work within 1800 seconds. This is because it incurs highly computational costs on quaternion matrix computations, making it difficult to compute the projection quaternion matrix in practice. To show these more precisely, we visualize some performances of Algorithm \ref{Alg4}, the multilinear methods and the linear methods on the {\tt GTDB} face dataset and the {\tt MNIST} handwritten digit dataset. Figure \ref{compic} and Figure \ref{comtimes} illustrate accuracy (ACC) together with CPU runtime for all competing algorithms across various dimensions, respectively. It is seen that the numerical performances of our proposed algorithm are satisfactory.

\begin{table}[htbp]
	\centering
	\caption{Example \ref{ex2}: Experimental outcomes obtained over the {\tt COIL100}, the dataset size is $128 \times 7200\times 128$, the reducing dimension $d = 99$.}\label{teble2.1}
	\begin{tabular}{lccc}
		\toprule
		\multirow{2}{*}{\makecell[c]{Algorithms\\(Methods)}} & \multicolumn{3}{c}{CPU Time / s ( Recognition rate $\pm$ Std-Dev\% )} \\
		\cmidrule(lr){2-4}
		& $m = 60\%N$ & $m = 70\%N$ & $m = 80\%N$ \\
		\midrule
		\bf{Algorithm \ref{Alg4}} & 2.12 (96.72$\pm$0.36\%) & 2.40(97.21$\pm$0.15\%) & 2.76 (97.69$\pm$0.39\%) \\
		CF-LDA & 2.41 (88.38$\pm$0.46\%) & 2.73 (88.91$\pm$0.73\%) & 3.04 (89.26$\pm$0.75\%) \\
		2D-LDA & 33.97 (89.07$\pm$1.13\%) & 39.04 (90.34$\pm$0.96\%) & 43.76 (91.64$\pm$1.07\%) \\
		2D-QLDA & -- & -- & -- \\
		\bottomrule
	\end{tabular}
	\vspace{0.5em}
	\begin{tablenotes}
		\footnotesize
		\item[1] Here, ``$-$'' refers to CPU computation time exceeding 1800 seconds.
	\end{tablenotes}
\end{table}
\begin{table}[htbp]
	\centering
	\caption{Example \ref{ex2}: Experimental outcomes obtained over the {\tt CFERET}, the dataset size is $192 \times 3528\times 128$, the reducing dimension $d = 90$.}\label{table2.2}
	\begin{tabular}{lccc}
		\toprule
		\multirow{2}{*}{\makecell[c]{Algorithms\\(Methods)}} & \multicolumn{3}{c}{CPU Time / s ( Recognition rate $\pm$ Std-Dev\% )} \\
		\cmidrule(lr){2-4}
		& $m = 60\%N$ & $m = 70\%N$ & $m = 80\%N$ \\
		\midrule
		\bf{Algorithm \ref{Alg4}} & 5.13 (34.89$\pm$0.63\%) & 5.65 (37.32$\pm$2.21\%) & 6.12 (39.60$\pm$1.66\%) \\
		CF-LDA & 1.39 (35.84$\pm$0.97\%) & 1.53 (39.81$\pm$1.22\%) & 1.64 (42.97$\pm$0.82\%) \\
		2D-LDA & 24.92 (41.36$\pm$33.47\%) & 28.55 (56.77$\pm$28.25\%) & 32.32 (73.54$\pm$2.04\%) \\
		2D-QLDA & -- & -- & -- \\
		\bottomrule
	\end{tabular}
	\vspace{0.5em}
	\begin{tablenotes}
		\footnotesize
		\item[1] Here, ``$-$'' refers to CPU computation time exceeding 1800 seconds.
	\end{tablenotes}
\end{table}
\end{example}

\begin{figure}[htbp]
    \centering
    \begin{subfigure}[b]{0.5\textwidth}
        \centering
        \includegraphics[width=6.5cm,height=5cm]{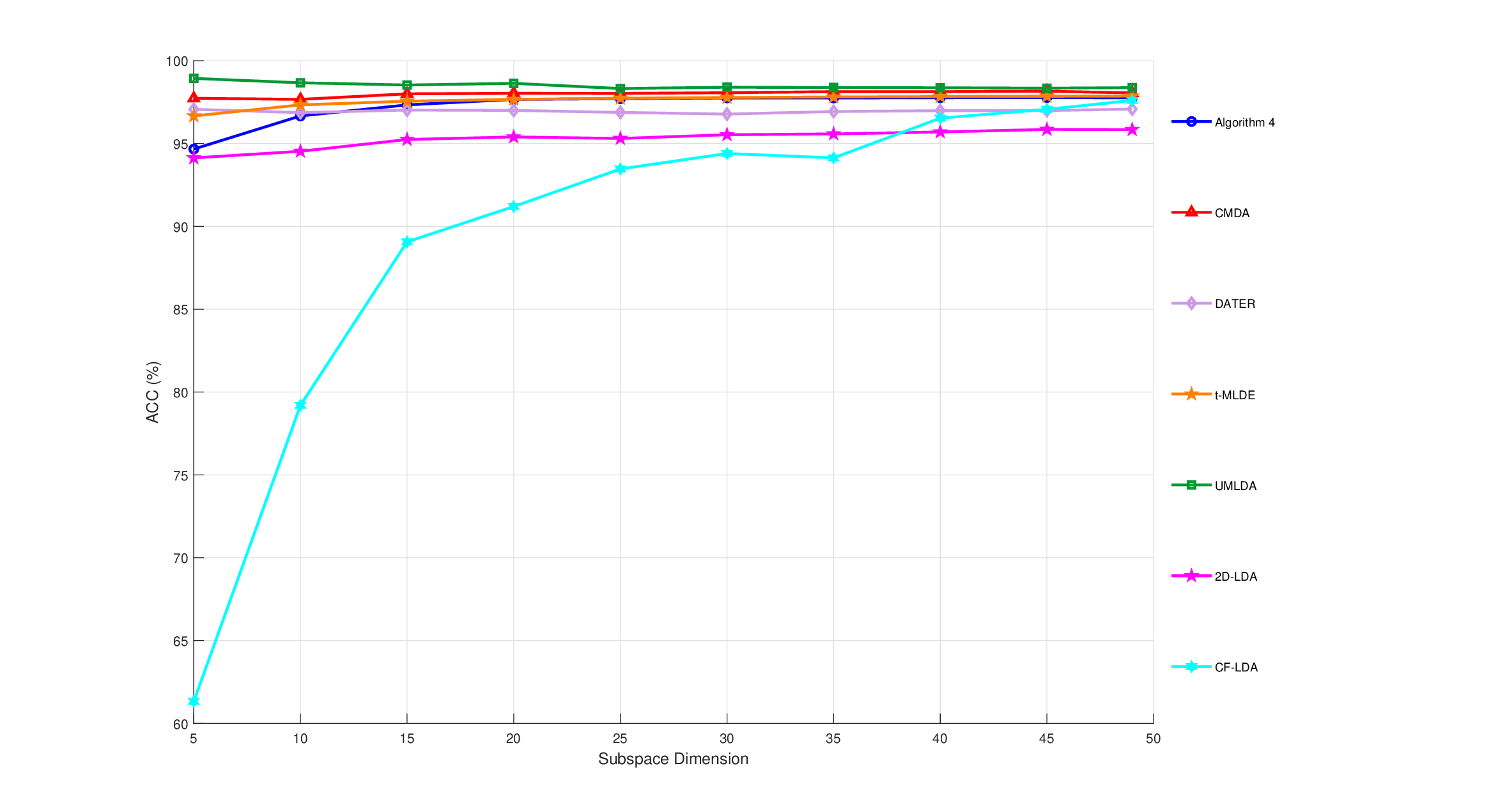}
        \caption{{\tt The GTDB dataset}}
    \end{subfigure}
    \hspace{-0.8cm}
    \begin{subfigure}[b]{0.5\textwidth}
        \centering
        \includegraphics[width=6.5cm,height=5cm]{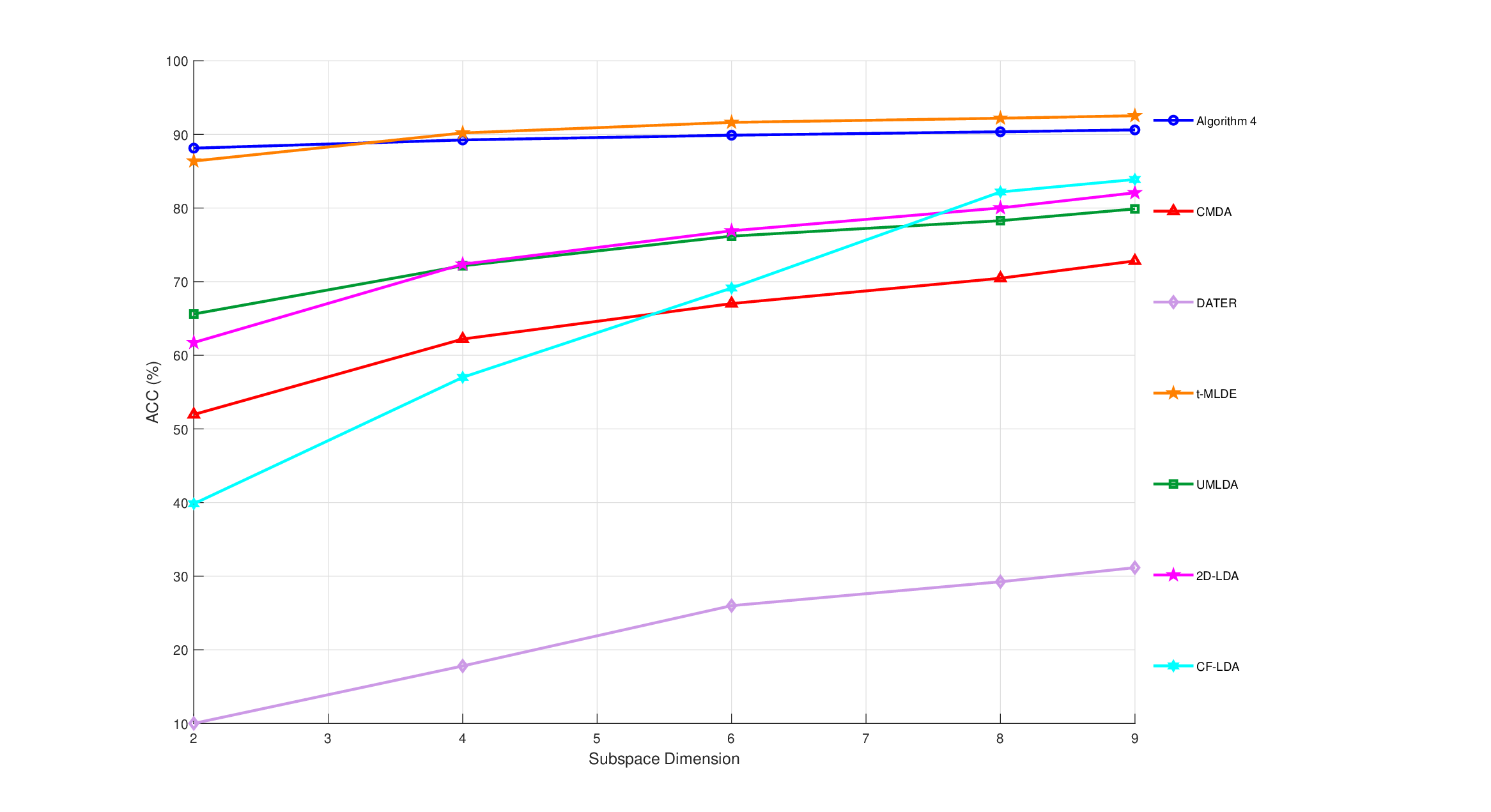}
        \caption{{\tt The MNIST dataset}}
    \end{subfigure}
    \caption{Example \ref{ex2}: A comparison of recognition rate (ACC) of the compared algorithms on the two datasets {\tt GTDB} and {\tt MNIST}.}
    \label{compic}
\end{figure}
\begin{figure}[htbp]
\centering
\begin{subfigure}[b]{1\textwidth}
  \centering
  \includegraphics[width=\linewidth,height=6cm]{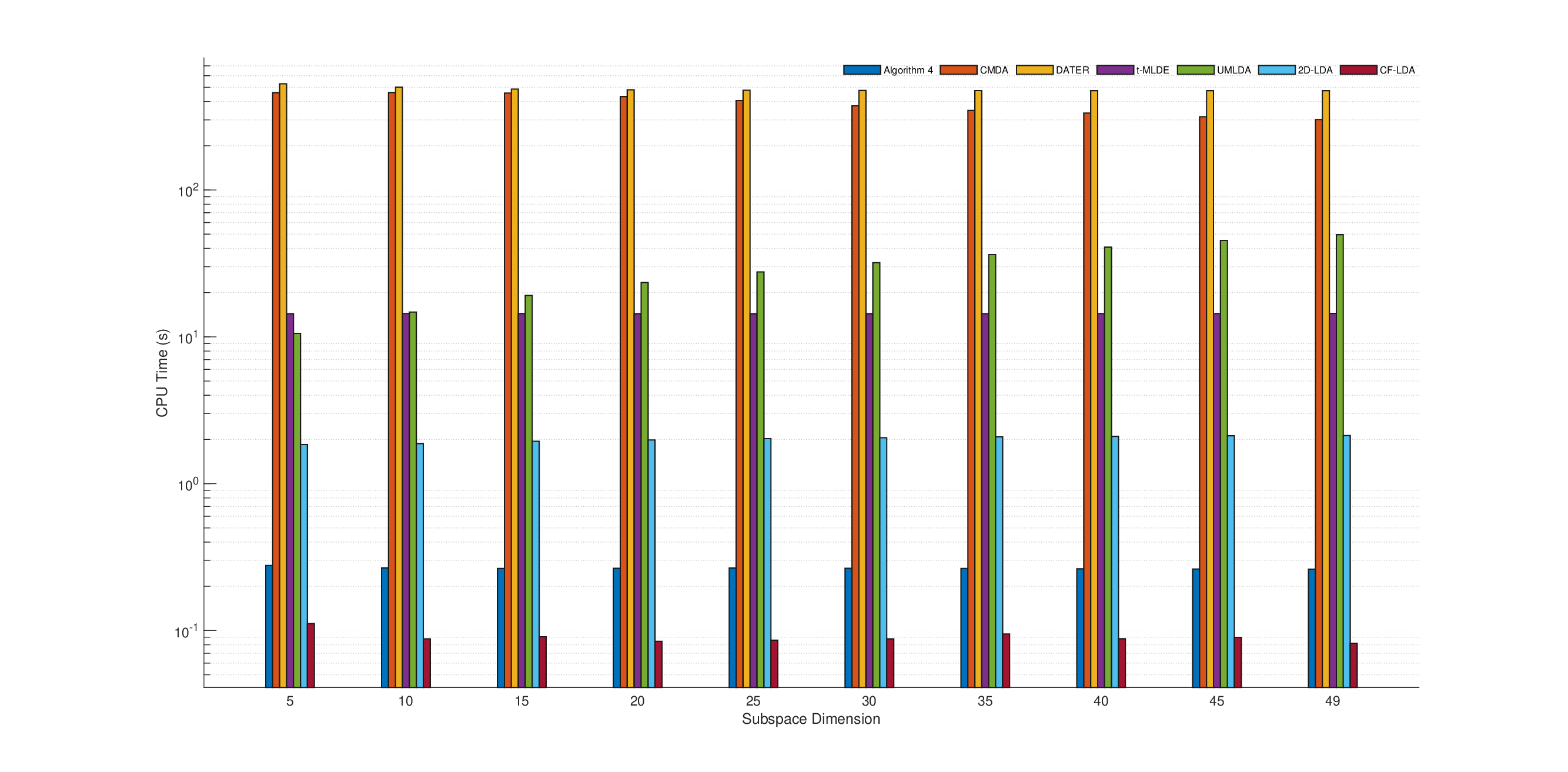}
  \caption{{\tt The GTDB dataset}}
\end{subfigure}
\hfill  
\begin{subfigure}[b]{1\textwidth}
  \centering
  \includegraphics[width=\linewidth,height=7cm]{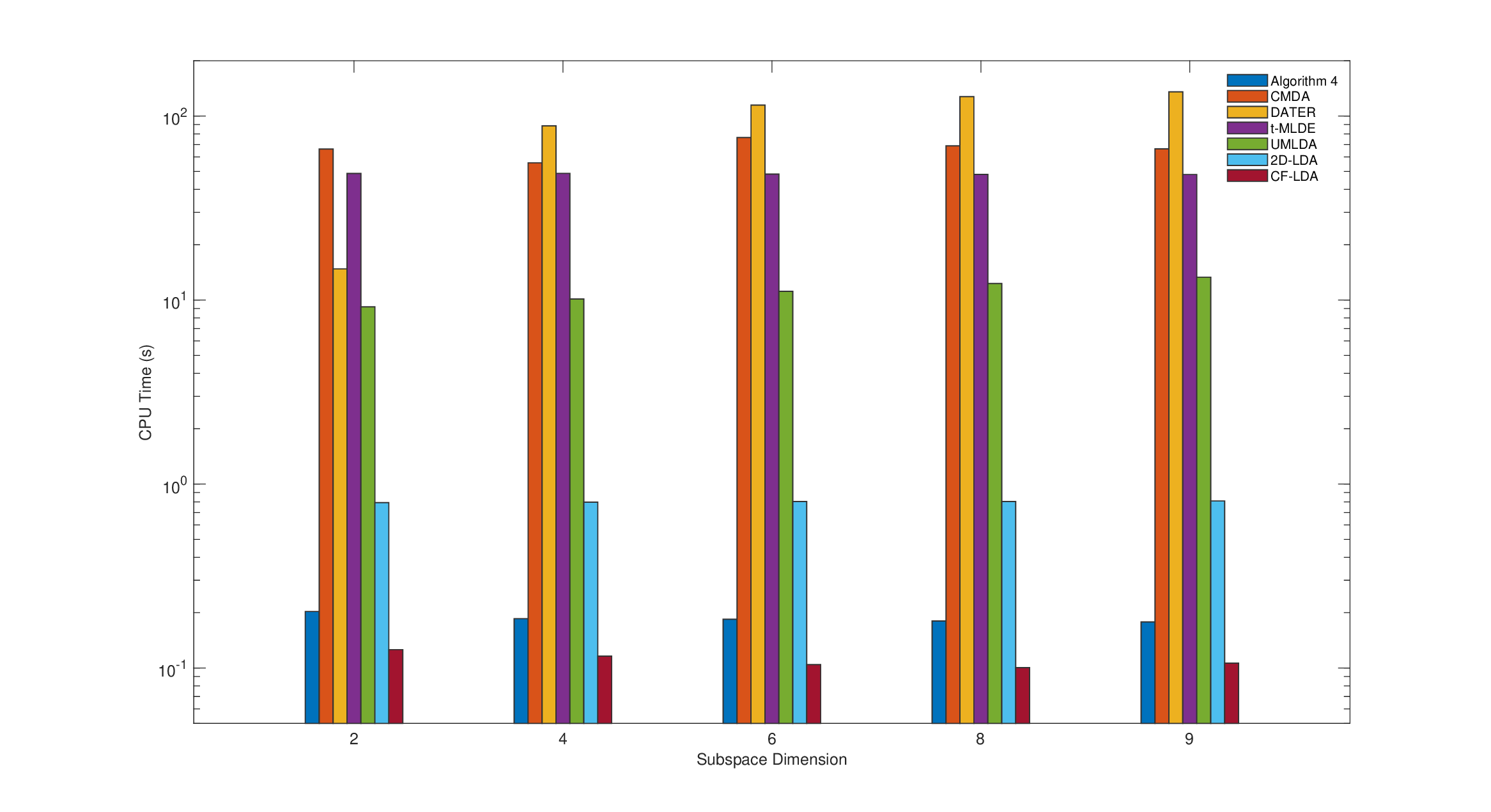}
  \caption{{\tt The MNIST dataset}}
\end{subfigure}
\caption{Example \ref{ex2}: Comparison on CPU runtime among competing algorithms using the {\tt GTDB} and {\tt MNIST} datasets.}
\label{comtimes}
\end{figure}

\begin{example}\label{ex3}
This numerical case illustrates the advantages delivered by Algorithm \ref{Alg4} for large-scale tensor datasets. We compare Algorithm \ref{Alg4} with the first seven multilinear and linear methods listed
in Table \ref{Tal3methods}. The datasets are two high-dimensional large-sample tensors. The first one is the {\tt CAS PEAL} dataset with dimension $360\times 21840\times 480$, and the second one is the {\tt Youtube} dataset with dimension $320\times 60000\times 320$. The corresponding numerical outcomes are summarized within Table \ref{Table3.1} and Table \ref{Table3.2}.

Some remarks are given. First, for the two large-scale problems, we see from the two tables that only Algorithm \ref{Alg4} and CF-LDA work, while the other six algorithms suffer from ``out of memory" or heavy computational costs. Second, Algorithm \ref{Alg4} outperforms CF-LDA in terms of recognition rate. For instance, as can be seen from Table \ref{Table3.2} on the {\tt Youtube} dataset, when $m=80\%$, the recognition rate of Algorithm \ref{Alg4} is approximately $99.23\%$, while that of CF-LDA stands at around $92.34\%$. Thus, the recognition rate of Algorithm \ref{Alg4} is approximately $7\%$ higher than that of CF-LDA. This is due to the fact that Algorithm \ref{Alg4} solves a tensor optimization model and retains the original structure of tensor data in some degree.

Again, our approach incurs larger computational overhead relative to the CF-LDA method. As illustrated in Table \ref{Table3.1} for the {\tt CAS PEAL} dataset, Algorithm \ref{Alg4} takes $264.52$ seconds when $m=70\%$, while CF-LDA consumes $86.31$ seconds. Combining recognition accuracy and workload as a whole, Algorithm \ref{Alg4} represents a favorable alternative for large-scale tensor-data  problems.

\begin{table}[htbp]
	\centering
	\caption{Example \ref{ex3}: Experimental outcomes obtained over the {\tt CAS PEAL}, the dataset size is $360 \times 21840\times 480$, the reducing dimension $d = 200$.} \label{Table3.1}
	\begin{tabular}{lccc}
		\toprule
		\multirow{2}{*}{\makecell[c]{Algorithms\\(Methods)}} & \multicolumn{3}{c}{CPU Time / s ( Recognition rate $\pm$ Std-Dev\% )} \\
		\cmidrule(lr){2-4}
		& $m = 60\%N$ & $m = 70\%N$ & $m = 80\%N$ \\
		\midrule
		\bf{Algorithm \ref{Alg4}} & 220.07 (48.21$\pm$0.42\%) & 264.52 (51.63$\pm$0.46\%) & 306.71 (55.08$\pm$0.56\%) \\
		t-MLDE & -- & -- & -- \\
		UMLDA & -- & -- & -- \\
		CMDA & O.M. & O.M. & O.M. \\
		DATER & O.M. & O.M. & O.M. \\
		CF-LDA & 80.94 (29.94$\pm$0.54\%) & 86.31(35.50$\pm$0.73\%) & 96.32 (40.11$\pm$0.58\%) \\
		2D-LDA & -- & -- & -- \\
		2D-QLDA & -- & -- & -- \\
		\bottomrule
	\end{tabular}
	\vspace{0.5em}
	\begin{tablenotes}
		\footnotesize
		\item[1] Here, ``O.M.'' indicates that the algorithm runs into an out-of-memory failure, while ``$-$'' refers to CPU computation time exceeding 1800 seconds.
	\end{tablenotes}
\end{table}
\begin{table}[htbp]
	\centering
	\caption{Example \ref{ex3}: Experimental outcomes obtained over the {\tt Youtube}, the dataset size is $320 \times 60000\times 320$, the reducing dimension $d = 200$.} \label{Table3.2}
	\begin{tabular}{lccc}
		\toprule
		\multirow{2}{*}{\makecell[c]{Algorithms\\(Methods)}} & \multicolumn{3}{c}{CPU Time / s ( Recognition rate $\pm$ Std-Dev\% )} \\
		\cmidrule(lr){2-4}
		& $m = 60\%N$ & $m = 70\%N$ & $m = 80\%N$ \\
		\midrule
		\bf{Algorithm \ref{Alg4}} & 513.50(99.00$\pm$0.07\%) & 998.38(99.14$\pm$0.08\%) & 1483.30(99.23$\pm$0.10\%) \\
		t-MLDE & -- & -- & -- \\
		UMLDA & -- & -- & -- \\
		CMDA & O.M. & O.M. & O.M. \\
		DATER & O.M. & O.M. & O.M. \\
		CF-LDA & 118.03 (91.98$\pm$0.26\%) & 131.1(92.13$\pm$0.21\%) & 146.42 (92.34$\pm$0.16\%) \\
		2D-LDA & -- & -- & -- \\
		2D-QLDA & -- & -- & -- \\
		\bottomrule
	\end{tabular}
	\vspace{0.5em}
	\begin{tablenotes}
		\footnotesize
		\item[1] Here, ``O.M.'' indicates that the algorithm runs into an out-of-memory failure, while ``$-$'' refers to CPU computation time exceeding 1800 seconds.
	\end{tablenotes}
\end{table}
\end{example}

\begin{example}\label{ex4}
This numerical example demonstrates the advantages offered by the sparse multilinear discriminant analysis method (Algorithm \ref{Alg5}) for tensor feature selection.
The datasets used are the {\tt COIL20} dataset of size $128 \times 1440\times 128$, and the {\tt AR2600} dataset of size $120 \times 2600\times 165$.
We run three feature selection algorithms including Algorithm \ref{Alg5}, STDA \cite{Lai_Sparse} and SLDA \cite{Nie2020Subspace}, where the stopping tolerance used in Algorithm \ref{Alg5}
is chosen as ${\it tol}=10^{-3}$. As a by-product, we also run eight feature extraction methods within this test case. The corresponding numerical outcomes are reported in Table \ref{Table4.1} and Table \ref{Table4.2}.

In contrast to these two feature-selection approaches, the results demonstrate that Algorithm \ref{Alg5} runs more efficiently than STDA and attains superior classification accuracy. However, the SLDA method fails to converge within the specified time. Moreover, among the feature extraction methods included in comparison, CMDA, DATER and 2D-QLDA also fail to work within 1800 seconds. On the other hand, compared with Algorithm \ref{Alg4}, the recognition accuracy of the sparse version Algorithm \ref{Alg5} shows a noticeable improvement after integrating feature selection. More precisely, Algorithm \ref{Alg5} achieves roughly 8\% higher recognition accuracy than Algorithm \ref{Alg4}. One reason is that Algorithm \ref{Alg5} significantly mitigates data redundancy and the curse of dimensionality, owing to our introduced feature-selection strategy.

In contrast to alternative feature-extraction approaches, we see from Table \ref{Table4.1} that, for the {\tt COIL20} dataset, the recognition accuracies of Algorithm \ref{Alg5} are comparable or a little higher than those of the others. For the {\tt AR2600} dataset, we see from Table \ref{Table4.2} that
Algorithm \ref{Alg4} achieves lower classification accuracy than some competing methods, and Algorithm \ref{Alg5} also performs worse than most feature extraction methods yet surpasses t-MLDE.
Furthermore, Algorithm \ref{Alg5} is often (much) faster than the others.
Thus, it remains hard to identify the optimal method, as the obtained outcomes are heavily problem-dependent. In addition, Figure \ref{relitu} visualizes the projector tensor $\mathcal{F}$, where blue cubes represent the zero-horizontal slices, i.e., the unselected features in the first feature dimension. This demonstrates that Algorithm \ref{Alg5} is a powerful candidate for high-dimensional tensor data feature selection.

\begin{figure}[htbp]
\makebox[\textwidth][c]{%
  \begin{subfigure}[b]{20cm}
    \centering 
    \includegraphics[width=20cm,height=6cm]{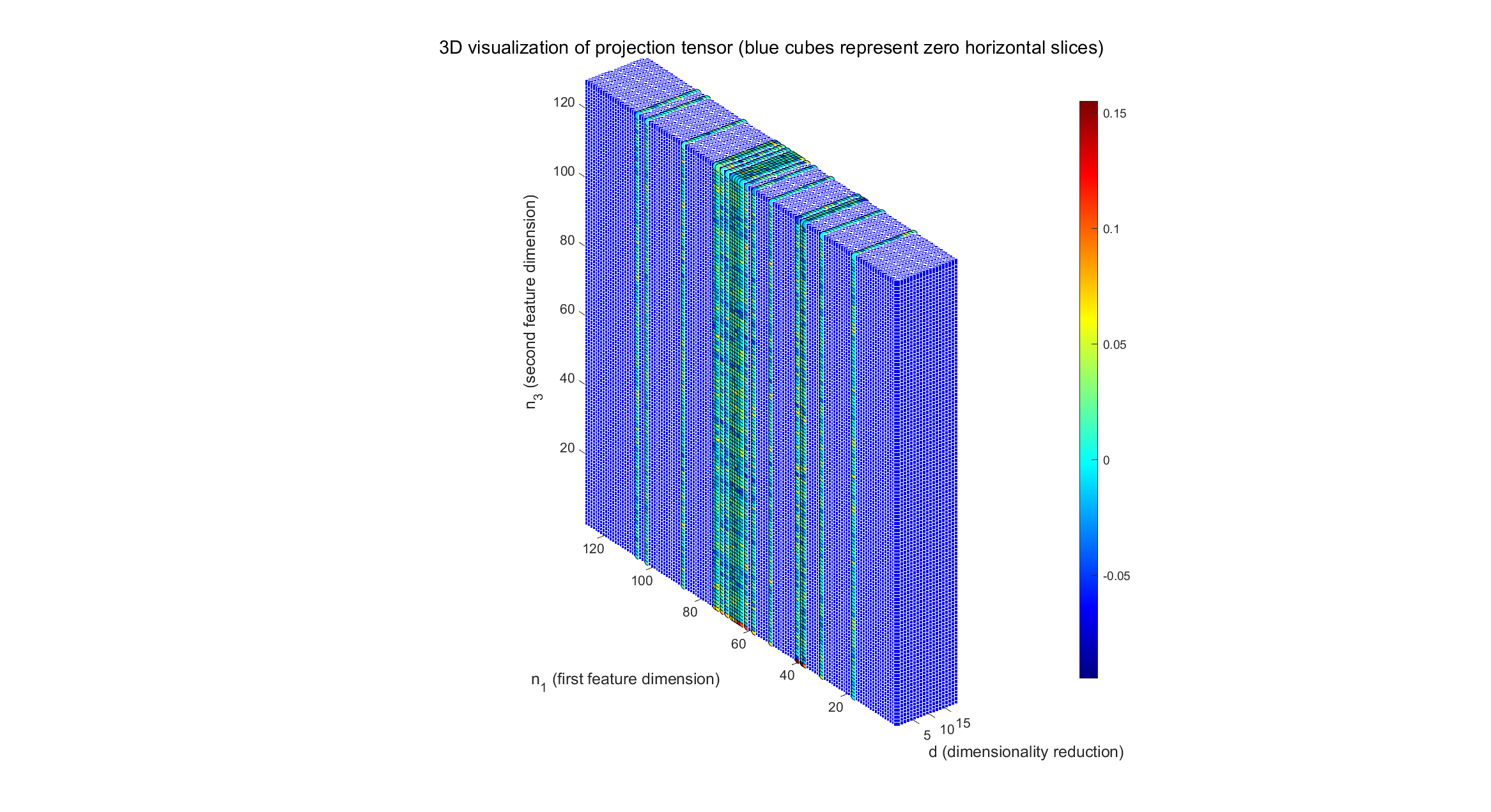}
    \caption{{\tt The COIL20 dataset}}
  \end{subfigure}
}

\vspace{0.5em} 

\makebox[\textwidth][c]{%
  \begin{subfigure}[b]{20cm}
    \centering
    \includegraphics[width=20cm,height=6cm]{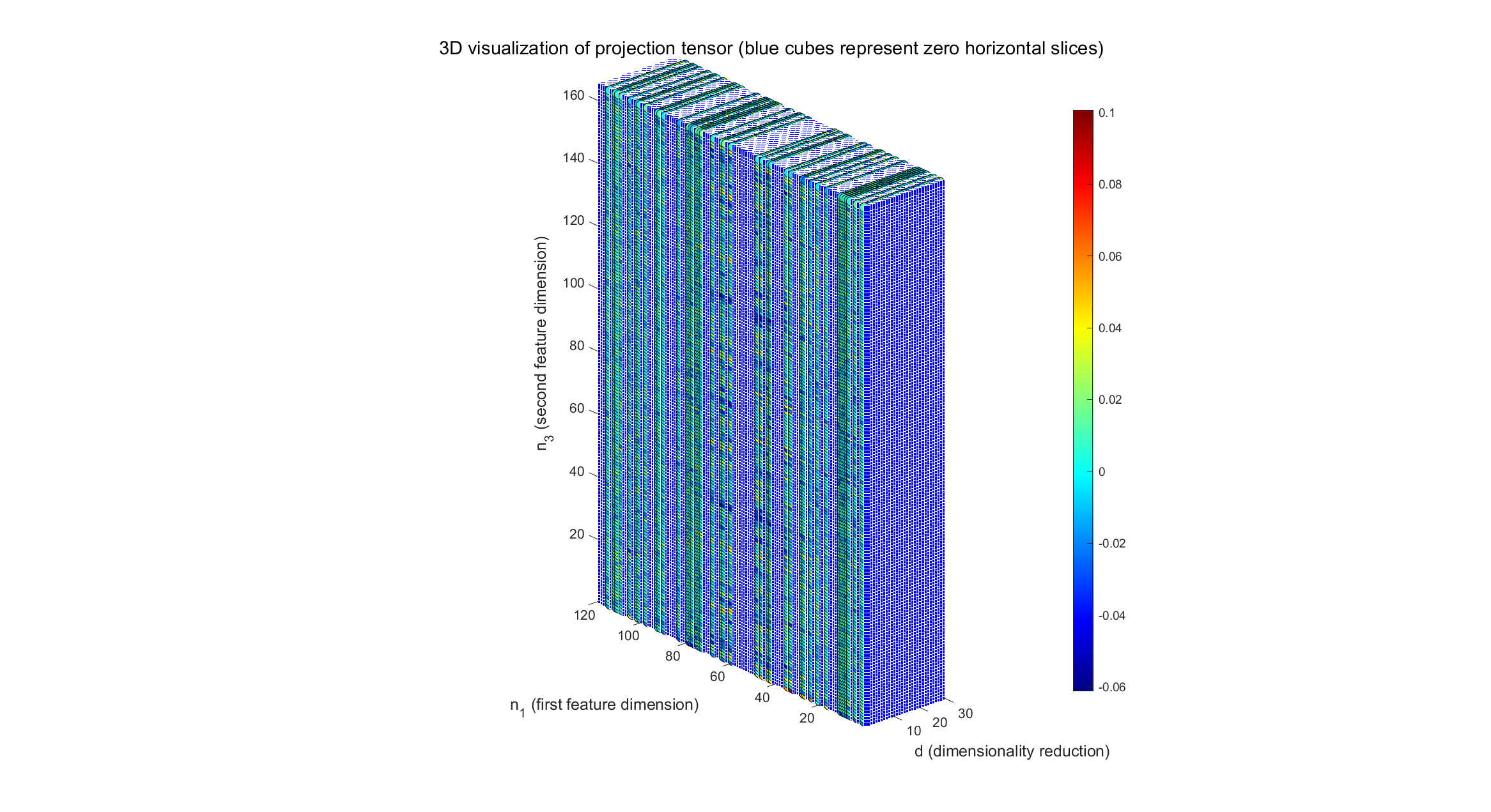}
    \caption{{\tt The AR2600 dataset}}
  \end{subfigure}
}

\caption{Example \ref{ex4}: Visualization of the projector tensor $\mathcal{F}$ on the two datasets {\tt COIL20} and {\tt AR2600}.}
\label{relitu}
\end{figure}
\begin{table}[h]
	\centering
	\caption{Example \ref{ex4}: Experimental outcomes obtained over the {\tt COIL20} dataset, whose size is $128 \times 1440\times 128$. The reducing dimension $d = 19$, while the count of chosen features is $k=20$.}  \label{Table4.1}
	\begin{tabular}{lccc}
		\toprule
		\multirow{2}{*}{\makecell[c]{Algorithms\\(Methods)}} & \multicolumn{3}{c}{CPU Time / s ( Recognition rate $\pm$ Std-Dev\% )} \\
		\cmidrule(lr){2-4}
		& $m = 60\%N$ & $m = 70\%N$ & $m = 80\%N$ \\
		\midrule
     \bf{Algorithm \ref{Alg5}} & 2.27(99.03$\pm$0.57\%)  & 2.52(99.35$\pm$0.23\%)  & 3.43(99.65$\pm$0.22\%)  \\
		STDA & 4.63(98.59$\pm$1.12\%)&4.84(99.36$\pm$0.06\%)  & 4.91(99.87$\pm$0.08\%)  \\
		SLDA & -- & -- & -- \\\hline
		\bf{Algorithm \ref{Alg4}} & 0.42(90.07$\pm$1.50\%)  & 0.47(90.95$\pm$1.97\%)  & 0.49(92.22$\pm$1.21\%)  \\
		t-MLDE & 71.60(99.06$\pm$0.63\%) & 82.27(99.49$\pm$0.23\%) & 87.86(99.58$\pm$0.05\%) \\
		UMLDA & 127.96(99.13$\pm$0.05\%) & 142.96(99.40$\pm$0.04\%) & 166.16(99.65$\pm$0.38\%) \\
		CMDA & -- & -- & -- \\
		DATER & -- & -- & -- \\
		CF-LDA & 0.28 (95.45$\pm$0.80\%) & 0.46(96.58$\pm$0.98\%) & 0.52 (96.11$\pm$1.27\%) \\
		2D-LDA & 6.16(95.73$\pm$2.39\%) &7.04(96.27$\pm$1.43\%)  &8.07(97.46$\pm$1.23\%)  \\
		2D-QLDA & -- & -- & -- \\
		\bottomrule
	\end{tabular}
	\vspace{0.5em}
	\begin{tablenotes}
		\footnotesize
		\item[1] Here, ``$-$'' refers to CPU computation time exceeding 1800 seconds.
	\end{tablenotes}
\end{table}
\begin{table}[htbp]
	\centering
	\caption{Example \ref{ex4}: Experimental outcomes obtained over the {\tt AR2600} dataset, whose size is $120 \times 2600\times 165$. The reducing dimension $d = 30$, while the count of chosen features is $k=40$.} \label{Table4.2}
	\begin{tabular}{lccc}
		\toprule
		\multirow{2}{*}{\makecell[c]{Algorithms\\(Methods)}} & \multicolumn{3}{c}{CPU Time / s ( Recognition rate $\pm$ Std-Dev\% )} \\
		\cmidrule(lr){2-4}
		& $m = 60\%N$ & $m = 70\%N$ & $m = 80\%N$ \\
		\midrule
 \bf{Algorithm \ref{Alg5}} & 6.49(56.40$\pm$2.54\%)  & 7.81(60.31$\pm$1.13\%)  & 8.40(64.35$\pm$4.25\%)  \\
		STDA & 8.08(54.27$\pm$3.47\%) &8.39(60.38$\pm$5.17\%)  & 8.49(62.07$\pm$3.79\%) \\
		SLDA & -- & -- & --\\\hline
		\bf{Algorithm \ref{Alg4}} & 1.56(48.17$\pm$3.59\%)  & 1.80(56.47$\pm$6.39\%)  & 1.97(57.88$\pm$4.61\%)  \\
		t-MLDE & 114.05(51.83$\pm$1.47\%) & 128.59(56.16$\pm$2.43\%) & 148.82(59.62$\pm$1.99\%) \\
		UMLDA & 411.35(95.40$\pm$0.66\%) & 475.58(96.62$\pm$1.17\%) & 533.24(96.82$\pm$0.55\%) \\
		CMDA & -- & -- & -- \\
		DATER & -- & -- & -- \\
		CF-LDA & 0.75 (82.50$\pm$0.39\%) & 0.84(84.46$\pm$1.12\%) & 0.93 (85.58$\pm$0.71\%) \\
		2D-LDA & 13.19(87.00$\pm$0.83\%) &15.19(88.78$\pm$0.47\%)  &17.07(90.38$\pm$1.29\%)  \\
		2D-QLDA & -- & -- & -- \\
		\bottomrule
	\end{tabular}
	\vspace{0.5em}
	\begin{tablenotes}
		\footnotesize
		\item[1] Here, ``$-$'' refers to CPU computation time exceeding 1800 seconds.
	\end{tablenotes}
\end{table}
\end{example}

\section{Conclusion}\label{sec6}
The tensor trace-ratio problem occupies a vital position within dimensionality reduction for high-dimensional large-sample tensor datasets. Its objective lies in mapping tensor samples from high-dimensional manifolds onto low-dimensional subspaces, whereby core information contained in raw inputs is retained to the greatest feasible extent. Some existing multilinear discriminant analysis methods proposed recently addressed the problem via an inner-outer iterative framework. However, when dealing with high-dimensional large-sample tensor datasets, the costs of these inner-outer iterative methods are extremely high, and the algorithms can even be infeasible.

Within the present study, our main investigation targets the tensor trace-ratio problem arising from multilinear discriminant analysis. For high-dimensional large-sample tensor data, a closed-form solution is derived assuming $\mathrm{bcirc}(\mathcal{X}^T)$ satisfies the full-row-rank condition. On this basis, a non-iterative algorithm is developed to address the tensor trace-ratio problem. Using the derived closed-form expression, the original tensor trace-ratio task gets converted to a tensor least-squares problem, for which the randomized t-SVD algorithm is adopted to calculate solutions. Corresponding theoretical analyses are carried out to validate both feasibility and reliability of the newly developed scheme. Furthermore, building upon the closed-form solution, this work develops a sparse multilinear discriminant analysis model (SMDA) regularized by the tensor $\ell_{2,0}$-norm, which can simultaneously achieve feature selection and dimensionality reduction, thereby mitigating the curse of dimensionality and information redundancy.
Numerical simulations reveal the computational characteristics of both developed methods.
The proposed methods demonstrate significant advantages in computational efficiency and scalability for large-scale problems, while achieving competitive classification accuracy on some datasets.

Several promising research directions merit deeper exploration. For example, how to apply the proposed strategy to other tensor dimensionality reduction approaches?
Furthermore, we mention that the closed-form solution presented in Section \ref{sec3} relies on the condition that ${\rm bcirc}\left(\mathcal{X}^T\right)$ is of full row rank. However, we find that both Algorithm \ref{Alg4} and Algorithm \ref{Alg5} still work very well even when this condition is unsatisfied. How to interpret it theoretically is interesting which constitutes an important direction for subsequent research.

\subsection*{Declarations}

\bmhead{Conflict of Interest}
The authors declared that they have no conflict of interest.

\bmhead{Acknowledgements}
The authors wish to extend their heartfelt gratitude to the anonymous referees and the managing editor for their constructive remarks and recommendations, which have substantially enhanced the quality of this manuscript.


\bmhead{Funding}
This work is supported by the National Natural Science Foundation of China under Grants 12671452 and 12271518, the Graduate Innovation Program of China University of Mining and Technology under Grant 2025WLKXJ140, and the Postgraduate Research \& Practice Innovation Program of Jiangsu Province under Grant KYCX25\_2854.

\bmhead{Author Contributions}
Yuhang Li: Writing the original draft; Gang Wu: Supervision, manuscript composition, revision and editorial preparation. Both authors approved the final manuscript.

\bmhead{Data availability}
Not applicable.

\end{document}